%  TeX
%  I add a %v  at the places where I modify things,
% to make them easy to find
% and  %vv  is a mark for larger modifications

%%% retractions? produce them on stripes, but compose again?
% Reflections too?

%%% connexit\'{e} des ARRF? Cutting $g$ in small
% bi-Lipschitz pieces (with the other proof?)
% But with a bi-H\"older topology, or in a special case?
% Teichmuller may be far, because of wrong topology
% and the other formula bay be better.

%%% connexe+ARRF donne juste 2 composantes?

%%% Curious that this is more complicated! (P crosses S)
% find another mapping, or a reason why the image is small??

\magnification = 1200
\input amssym.def
\input amssym.tex
\input epsf.tex

\def \ms {\medskip}
\def \msi {\medskip\noindent}
\def \bsi {\bigskip\noindent}

\def \dsp {\displaystyle}
\def \sm { \setminus}
\def \d {\partial}
\def \wt {\widetilde}

\def \qed {\hfill $\square$}
\def \R {\Bbb R}

\def \N {\Bbb N}

\overfullrule=0pt

\def \dist{\mathop{\rm dist}\nolimits}

\def \Max{\mathop{\rm Max}\nolimits}

\centerline{REIFENBERG PARAMETERIZATIONS
FOR SETS WITH HOLES}

\vskip 0.5cm
\centerline{ Guy David and 
Tatiana Toro\footnote * {The second author was partially supported by 
DMS grants 0600915 and 0856687. Part of the work was carried out
while the authors were visiting IPAM.}}
\vskip 1cm

\noindent
{\bf R\'{e}sum\'{e}.} 
On g\'{e}n\'{e}ralise la d\'{e}monstration du th\'{e}or\`{e}me du
disque topologique de Reifenberg pour inclure le cas d'ensembles
ayant des trous, et on donne des conditions suffisantes sur l'ensemble
$E$ pour l'existence de param\'{e}trage de $E$ par un plan affine 
ou une vari\'{e}t\'{e} de dimension $d$. L'une de ces conditions
porte sur la sommabilit\'{e} des carr\'{e}s des nombres de 
P. Jones $\beta_1(x,r)$, et s'applique en particulier 
aux ensembles localement Ahlfors-r\'{e}guliers et \`{a} l'existence %v
de tr\`{e}s grand morceaux d'images bi-Lipschitziennes de $\R^d$. %v

\bigskip \noindent
{\bf Abstract.}  We extend the proof of Reifenberg's Topological Disk
Theorem to allow the case of sets with holes, and give sufficient
conditions on a set $E$ for the existence of a bi-Lipschitz %v
parameterization of $E$ by a $d$-dimensional plane or smooth 
manifold. Such a condition is expressed in terms of square summability 
for the P. Jones numbers $\beta_1(x,r)$. In particular, it applies in the
locally Ahlfors-regular case to provide very big pieces of 
bi-Lipschitz images of $\R^d$.


\bsi
{\bf 1. Introduction}
\ms

In this paper we take the usual proof of Reifenberg's well-known 
topological disk theorem, and make it work in several  
different contexts. Let us give a local statement for 
Reifenberg's theorem before we discuss it further.

\ms\proclaim Theorem 1.1 [R1]. %% acta 60
For all choices of integers
$0 < d < n$, and $0 < \tau < 10^{-1}$, we can find $\varepsilon > 0$
such that the following holds.
Let $E \i \R^n$ be a closed set that contains the origin,
and suppose that for $x\in E \cap B(0,10)$ and $0 \leq r \leq 10$,
we can find a $d$-dimensional affine subspace $P(x,r)$ of $\R^n$ that
contains $x$ and such that 
$$\eqalign{
&\dist(y,P(x,r)) \leq \varepsilon r 
\hbox{ for $y\in E\cap B(x,r)$ and }
\cr& \hskip 3cm
\dist(y,E) \leq \varepsilon r 
\hbox{ for $y\in P(x,r)\cap B(x,r)$.}
}\leqno (1.2)
$$
Then there is a bijective mapping $g : \R^n \to \R^n$
such that
$$
|g(x)-x| \leq \tau \ \hbox{ for } x\in \R^n,
\leqno (1.3)
$$
$$
{1 \over 4} |x-y|^{1+\tau} \leq |g(x)-g(y)| \leq
3 |x-y|^{1-\tau}
\leqno (1.4)
$$
for $x,y\in \R^n$ such that $|x-y| \leq 1$ and, if we set 
$P = P(0,10)$,
$$
E \cap B(0,1) = g(P) \cap B(0,1).
\leqno (1.5)
$$

\ms
Thus we can get $g$ to be bi-H\"older with any exponent
smaller that $1$, if $\varepsilon$ is assumed to be 
accordingly small. Of course the constant $10$ is far from 
optimal here, and usually we do not need to know $g$ far away %v
from the origin.
Also, we could replace $1/4$ and $3$ in (1.4) with constants that
are arbitrarily close to $1$ (see for instance [DDT]), %% Proppen
but we shall make no serious attempts in this paper to prove 
(1.4) with such constants.
Often the existence of a bi-H\"older parameterization defined
on $P$ is enough, but in some cases it is good to know that
it comes from a globally defined bi-H\"older mapping.

Recall that we cannot hope to get a bi-Lipschitz mapping
$g$ in general, because very flat snowflake curves
in $\R^2$ can satisfy (1.2) with arbitrarily small values
of $\varepsilon$, but do not have finite length (and 
even have Hausdorff dimensions larger than 1). We cannot hope
to always have a quasisymmetric parameterization either,
this time because the product in $\R^3$ of a snowflake
(in $\R^2$) and a line admits no a quasisymmetric
parameterization [V].  %% Vaisala

We shall give a more global statement later (Theorem 12.3),
where $E$ is assumed to be close to some smooth
subvariety $\Sigma_0$, and we get that $E =g(\Sigma_0)$
on a more general set $U$.
But the two statements are very similar.

\ms
First we want to extend Theorem 1.1 to situations where we 
only assume that for $x\in E$ and $0 < r \leq 10$, 
we can find $P(x,r)$ such that
$$
\dist(y,P(x,r)) \leq \varepsilon r 
\hbox{ for } y\in E\cap B(x,r).
\leqno (1.6)
$$
Of course in this case we can only hope to get that
$E \cap B(0,1) \i g(P) \cap B(0,1)$ instead of (1.5),
but even so we shall give in Counterexample 12.28
an example that shows that some additional compatibility
condition between the $P(x,r)$ is needed. Let us state this
in terms of the following normalized local Hausdorff distances. Set
$$
d_{x,r}(E,F) = {1 \over r} \Max\Big\{ 
\sup_{y\in E\cap B(x,r)} \dist(y,F) \, ; \, 
\sup_{y\in F\cap B(x,r)} \dist(y,E) \Big\}
\leqno (1.7)
$$
for $x \in \R^n$ and $r > 0$, and when $E, F \i \R^n$ both
meet $B(x,r)$. [We shall not need the other case.] %v
We shall assume that
$$
d_{x,10^{-k}}(P(x,10^{-k}),P(x,10^{-k+1})) \leq \varepsilon
\leqno (1.8)  
$$
for $x\in E$ and $k \geq 0$, and that
$$
d_{x,10^{-k+2}}(P(x,10^{-k}),P(y,10^{-k})) \leq \varepsilon
\leqno (1.9)   
$$
when $k \geq 0$ and $x,y \in E$ are such that
$|x-y| \leq 10^{-k+2}$.

Here is the (local) analogue Theorem  1.1 in this
context; see Theorem  12.18 for a more general statement. 

\ms\proclaim Theorem 1.10. 
For all choices of integers $0 < d < n$ and $0 < \tau < 10^{-1}$, 
we can find $\varepsilon > 0$ such that, if $E \i B(0,1) \i\R^n$ 
is a closed set that contains the origin, and if for $x\in E$ and 
$0 \leq r \leq 10$, we can find a $d$-dimensional affine subspace 
$P(x,r)$ that contains $x$ and such that (1.6), (1.8), 
and (1.9) hold for each $k\ge 0$, 
then there is a bijective mapping $g : \R^n \to \R^n$ such that
(1.3) and (1.4) hold, and $E \i g(P(0,10))$.
Moreover, $\Sigma = g(P(0,10))$ is Reifenberg-flat, in the sense that
for $x \in \Sigma$ and $r > 0$, there is an affine $d$-plane $Q(x,r)$ 
through $x$ such that $d_{x,r}(\Sigma,Q(x,r)) \leq C \varepsilon$.

\ms
We are also interested in estimating the distortion
of $g$, and in particular getting sufficient conditions on $E$
that allow us to get a bi-Lipschitz mapping $g$. First suppose
that we are in the context of Theorem 1.1, and express a sufficient
condition in terms of the Jones numbers
$$
\beta_\infty(x,r) = {1 \over r} \, \inf_{P} \,
\Big\{ \sup \big\{\dist(y,P) \, ; \, y\in E \cap B(x,r) \big\}\Big\},
\leqno (1.11)
$$
where the infimum is taken over all $d$-planes $P$ through $x$.
[In the most usual variants, we do not require $P$ to contain
$x$, but in the present context $x$ will always lie on $E$,
so there is only a small difference.]
Then set 
$$
J_\infty(x) = \sum_{k \geq 0} \beta_\infty(x,10^{-k})^2
\leqno (1.12)
$$
for $x\in E$. This type of function was introduced by
C. Bishop and P. Jones in [BJ2] %% Schwarzian
and used a lot by Bishop, Jones, Lerman, and others in the context of 
Lipschitz or nearly Lipschitz parameterizations, so it is not 
surprising that it shows up here too. 
See for instance [BJ], [DS1,3], [J1], [J2], [JL], [L\'{e}], [Lr2], [P1]. %%

\ms \proclaim Theorem 1.13.
Let $n$, $d$, and $E$ be as in Theorem 1.1, and in
particular assume that for $x\in E \cap B(0,10)$ and 
$0 < r \leq 10$, (1.2) holds for some $d$-plane $P(x,r)$
through $x$.
Assume in addition that $J_\infty$ is bounded on $E \cap B(0,10)$.
Then the mapping provided by Theorem~1.1 is also bi-Lipschitz: there
exists $C \geq 1$, that depends only on $n$, $d$, and a bound for 
$J_\infty$, such that 
$$
C^{-1}|x-y| \leq |g(x)-g(y)| \leq C |x-y|
\ \hbox{for } x,y \in \R^n.
\leqno (1.14)
$$

\ms
The reader will probably have a good idea of what happens
in this paper by considering sets in the plane that are obtained from a
line segment by a snowflake construction, where each
segment of the $k$th generation is replaced with four
shorter segments. We allow the angles that govern the 
construction to depend on the generation and even on the 
segment in a given generation, but demand that all these
angle be small (so that we get the Reifenberg condition (1.2)),
and even depend gently on $k$ and the segment.

With no further constraint on the angles, the limit set $E$ is 
Reifenberg-flat, and we merely have a bi-H\"older parameterization.
But we can also choose the angles so that %v
the sum of the squares of the angles used in the construction
of the ancestors of any given segment be a bounded function.
In this case, Theorem 1.13 applies. In fact, the limit curve is 
chord-arc, and the existence of a bi-Lipschitz mapping
of $\R^2$ that sends the unit circle to the curve would also 
follow from an extension result of [Tu] or [JeK]. %% 

We can also think about similar constructions in higher
dimensions, and get sufficient conditions (that are now
further from being necessary). In a way, our proofs will say that
this type of example gives a good idea of the general situation,
because we shall rely on successive approximations of $E$ by
$d$-planes, and worry about the square summability of the normalized
distance between them. The fact that the planes do not correspond
exactly to faces of intermediate objects will not matter much.

We can also use the possibly smaller numbers
$$
\beta_q(x,r) = \inf_{P} \Big\{ 
r^{-d} \int_{y \in E \cap B(x,r)} {\dist(y,P)^q \over r^q} 
\, dH^d(y) \Big\}^{1/q},
\leqno (1.15)
$$
defined for $x\in E \cap B(0,10)$, $0 < r \leq 10$, and 
$1 \leq q < +\infty$, and where this time the infimum is taken 
over all $d$-planes $P$ through $B(x,r)$. [We do not want
to force $P$ through $x$, because $x$ may exceptionally
be a little far from the best plane.] Then set
$$
J_q(x) = \sum_{k \geq 0} \beta_q(x,10^{-k})^2
\leqno (1.16)  %vv Not understood what happenned there
$$ % but I checked the next occurences
for $x\in E\cap B(0,10)$.
In the special case when $E$ is locally Ahlfors-regular, i.e.,
when there is a constant $C \geq 1$ such that
$$
C^{-1} r^d \leq H^d(E\cap B(x,r)) \leq C r^d
\ \hbox{for } x \in E \hbox{ and } 0 < r \leq 10,
\leqno (1.17)
$$
$\beta_q(x,r) \leq C\beta_p(x,r) \leq C'\beta_\infty(x,r)$ 
when $1 \leq q < p < +\infty$, by H\"older's inequality and (1.17).

Here is the analogue of Theorem 1.13 in this context.

\ms \proclaim Theorem 1.18.
Let $n$, $d$, and $E$ be as in Theorem 1.1, and in particular assume 
that for $x\in E \cap B(0,10)$ and  $0 < r \leq 10$, (1.2) holds 
for some $d$-plane $P(x,r)$ through $x$.
Assume in addition that $J_1$ is bounded on $E \cap B(0,10)$.
Then the mapping provided by Theorem~1.1 is also bi-Lipschitz: there
exists $C \geq 1$, that depends only on $n$, $d$, and a bound for $J_1$,
such that (1.14) holds.

\ms
Note that the Ahlfors-regularity property (1.17)
is not used in Theorem 1.18 (although the lower bound
will be proved in Lemma~13.6), but as long as we do not know 
that (1.17) holds, we cannot be sure that 
$\beta_1(x,r) \leq C \beta_\infty(x,r)$ and $J_1\leq C J_{\infty}$. 
[On the other hand, a posteriori, (1.17) holds if $E$ is the 
bi-Lipschitz image of a $d$-plane.]

The boundedness of $J_1$ is not necessary for $E$ to be
the bi-Lipschitz image of a $d$-plane, but it is not too far off:
some sort of BMO condition is needed. See Remark 15.38.

See Corollaries 12.44 and 13.4 for more general analogues 
of Theorems 1.13 and 1.18. 

In codimension $1$, we shall also give in Corollary 13.46 
a sufficient condition for the boundedness of $J_1$ 
(if $E$ is locally Reifenberg-flat, as in Theorem 1.1 or 12.3),
and hence for the existence of a bi-Lipschitz parameterization
as above. This sufficient condition is expressed in terms 
of the unit normal to $E$. It is reminiscent of one of the
equivalent definitions of the Chord-Arc Surfaces with Small 
Constants introduced and studied by Semmes [Se1,2,3]. %%

\ms
There is also a sufficient condition for $g$ in Theorem 1.10
to be bi-Lipschitz, which is expressed in terms of the squares
of the distances implicit in (1.8) and (1.9). That is, set
$$\eqalign{
\varepsilon_k(x) &= d_{x,10^{-k}}(P(x,10^{-k}),P(x,10^{-k+1})) 
\cr& \hskip 2cm
+ \sup_{y\in E\cap B(x,10^{-k-2})} d_{x,10^{-k+1}}(P(x,10^{-k}),P(y,10^{-k})) 
}\leqno (1.19)
$$
for $x\in E$ and $k \geq 0$, and then
$$
J(x) = \sum_{k \geq 0} \varepsilon_k(x)^2.
\leqno (1.20)
$$

\ms \proclaim Theorem 1.21.
Let $n$, $d$, and $E$ be as in Theorem 1.10, and assume in 
addition that $J$ is bounded on $E$. 
Then the mapping provided by Theorem~1.10 is %v
also bi-Lipschitz, i.e., (1.14) holds with a constant that 
depends only on $n$, $d$, and a bound for $J$.

\ms
See Corollary 12.33 for a more general statement.

\ms
The various statements above, and their generalized
counterparts, are all derived with the same algorithm, and
slightly different parameters. The algorithm also allows a 
(single in our case) stopping time, which for instance %v
allows us, when the functions $J_1$ or $J_\infty$ above  
are unbounded, to get a bi-Lipschitz mapping $g$ such that
$g(P(0,10))$ contains the points of $x\in E$ such that 
$J_1(x)$, or $J_\infty(x)$, is less than a constant.
[See Remark 14.13.]
This works best if we have a good control on 
$J_1$ or $J_\infty$, as in the following.

In the case of Reifenberg-flat sets (i.e., as in Theorem 1.1) 
that are also locally Ahlfors-regular (as in (1.17)), the mapping 
$g$ provides a bi-Lipschitz image 
of a $d$-plane that covers most of $E \cap B(0,1)$. 
This is not too surprising, because we could expect
(and indeed show in Theorem 15.4 ), with such a strong 
assumption as local Reifenberg flatness, that $E$ is locally 
uniformly rectifiable, and even contains big pieces of
Lipschitz graphs. The fact that $E \cap B(0,1)$ is 
almost covered by a bi-Lipschitz image of $\R^d$ (but maybe
in a larger $\R^m$) then follows from [DS1]. %% asterisques&perils

Anyway, we get a suitable control of $J_1$ from
the local uniform rectifiability of $E$, and this allows
us to find a bi-Lipschitz mapping $g$ such that 
$g(P(0,2))$ contains most of $E\cap B(0,1)$.
See Theorem 14.1.

What is new here is the fact that we do not need to enlarge
the ambient space to construct bi-Lipschitz mappings
defined on $\R^d$, and that the bi-Lipschitz mappings
in question even have a bi-Lipschitz extension to $\R^n$.

Note that for Chord-Arc Surfaces with Small Constants (CASSC),
this was known, and $E$ even contains very big pieces of
Lipschitz graphs [Se1,2,3]. %%%  or refer to the right Semmes?
Our case is somewhere in the middle (the CASSC are known to be
locally Reifenberg-flat).

When $d = n-1$, and when we still assume (1.2) and (1.17),
our construction will also give (disjoint) domains 
$\Omega_1$ and $\Omega_2$, which are bi-Lipschitz images of half 
spaces, do not meet $E \cap B(0,1)$, but are such that 
$\partial \Omega_1 \cap \partial \Omega_1$ is contained
in $E$ and contains most of $E \cap B(0,1)$. 
The difference with a similar result of [DJ] %% cond B for harmonic
that applies in the more general context of Condition B domains,
is that we get very big pieces and a simultaneous approximation
from both sides of the set. See Propositions 14.16 and 15.45.
Also see Section 15 for other comments about the properties
of locally Ahlfors-regular Reifenberg-flat sets and uniformly
rectifiable sets. %v

When $E$ is locally Reifenberg-flat and Ahlfors-regular,
the mapping $g$ that we get from Theorems 1.1 and 12.3
is not bi-Lipschitz in general, but it is controlled by
$J_1$, which is not so large (it satisfies BMO-type
estimates, by local uniform rectifiability and the
so-called geometric lemma from [DS1]). %% Astersique 
So our mapping $g$ may potentially be useful for some
problems. This is an interesting question which %v
we do not pursue further at this moment. 
That is, we shall not try to follow up on the distortion 
estimates that we could deduce from Sections~7-11 and 
the geometric lemma.

Recall that for uniformly rectifiable sets, we have
reasonable parameterizations that are not Lipschitz,
but controlled by an $A_1$ weight (Condition C7 on
page 14 of [DS1]). These parameterizations %% Astersique 
are neither injective nor surjective, though. On the opposite 
end, Semmes [Se2] gives good parameterizations of the CASSC, %% 
with an $L^p$ control on the derivative, with $p <+\infty$
as large as we want. The authors did not check whether
his parameterization extends nicely to $\R^n$.

We return to the general case of Theorem 1.1 or 12.3.
Observe also that $g$ may be used to find a H\"older 
retraction from a neighborhood of $E$ onto $E$, or a
H\"older reflection across $E$, and again
these mappings should not be far from Lipschitz when
$E$ is locally Ahlfors-regular. 
In the slightly more general context of [DDT] 
(where $E$ is uniformly close to minimal cones, and not
just planes) for two-dimensional sets in $\R^3$, 
retractions onto $E$ may eventually be used to prove existence
results for two-dimensional almost minimal sets. 

\ms
There is essentially one construction in this paper,
and various estimates on derivatives and distances
between planes.
The idea of the main construction comes from
Reifenberg [R1], revised by Morrey [Mo] %% %%
and many others (e.g [Si], [To]). %%
The systematic use of the functions $\beta(x,r)$ and $J(x)$ was 
introduced by P. Jones and C. Bishop (starting from [J1] and [BJ2]) %%
and used successfully by many others ([BJ], [DS], [J1, [JL], %%
[Lr], [L\'{e}], [P1], [P2], [Se1], [Se2], [Sc], [To]), %%
in the context of parameterizing sets
in a Lipschitz or almost Lipschitz way.
The present argument is a mixture of both techniques,
but we do very little in terms of stopping time regions.

One of the reasons why we think the Reifenberg construction
is very powerful is that it is a top down algorithm
which allows us to move the points little by little 
(so that they land in $E$ at the end); this is a little
more flexible than the standard stopping time arguments
that tend to project points on a single Lipschitz graph
and then stop. Here we work more, but we stop less
(and hence need to glue less). Probably there is a way
to incorporate the present paper in a stopping time argument,
as in the work of [JL]. %%

A minor difference between this paper and the previous ones is
that we decided that the main thing that governs the construction
is the choice of approximating $d$-planes (at all scales and %v
locations), rather than the set itself. This is why we still
can construct $g$ when $E$ has big holes, provided that we
can choose planes in a coherent way. Trying to extend $E$
instead seems more complicated.

Before we started to write this up, we thought
this would be a good opportunity to write down a simple 
proof of Reifenberg's result. This hopefully worked up to
Section 5, but maybe not later. This is not bad, %v
because Section 4 gives a very good idea of the the algorithm, 
which is the most important part of the proof. The reader will 
probably think that after this, things get a little technical, 
but we could not help it. We tried to cut the proof in somewhat 
independent pieces too.

The plan for the rest of this paper is as follows.
Section 2 contains our basic assumptions on a model manifold 
$\Sigma_0$ (such as a plane), collections of 
$10^{-k}$-nets $\{ x_{j,k} \}$, $j \in J_k$,
and families $\{ P_{j,k} \}$, $j \in J_k$, of affine
$d$-planes through the $x_{j,k}$. These provide the initial
setting that will govern the construction. We also give 
there the two main technical statements (Theorem 2.15
and its complement Theorem 2.23 for the bi-Lipschitz estimates).

We construct a partition of unity $\{ \theta_{j,k} \}$
in Section 3, and use it in Section 4 to define our main
mapping $f$, whose goal is to send a small neighborhood
of $\Sigma_0$ to our final set (typically, the set $E$ in 
Theorem 1.1). We will obtain $f$ as the limit of composed 
functions $f_k =\sigma_{k-1} \circ \cdots \circ \sigma_{0}$,
where each $\sigma_k$ moves points near $E$ (or rather, near %v
the planes $P_{j,k}$ of the $k$-th generation) at the scale $10^{-k}$. 

The proof of Theorems 2.15 and 2.23 will keep us busy for
Sections 5-12. In Section~5 we show by induction that near
each $x_{j,k}$, $\Sigma_k = f_k(\Sigma_0)$ coincides with a 
small Lipschitz graph over $P_{j,k}$ (see Proposition 5.4).
We use this in Section 6 to show that the limit
$\Sigma = f(\Sigma_0)$ is Reifenberg-flat (recall that when
we start with a set $E$ with holes, $\Sigma$ will
be larger than $E$). In Section 7, we estimate 
the differential of $\sigma_k$ in terms of distances between
the $P_{j,k}$, and we use this in Section 8 to prove the
desired bi-H\"older or bi-Lipschitz estimates on the
restriction of $f$ to $\Sigma_0$. At this point we have
a good parameterization of $\Sigma$, which we still need
to extend to $\R^n$.

First we give $C^2$ estimates on the intermediate
surfaces $\Sigma_k$, which we use to construct isometries
$R_k(z)$, $z\in \Sigma_0$, that map the tangent plane
to $\Sigma_0$ at $z$ to the tangent plane to $\Sigma_k$
at $f_k(z)$. These play the same role as the continuous
choice of orthonormal basis for the tangent plane to $\Sigma_k$
that was used in [Mo] or [To]. %% %%
We finally define the full $g$ in Section 10, and prove the 
desired bi-H\"older or bi-Lipschitz estimates in Section 11.

Sections 12-15 contain various applications of the previous
construction. In Section~12 we give the main generalization
of Theorem 1.1 (namely, Theorem 12.3), its variant with holes
(Theorem 12.18), and their bi-Lipschitz variants
(Corollaries 12.33 and 12.44). These are mostly
expressed in terms of numbers $\beta_\infty(x,r)$.

Section 13 contains a variant of Corollary 12.44
(a bi-Lipschitz statement) expressed in terms of
numbers $\beta_1(x,r)$, and a sufficient condition
for sets of codimension one to be contained in a bi-Lipschitz
image of $\R^d$ or $\Sigma_0$, expressed in terms of the 
(continuous) unit normal to $E$.

In Section 14 we show that if in addition to the assumptions
of Theorem 1.1 or Theorem 12.3, $E$ is locally Ahlfors-regular,
then we can find a bi-Lipschitz mapping $g$ as above, such that
$g(\Sigma_0)$ contains most of $E$. See Theorem 14.1,
and Proposition 14.16 for the description of
saw-tooth domains in $\R^n \sm \Sigma_0$
whose images by $g$ do not meet $E$, but 
have a big piece of $E$ in their boundary.
The proof of Theorem 14.1 is completed in Section~15,
where we also discuss the uniform rectifiability properties
of the locally Ahlfors-regular Reifenberg-flat sets.

The authors whish to thank Raanan Schul for interesting conversations
on prameterizations.
%%% others? %vv

\bsi
{\bf 2. Coherent families of balls and planes}
\ms

Let us first describe the simplest situation where we can 
create a $d$-dimensional Reifenberg-flat set $\Sigma$ and a 
parameterization of $\Sigma$ by a $d$-plane or a smooth surface. 
In the case of the standard Reifenberg theorem 
(Theorem 1.1 for instance), $\Sigma$ will coincide
with $E$ on $B(0,1)$. In the situation of Theorem 1.10, 
$E$ will be contained in $\Sigma$.

At the end of this section, we shall give our main technical
statements, which will be proved in Sections 3-11, and
made cleaner or applied only in later sections.

First set $r_k = 10^{-k}$ for $k\in \N$, and choose 
a collection $\{ x_{j,k} \}$, $j\in J_k \,$, of points in $\R^n$,
so that
$$
|x_{i,k}-x_{j,k}| \geq r_k
\ \hbox{ for } i,j\in J_k, i\neq j.
\leqno (2.1)
$$
In our applications, the points $x_{j,k}$ will lie on a given set $E$,
but we do not need to know this.
Also set $B_{j,k} = B(x_{j,k},r_k)$ and, for $\lambda > 1$,
$$
V_k^\lambda = \bigcup_{j\in J_k} \lambda B_{j,k} = 
\bigcup_{j\in J_k} B(x_{j,k}, \lambda r_k)
\leqno (2.2)
$$
(we usually denote by $\lambda B$ the ball with the same center 
as $B$ and $\lambda$ times the radius). We shall
assume that for $k\geq 1$ and $j \in J_k$, 
$$
x_{j,k} \in V_{k-1}^2.
\leqno (2.3)
$$

The most standard way to produce the $x_{j,k}$ is the following.
We start from a set $E_0 \i \R^n$ (typically, a subset of 
$E$ from the previous section), and then choose a nonincreasing
sequence $\{ E_k \}$ of subsets of $E_0$ (to allow for a stopping
time argument). For each $k \geq 0$, we let the $x_{j,k}$, $j\in J_k$
be a maximal collection of points of $E_k$ such that (2.1) holds.
Then of course $E_k \i \cup_{j \in J_k} \overline B_{j,k}$, and
(2.3) follows because $x_{j,k} \in E_k \i E_{k-1}$.

Again we shall not need to know that the $x_{j,k}$ were produced that
way, but (2.3) is nonetheless a stopping time coherence condition, or
a way of asking that we do not resume the construction in places where 
we stopped it at a larger scale.

\ms
We shall also assume that the initial points $x_{j,0}$ are close to 
some smooth $d$-dimensional surface $\Sigma_0$. For the results mentioned
in Section 1, $\Sigma_0$ will be a plane, but it will not disturb us 
too much to allow slightly more general surfaces, very flat at the 
unit scale, but with a potentially complicated behavior at larger
scales. We assume that for each $x\in \Sigma_0$, there is a $d$-plane $P_x$ 
through $x$ and a $C^2$ function $F_x : P_x \to P_x^{\perp}$
(the $(n-d)$-dimensional vector space of $\R^n$ which is orthogonal
to $P_x$) such that
$$
|F_x(y)|+|DF_x(y)|+|D^2F_x(y)| \leq \varepsilon 
\ \hbox{ for } y \in  P_x,  
\leqno (2.4)    
$$
where we denote by $DF_x(y)$ the differential of $F_x$ at $y$
and by $D^2F_x$ the differential of $DF_x$; the choice
of norm for $|DF_x(y)|$ would not matter much, but let us take
the operator norm as acting on vectors, and
$$
\Sigma_0 \cap B(x,200) = \Gamma_{F_x} \cap B(x,200),
\leqno (2.5) 
$$ 
where we denote by $\Gamma_{F_x} = \{ y + F_x(y) \, ; \, y\in P_x \}$
the graph of $F_x$ over $P_x$.

We may also assume that $\Sigma_0$ is smooth up to some order
$m_0 > 2$, and precisely that for $2 < m \leq m_0$,
there exists $M_m \geq 0$ such that, in the description above,
$$
|D^m F_x| \leq M_m 
\ \hbox{ on } P_x,
\leqno (2.6)
$$
because this assumption essentially costs us nothing (see
Remark 2.13); then the intermediate mappings $f_k$ that we construct
will be smooth of the same order, and so will be the mapping $g$ away
from $\Sigma_0$. But the precise estimates will mostly be left to the 
reader.

The relation with our initial net $\{ x_{j,0} \}$ is that %v
we assume that
$$
\dist(x_{j,0},\Sigma_0) \leq \varepsilon
\hbox{ for } j\in J_0.
\leqno (2.7)  
$$            

\ms
The last part of our structure is a coherent family of $d$-planes.
For each $k \geq 0$ and $j \in J_k$, we assume that we are
given a $d$-plane $P_{j,k}$ through $x_{j,k}$, and we shall require
some compatibility conditions to hold. Let us use the 
normalized local Hausdorff distance $d_{x,r}(E,F)$ defined in (1.7);
we demand that 
$$
d_{x_{j,k},100r_k}(P_{i,k},P_{j,k}) \leq \varepsilon
\ \hbox{ for $k \geq 0$ and $i,j \in J_k$ such that }
|x_{i,k}-x_{j,k}| \leq 100 r_k \, ,
\leqno (2.8)
$$
that
$$
d_{x_{i,0},100}(P_{i,0},P_{x}) \leq \varepsilon
\ \hbox{ for $i \in J_0$ and $x\in \Sigma_0$ such that }
|x_{i,0}-x| \leq 2,  
\leqno (2.9)
$$
and, for $k \geq 0$, that
$$
d_{x_{i,k},20r_k}(P_{i,k},P_{j,k+1}) \leq \varepsilon
\ \hbox{ for $i \in J_k$ and $j\in J_{k+1}$ such that }
|x_{i,k}-x_{j,k+1}| \leq 2 r_k.
\leqno (2.10)
$$

\ms\proclaim Definition 2.11.
A \underbar{coherent collection of balls and planes} (in short
a CCBP) is a triple $(\Sigma_0, \{ B_{j,k}\}, \{ P_{j,k}\})$,
with the properties that we just described 
(see (2.1), (2.3), (2.7), (2.8), (2.9) and (2.10)) . We shall always
assume that $\varepsilon > 0$ is small enough, depending on
$n$ and $d$.  

\msi {\bf Remark 2.12.}
In the standard Reifenberg case, $\{ x_{j,k} \}$ will be 
an $r_k$-dense collection chosen in a Reifenberg-flat set $E$,
and $P_{j,k}$ will be chosen such that 
$d_{x_{j,k},110r_k}(E,P_{j,k}) \leq \varepsilon$.
Then (2.8)-(2.10) (with the constant $C\varepsilon$)
will follow from elementary geometry. This will be checked %v
near (12.14), for the proof of Theorem 12.3.

Even when we want to study a compact set, we may
find it more convenient to use an unbounded set $\Sigma_0$,
such as a plane. We will just need to choose points $x_{j,k}$ that 
lie in a compact set, and our construction will simply leave the 
faraway part of $\Sigma_0$ alone. We find it amusing to allow
sets $\Sigma_0$ that are different from planes, and even that may
not be orientable, and get a global statement anyway. But the
construction is essentially local: our mapping $g$ will coincide
with the identity away from $\Sigma_0$.

\msi {\bf Remark 2.13.}
It would seem more natural to assume only that
$$
|F_x(y)| \leq \varepsilon
\ \hbox{ for } y \in  P_x
\leqno  (2.14)
$$
instead of (2.4), but since the only relation between $\Sigma_0$ 
and the points $x_{j,k}$ is through the proximity relation (2.7), 
it would be easy to check that if we are given $\Sigma_0$ with
(2.7) and the weaker property (2.14), we can replace it with a smoother
$\Sigma'_0$ that satisfies (2.4) and (2.7) with the constant 
$C\varepsilon$. So we decided to require (2.4) directly, and avoid
the smoothing argument. For the same reason, we will not really
lose generality by assuming (2.6).

Let us now state our two main technical results, to be proved in 
Sections 3-11 and improved or applied later.

\ms\proclaim Theorem 2.15.
Let $(\Sigma_0, \{ B_{j,k}\}, \{ P_{j,k}\})$ be a CCBP (as in
Definition 2.11), and assume that $\varepsilon$ is small enough,
depending on $n$ and $d$. Then there is a bijection 
$g : \R^n \to \R^n$, with the following properties:
$$
g(z) = z \ \hbox{ when } \dist(z,\Sigma_0) \geq 2,
\leqno (2.16)
$$
$$
|g(z)-z| \leq C \varepsilon
\ \hbox{ for } z\in \R^n,
\leqno (2.17)
$$
$$
{1 \over 4} |z'-z|^{1+C\varepsilon} \leq |g(z')-g(z)| 
\leq 3 |z'-z|^{1-C\varepsilon}  
\leqno (2.18)
$$
for $z,z'\in \R^n$ such that $|z'-z| \leq 1$,
and $\Sigma = g(\Sigma_0)$ is a $C\varepsilon$-Reifenberg
flat set that contains the accumulation set
$$\eqalign{
E_\infty &= \big\{ x\in \R^n \, ; \, 
x \hbox{ can be written as }
x = \lim_{m\to +\infty} x_{j(m),k(m)}, 
\hbox{ with } k(m) \in \Bbb N 
\cr&\hskip 2.5cm\hbox{ and } 
j(m) \in J_{k(m)} \hbox{ for } m \geq 0,  \hbox{ and }
\lim_{m \to  +\infty} k(m) = +\infty \big\}.
}\leqno (2.19)
$$

\ms
Here and below, $C$ is used to denote constants that may depend
on $n$ and $d$, but not on $\varepsilon$ (and even less on $z$ or 
$z'$); the actual value of $C$ may vary a lot from one line to the next.

By $C\varepsilon$-Reifenberg flat set, we mean that
for $x\in \Sigma$ and $0 < t \leq 1$, there is an affine $d$-plane
$P(x,t)$ through $x$ such that
$$
d_{x,t}(\Sigma,P(x,t)) \leq C \varepsilon.
\leqno (2.20)
$$

More information on $g$ and $\Sigma$ will be obtained during the
proof, but these are the main properties. If we want $g$ to be 
bi-Lipschitz, we need additional information on the speed at which
the $P_{j,k}$ may change. Set
$$\eqalign{
\varepsilon'_k(y) = \sup\big\{ &
d_{x_{i,l},100r_{l}}(P_{j,k},P_{i,l}) \, ; \,
j\in J_k, \, l \in \{ k-1, k \}, 
\cr& \hskip 3.5cm
i\in J_{l}, \hbox{ and } y \in 10 B_{j,k} \cap 11 B_{i,l} 
\big\}  
}\leqno (2.21) 
$$
for $k\geq 1$ and $y \in  V_k^{10}$, and
$$
\varepsilon'_k(y) = 0
\ \hbox{ when } y\in \R^n \sm V_k^{10},
\leqno (2.22) 
$$
i.e., when there are no pairs $(j,k)$ as above.

\ms\proclaim Theorem 2.23.
Still assume that $(\Sigma_0, \{ B_{j,k}\}, \{ P_{j,k}\})$ is a CCBP,
with $\varepsilon$ is small enough (depending on $n$ and $d$).
Assume in addition that for some $M < +\infty$
$$
\sum_{k \geq 0} \varepsilon'_k(g_k(z))^2 \leq M
\ \hbox{ for all } z \in \Sigma_0.
\leqno (2.24)
$$
Then $g : \R^n \to \R^n$ is bi-Lipschitz
: there is a
constant $C(n,d,M) \geq 1$ such that
$$
C(n,d,M)^{-1} |z'-z| \leq |g(z')-g(z)| \leq C(n,d,M) |z'-z|
\ \hbox{ for } z,z'\in \R^n.
\leqno (2.25)
$$

\ms
Note that $\varepsilon$  is not required to depend on $M$.
Condition (2.24) looks quite cumbersome as it depends on $g$.
In Sections 12 and 13 we provide sufficient conditions that imply (2.24) but do not
depend on $g$.

Our plan is to prove these two theorems in Sections 3-11,
and then only return to the statements and applications.

\bsi
{\bf 3. A partition of unity}
\ms  

From now on, we are given a CCBP
$(\Sigma_0, \{ B_{j,k}\}, \{ P_{j,k}\})$, as in Definition 2.11, 
and we start to construct a Reifenberg parameterization 
of some set $\Sigma$. This short section is  devoted to the 
construction of partitions of unity adapted to the $\{ B_{j,k}\}$.

Select a basic $C^{\infty}$ function $\theta$, with compact
support in $B(0,10)$, such that $0 \leq \theta(y) \leq 1$ 
everywhere, and $\theta(y)=1$ on $B(0,9)$. Then set
$$
\wt \theta_{j,k}(y) = \theta(r_k^{-1}(y-x_{j,k}))
\leqno (3.1)
$$
for $j\in J_k$ and $y\in \R^n$. Observe that
$$
\sum_{j\in J_k} \wt\theta_{j,k}(y) \geq 1
\ \hbox{ for } y\in V_k^9 = \bigcup_{j\in J_k} 9B_{j,k} \, .
\leqno (3.2)
$$
Since we also want to cover $\R^n \sm V_k^9$,
we shall need additional balls and function.
Choose a maximal collection $\{ x_{l,k} \}$, $l\in L_k$,
of points of $\R^n \sm V_k^9$, such that
$$
|x_{l,k}-x_{m,k}| \geq {r_k \over 2} \ \hbox{ for $l,m\in L_k$
such that $m\neq l$},
\leqno (3.3)
$$
and then set $B_{l,k} = B(x_{l,k},{r_k \over 10})$
for $l\in L_k$. Obviously
$$
\hbox{the $9B_{l,k}$, $l \in L_k$, cover $\R^n \sm V_k^9$}
\leqno (3.4)
$$
by maximality, so that if we set
$$
\wt\theta_{l,k}(y) = \theta\Big({10 (y-x_{l,k})\over r_k}\Big)
\leqno (3.5)
$$
for $l\in L_k$ and $y\in \R^n$, we get that
$$
 \Theta(y) := \sum_{j\in J_k \cup L_k} \wt\theta_{j,k}(y) \geq 1
\ \hbox{ for } y\in \R^n,
\leqno (3.6)
$$
by (3.2) and because $\theta_{l,k}(y) = 1$ on $9B_{l,k}$. 
Of course $\Theta(y) \leq C$ because the $10B_{j,k}$,
$j\in J_k \cup L_k$, have bounded overlap. Now set
$$
\theta_{j,k}(y) = \wt\theta_{j,k}(y)/\Theta(y)
\hbox{ for } j\in J_k \cup L_k.
\leqno (3.7)
$$
Then
$$
\sum_{j\in J_k \cup L_k} \theta_{j,k}(y) = 1.
\leqno (3.8)
$$
By construction,
$$
\theta_{j,k} \hbox{ is nonnegative and compactly supported in } 
10B_{j,k}
\leqno (3.9)
$$
and
$$
|\nabla^m \theta_{j,k}(y)| \leq C_m r_k^{-m}
\leqno (3.10)
$$
for $y\in \R^n$, $j\in J_k \cup L_k$ and $m \geq 0$, because we have
similar estimates on the $\wt\theta_{j,k}$, by (3.6), and 
because the $10B_{j,k}$ have bounded overlap.

We now put all the $\theta_{l,k}$, $l\in L_k$,  
together in a single function $\psi_k$. After this, we will not
 need
to mention the $\theta_{l,k}$, $l\in L_k$, any more. Set
$$
\psi_k = \sum_{l\in L_k} \theta_{l,k}.
\leqno (3.11)
$$
Observe that
$$
\psi_k(y) = 0 \hbox{ on }
V^{8}_k = \bigcup_{j \in J_k} 8B_{j,k}
\leqno (3.12)
$$
because (3.9) says that for $l \in L_k$, $\theta_{l,k}$ is 
supported in $10B_{l,k} = B(x_{l,k},r_k)$ with an $x_{l,k}$ 
that lies out of $V^9_k$ by definition of $L_k$. Then, since
$$
\psi_k(y) + \sum_{j\in J_k} \theta_{j,k}(y) = 1
\leqno (3.13)
$$
by (3.8), we deduce from (3.12) that
$$
\sum_{j\in J_k} \theta_{j,k}(y) = 1
\ \hbox{ for } y \in V^8_k.
\leqno (3.14)
$$
Finally observe that
$$
|\nabla^m \psi_{k}(y)| \leq C_m r_k^{-m}
\leqno (3.15)
$$
for $y\in \R^n$ and $m \geq 0$, by (3.10) and because the
$10B_{l,k}$ have bounded overlap.

\bsi  
{\bf 4. Definition of a mapping $f$ on $\Sigma_0$}
\ms  

Our plan is to define first a bi-H\"older mapping $f$ on $\Sigma_0$.
In fact, we define $f$ on the whole $\R^n$, but later on, we
only care about the values of $f$ on $\Sigma_0$.
The function $f$ appears as the limit of a sequence of functions
$f_k$, where $f_k$ is defined by induction by 
$$
f_0(y)=y \ \hbox{ and } \,
f_{k+1} = \sigma_k \circ f_k
\leqno (4.1)
$$
for some $\sigma_k$ whose main goal is to push
points in the direction of the $d$-planes $P_{j,k}$
wherever they are defined. We set 
$$
\sigma_k(y) = y + \sum_{j \in J_k} \theta_{j,k}(y) \, [\pi_{j,k}(y)-y]
=\psi_k(y) y + \sum_{j \in J_k}\theta_{j,k}(y) \, \pi_{j,k}(y),
\leqno (4.2)
$$
where $\pi_{j,k}$ denotes the orthogonal projection
from $\R^n$ onto $P_{j,k}$ and by (3.13). 

We want to show (in the next few sections) that the $f_k$ tend to a 
limit $f$ which is bi-H\"older on $\Sigma_0$, and that 
$\Sigma = f(\Sigma_0)$ is a Reifenberg-flat set.
Later on, we shall extend the restriction of $f$ to $\Sigma_0$
so that it is defined and bi-H\"older on $\R^n$.
Note that
$$
|\sigma_k(y) - y| \leq 10 r_k  
\ \hbox{ for } y\in \Bbb R^n,
\leqno (4.3)
$$ 
because $\sum_{j \in J_k} \theta_{j,k}(y) \leq 1$, and 
$|\pi_{j,k}(y)-y| \leq 10r_k$ when $\theta_{j,k}(y) \neq 0$ 
(because then $y\in 10B_{j,k}$ by (3.9)). So 
$$
||f_{k+1}-f_k||_\infty \leq 10r_k,
\leqno (4.4)
$$
and the sequence $\{ f_k \}$ converges uniformly on $\R^n$ 
to a continuous function $f$. 

Clearly each $\sigma_k$ is smooth. We shall need
estimates on its derivative. First note that
$$
\sigma_k(y) = y
\ \hbox{ and } \, 
D\sigma_k (y) = I
\ \hbox{ for } y \in \R^n \sm V_{k}^{10},
\leqno (4.5)
$$
where $I$ denotes the identity map and we set
$V_{k}^{10} = \cup_{j\in J_k} 10B_{j,k}$, just by (3.9) and (4.2).
Next we record what we get when $y\in V_{k}^{10}$.

Denote by $\pi_{j,k}$ the orthogonal projection onto
$P_{j,k}$, by $D\pi_{j,k}$ its differential (which is also
the orthogonal projection onto the vector space parallel
to $P_{j,k}$), by $P_{j,k}^\perp$ the $(n-d)$-vector space 
orthogonal to $P_{j,k}$, and by $\pi_{j,k}^\perp$ the orthogonal 
projection onto $P_{j,k}^\perp$. Notice that $D\pi_{j,k}^\perp = 
\pi_{j,k}^\perp$ because $\pi_{j,k}^\perp$ is linear.

\ms \proclaim Lemma 4.6. 
For $y\in V_{k}^{10}$, we have that
$$
D\sigma_k (y) 
= \psi_k(y) I + \sum_{j \in J_k}\theta_{j,k}(y) D\pi_{j,k}
+ y D\psi_k(y) + \sum_{j \in J_k} \pi_{j,k}(y) D\theta_{j,k}(y).
\leqno (4.7)
$$
In addition, choose any $i = i(y) \in J_k$ such that 
$y\in 10B_{i,k}$, and set 
$$
L(y) 
= \psi_k(y) D\pi_{i,k}^\perp + D\pi_{i,k} + [y - \pi_{i,k}(y)] D\psi_k(y).
\leqno (4.8)
$$
Then 
$$
|D\sigma_k (y) - L(y)| \leq C \varepsilon.
\leqno (4.9)
$$

\ms
First, (4.7) comes directly from the second part of (4.2).
Then notice that when we replace the various $\pi_{j,k}$ in 
the right-hand side of (4.7) with $\pi_{i,k}$, we get the quantity

$$\eqalign{
\psi_k(y) I &+ \sum_{j \in J_k}\theta_{j,k}(y) D\pi_{i,k}
+ y D\psi_k(y) + \sum_{j \in J_k} \pi_{i,k}(y) D\theta_{j,k}(y)
\cr=& \, \psi_k(y) I + (1-\psi_k(y)) D\pi_{i,k}
+ y D\psi_k(y) - \pi_{i,k}(y) D\psi_k(y)
\cr
=& \, \psi_k(y) D\pi_{i,k}^\perp + D\pi_{i,k} + [y - \pi_{i,k}(y)] D\psi_k(y)
= L(y),
}\leqno (4.10)
$$
by (3.13).
But $y\in 10B_{j,k}$ for all $j$ such that $\theta_{j,k}(y) \neq 0$ 
or $D\theta_{j,k}(y) \neq 0$, so (2.8) says that 
$$
d_{x_{i,k},100r_k}(P_{i,k},P_{j,k}) \leq \varepsilon,
\leqno (4.11)
$$
which implies that
$$
|\pi_{j,k}(y)-\pi_{i,k}(y)| \leq 500 \varepsilon r_k
\ \hbox{ and } \, 
|D\pi_{j,k}-D\pi_{i,k}| \leq 500 \varepsilon. 
\leqno (4.12)
$$
Thus (4.7) yields
$$\eqalign{
|D\sigma_k (y) - L(y)| 
&\leq \sum_{j \in J_k}\theta_{j,k}(y) | D\pi_{j,k}- D\pi_{i,k}|
+ \sum_{j \in J_k} |\pi_{j,k}(y)- \pi_{i,k}(y)| |D\theta_{j,k}(y)|
\cr&
\leq C \varepsilon,
}\leqno (4.13)
$$
by (3.10), and Lemma 4.6 follows. \qed

\ms
The situation for $y\in V_{k}^{8}$ is much simpler,
because then $\psi_k(y)=0$ and $D\psi_k(y)=0$ by (3.12),
so (4.7) and (4.8) simplify. Incidentally, in the standard 
Reifenberg case this would happen for all $y\in \Sigma_k = f_k(\Sigma_0)$.

\ms\proclaim Corollary 4.14. For $y\in V_{k}^{8}$, 
let $i=i(y) \in J_k$ be such that $y\in 10B_{i,k}$.
Then
$$
|\sigma_k (y) - \pi_{i,k}(y)| \leq C \varepsilon r_k,
\leqno (4.15)
$$
$$
D\sigma_k (y) 
= \sum_{j \in J_k}\theta_{j,k}(y) D\pi_{j,k}
+ \sum_{j \in J_k} \pi_{j,k}(y) D\theta_{j,k}(y),
\leqno (4.16)
$$
and 
$$
|D\sigma_k (y) - D\pi_{i,k}| \leq C \varepsilon.
\leqno (4.17)
$$

\ms
Indeed, since $\psi_k(y)=0$ and $D\psi_k(y)=0$ by (3.12), (4.16)
follows from (4.7) and (4.17) holds because $L(y) = D\pi_{i,k}$.
Finally, $\sum_{j\in J_k} \theta_{j,k}(y) =1$ by (3.14), so
second part of (4.2) yields
$$\eqalign{
\sigma_k (y) - \pi_{i,k}(y)
&= \psi_k(y) y + \sum_{j \in J_k}\theta_{j,k}(y) \,
[\pi_{j,k}(y) -\pi_{i,k}(y)] 
\cr&
= \sum_{j\in J_k} \theta_{j,k}(y) [\pi_{j,k}(y) -\pi_{i,k}(y)] 
}\leqno (4.19)
$$
(recall that $\psi_k(y) = 0$ by (3.12)).
If $\theta_{j,k}(y) \neq 0$, then $y\in 10B_{j,k}$ by (3.9), 
and so $|x_{j,k}-x_{i,k}| \leq 20 r_k$. For each such $j$,
(2.8) says that 
$d_{x_{j,k},100r_k}(P_{i,k},P_{j,k}) \leq \varepsilon$ (as in (4.11)), 
hence $|\pi_{i,k}(y) - \pi_{j,k}(y)| \leq C \varepsilon r_k$ (as in 
(4.12)); then (4.15) follows from (4.19).
\qed

\bsi
{\bf 5. Local Lipschitz graph descriptions of the $\Sigma_k$}
\ms  

In this section we study the local Lipschitz regularity
of the sets $\Sigma_k = f_k(\Sigma_0)$, $k \geq 0$, where
$\Sigma_0$ is our smooth initial comparison surface,
and the $f_k$ are as in (4.1). Thus
$$
\Sigma_{k+1} = \sigma_k(\Sigma_k)
= \sigma_{k} \circ \cdots \circ \sigma_0(\Sigma_0)
\ \hbox{ for } k \geq 0.
\leqno (5.1)
$$
Our main result will be a small Lipschitz graph description 
of $\Sigma_k$ near the $x_{j,k}$. 
It will be easier to state with the following notation 
for boxes. When $x\in \R^n$, $P$ is a $d$-plane through $x$, 
and $R>0$, we define the box $D(x,P,R)$ by
$$\eqalign{
D(x,P,R) &= \big\{ z + w \, ; \, z\in P \cap B(x,R)
\hbox{ and } w\in P^\perp \cap B(0,R) \big\}
\cr&\simeq [P \cap B(x,R)] \times [P^\perp \cap B(0,R)],
}\leqno (5.2)
$$
where $P^\perp$ denotes the $(n-d)$-dimensional
vector space orthogonal to $P$. Also recall that
when $A : P \to P^\perp$ is a Lipschitz mapping,
the graph of $A$ over $P$ is
$$
\Gamma_A = \big\{ z + A(z) \, ; \, z \in P \big\}.
\leqno (5.3)
$$

We put a lot of information together in the next statement,
so that we can prove everything at the same time by induction.

\ms\proclaim Proposition 5.4. There exist constants 
$C_i$, $1 \leq i \leq 7$, such that the 
following holds for all $k\geq 0$ and  $j\in J_k$. 
First, there is a function 
$A_{j,k} : P_{j,k} \cap 49B_{j,k} \to P_{j,k}^{\perp}$, 
of class $C^2$, such that 
$$
|A_{j,k}(x_{j,k})| \leq C_1 \varepsilon r_k,  
\leqno (5.5)
$$
$$
|DA_{j,k}(z)| \leq C_2 \varepsilon
\ \hbox{ for } z\in P_{j,k} \cap 49 B_{j,k}, 
\leqno (5.6)
$$
and, if we denote by $\Gamma_{A,j,k}$ the graph of $A_{j,k}$
over $P_{j,k}$,
$$
\Sigma_{k} \cap D(x_{j,k},P_{j,k},49r_k) = 
\Gamma_{A,j,k} \cap D(x_{j,k},P_{j,k},49r_k). 
\leqno (5.7)
$$
Next, there is a function 
$F_{j,k} : P_{j,k} \cap 40 B_{j,k} \to P_{j,k}^{\perp}$,
of class $C^2$, such that
$$
|F_{j,k}(x_{j,k})| \leq C_3 \varepsilon r_k,  
\leqno (5.8)
$$
$$
|DF_{j,k}(z)| \leq C_4 \varepsilon
\ \hbox{ for } z\in P_{j,k} \cap 40 B_{j,k}, 
\leqno (5.9)
$$
$$
|DF_{j,k}(z)| \leq C_5 \varepsilon
\ \hbox{ for } z\in P_{j,k} \cap 7B_{j,k}, 
\leqno (5.10)
$$
and, if we denote by $\Gamma_{F,j,k}$ the graph of $F_{j,k}$
over $P_{j,k}$,
$$
\Sigma_{k+1} \cap D(x_{j,k},P_{j,k},40r_k) = 
\Gamma_{F,j,k} \cap D(x_{j,k},P_{j,k},40r_k). 
\leqno (5.11)
$$
Finally,
$$
|\sigma_k(y)-y| \leq C_6 \varepsilon r_k
\ \hbox{ for } y \in \Sigma_k
\leqno (5.12)
$$
and
$$
|D\sigma_k(y)-D\pi_{j,k} - \psi_k(y) D\pi_{j,k}^\perp|
\leq C_7 \varepsilon
\ \hbox{ for } y \in \Sigma_k \cap 45B_{j,k}. 
\leqno (5.13)
$$

\ms
The statement is a little more complicated than what we will use,
because we want to be able to prove it easily by induction.
[After this, we will not
 care any more about (5.10) and the
difference between the $C_i$'s.]
We gave different names to the constants $C_i$ so that
we can track more easily their mutual dependences. For this 
reason, we will
 also make sure that the constants $C$ written 
in the proof do not depend on the $C_i$. 
Let us announce now that we will
 be able to choose $C_5$ first,
then $C_2$ and $C_3$, $C_1$, $C_6$ and $C_7$, and then $C_4$,
all large enough depending on the previous ones; we shall
also display the mutual dependence relations between the $C_i$ as 
they show up. Note that we need to take $\varepsilon$ small enough, 
depending on all the $C_i$.

Proposition 5.4 only gives information on $\Sigma_k$ and 
$\Sigma_{k+1}$ near the $x_{j,k}$; however, away from the $x_{j,k}$,
the $\Sigma_k$ stay the same because $\sigma_k = I$ on 
$\R^n \sm V_k^{10}$ by (4.5), so it will be easy to get some control 
there too. See Lemma 6.12 and Proposition 6.15.

We cut the proof of Proposition 5.4 into four smaller steps.

\msi{\bf Step 1.}
We first verify (5.5)-(5.7) for $k = 0$. Let $j \in J_0$ be given, 
and use (2.7) to choose 
$x \in \Sigma_0$ such that $|x-x_{j,0}| \leq \varepsilon$.
Let $F_x$ and $\Gamma_{F_x}$ be as in the description of $\Sigma_0$
as a small Lipschitz graph over $P_x$ (see near (2.4)). We just
want to write $\Gamma_{F_x}$ as a Lipschitz graph over $P_{j,0}$. 
Observe that
$$
d_{x_{i,0},100}(P_{j,0},P_{x}) \leq \varepsilon
\leqno (5.14)
$$
by (2.9). In particular $P_{j,0}$ and $P_x$ make a small angle,
and by (2.4) $\Gamma_{F_x}$ is also the graph of some
function $A = A_{j,0} : P_{j,0} \to P_{j,0}^\perp$. 
Of course $A$ is $C^2$ because $F_x$ is $C^2$, and it is
$3\varepsilon$-Lipschitz, as needed for (5.6).
[See the proof of (5.18) if you are not sure about this.]
Next, $D(x_{j,0},P_{j,0},49) \i B(x,200)$, so (5.7) follows
from (2.5) because $\Gamma_{A,j,0} = \Gamma_{F_x}$.
Finally, note that by (2.4), $\Gamma_{F_x}$ passes within
$\varepsilon$ of $x$, hence within $2\varepsilon$ of $x_{j,0}$;
then $|A(x_{j,0})| \leq C \varepsilon$, because we already know
that $A$ is $3 \varepsilon$-Lipschitz.

\msi{\bf Step 2.}
Next we show that (5.8)-(5.11) for $k$ imply (5.5)-(5.7) for
$k+1$, as it is essentially the same proof.
Let $j \in J_{k+1}$ be given, and use (2.3) to find
$i \in J_k$ such that 
$$
x_{j,k+1} \in 2B_{i,k}.
\leqno (5.15)
$$
By induction assumption, we have a $C^2$ and 
$C_4\varepsilon$-Lipschitz 
function $F_{i,k} : P_{i,k} \cap 40B_{i,k} \to P_{i,k}^\perp$,
with the properties (5.8)-(5.11). Let us use Whitney's theorem to
extend $F_{i,k}$ as a function from $ P_{i,k}$ to  $P_{i,k}^\perp$, 
so that it is still $C^2$ and $CC_4 \varepsilon$-Lipschitz. %vv
%% I added $C$ before $C_4$
%% comment: I prefer not to use the constant 1 when I can avoid
%% and also I don't know whether I can arrange C^2 too.
The truth is that we shall never care about the values
of $F_{i,k}$ out of $P_{i,k} \cap 40B_{i,k}$, but the argument
is simpler to write that way. Since
$$
d_{x_{i,k},20r_k}(P_{i,k},P_{j,k+1}) \leq \varepsilon
\leqno (5.16)
$$
by (5.15) and (2.10), we can find a $C^2$ function 
$A = A_{j,k+1} : P_{j,k+1} \to P_{j,k+1}^\perp$ such that
$$
\Gamma_{A,j,k+1} = \Gamma_{F,i,k}.
\leqno (5.17)
$$
Let us check that 
$$
|DA| \leq (C C_4+C) \varepsilon.
\leqno (5.18)
$$
Let $x\in P_{j,k+1}$ be given, and observe that
if $v$ is a unit vector in the vector space $P_{j,k+1}'$ 
parallel to $P_{j,k+1}$, $|DA(x) \cdot v| = \tan\alpha$,
where $\alpha = {\rm Angle}(v+DA(x)\cdot v,P_{j,k+1}')$.
Thus 
$$\eqalign{
|DA(x)| 
&= \sup_{v\in P_{j,k+1}' \, ; \, |v|=1} |DA(x) \cdot v|
\cr&
= \sup_{v\in P_{j,k+1}' \, ; \, |v|=1} 
\tan{\rm Angle}(v+DA(x)\cdot v,P_{j,k+1}')
\cr&
= \sup_{w\in T\Gamma(y)} \tan{\rm Angle}(w,P_{j,k+1}'),
}\leqno (5.19)
$$
where $T\Gamma(y)$ denotes the tangent plane to $\Gamma_{A,j,k+1}$
at $y = x+A(x)$. Next let $x' \in P_{i,k}$ be such that
$y = x'+F_{i,k}(x')$; then
$$\eqalign{
\sup_{w\in T\Gamma(y)} {\rm Angle}(w,P_{j,k+1}')
&\leq C \varepsilon +
\sup_{w\in T\Gamma(y)} {\rm Angle}(w,P_{i,k}')
\cr&
\leq C \varepsilon + \sup_{v'\in P_{i,k}' \, ; \, |v'|=1} 
\tan{\rm Angle}(v'+DF_{i,k}(x')\cdot v',P_{i,k}')
}\leqno (5.20)
$$
by (5.16), because $T\Gamma(y)$ is also the tangent plane to 
$\Gamma_{F,i,k}$ at $y$ (by (5.17)), and by the same computation 
as for (5.19). Set
$\beta = {\rm Arctan}(C_4 \varepsilon)$; then 
${\rm Angle}(v'+DF_{i,k}(x')\cdot v',P_{i,k}') \leq \beta$
for unit vectors $v'\in P_{i,k}'$, by (5.9) (and maybe because 
we extended $F_{i,k}$ in a $C C_4 \varepsilon$-Lipschitz way), so
$$\eqalign{
|DA(x)| &\leq \tan(C\varepsilon + \beta)
= \tan \beta + \int_{\beta}^{\beta+C\varepsilon} {dt \over \cos^2(t)}
\cr&
\leq \tan \beta + {C \varepsilon \over \cos^2(\beta+C\varepsilon)}
\leq C C_4 \varepsilon + C \varepsilon,
}\leqno (5.21)
$$
as needed for (5.18).

Notice that
$$
D(x_{j,k+1},P_{j,k+1},35r_{k}) \i D(x_{i,k},P_{i,k},40r_{k})
\leqno (5.22)  
$$
by (5.15) and (5.16), so
$$\eqalign{
\Gamma_{A,j,k+1} \cap D(x_{j,k+1},P_{j,k+1},35r_{k})
&= \Gamma_{F,i,k} \cap D(x_{j,k+1},P_{j,k+1},35r_{k})
\cr&= \Sigma_{k+1} \cap D(x_{j,k+1},P_{j,k+1},35r_{k})
}\leqno (5.23)
$$
by (5.17) and (5.11). This is better than (5.7), which
only requires a control on the smaller 
$D(x_{j,k+1},P_{j,k+1},49r_{k+1})$.

Next we check (5.5) and (5.6). We first need to control
$A$ at one point. Set $B = \overline B(x_{i,k},C_3\varepsilon r_{k})$.
By (5.8), $\Gamma_{F,i,k}$ meets $B$; by (5.17), so does
$\Gamma_{A,j,k+1}$. Choose $y_0\in B \cap \Gamma_{A,j,k+1}$
and $z_0\in P_{j,k+1}$ such that $y_0=z_0+A(z_0)$. Then
$$
|A(z_0)| = \dist(y_0,P_{j,k+1}) \leq \dist(y_0,P_{i,k}) 
+ 20 \varepsilon r_k
\leq C_3\varepsilon r_{k} + 20 \varepsilon r_k
\leqno (5.24)
$$
by (5.16), and because $y_0\in B$ and $x_{i,k} \in P_{i,k}$.
Now we check (5.6). Let $x\in P_{j,k+1}\cap 49B_{j,k+1}$ 
be given, set $y = x + A(x) \in \Gamma_{A,j,k+1}$, and use (5.17) 
to find $x'\in P_{i,k}$ such that $y = x' + F_{i,k}(x')$. 
Note that
$$\leqalignno{
|x'-x_{i,k}| &\leq |y-x_{i,k}| \leq 2r_k + |y-x_{j,k+1}|
= 2r_k + |x+A(x)-x_{j,k+1}|
\cr& \leq 2r_k + 49r_{k+1} + |A(x)| 
\leq 2r_k + 49r_{k+1} + |A(z_0)|+|A(x)-A(z_0)|
< 7 r_k  
&(5.25)
}%\leqno (5.25)
$$
because $x_{i,k} \in P_{i,k}$, by (5.15), and because
$|A(x)| < r_{k+1}/10$ by (5.24) and (5.18). Thus $x'\in 7B_{i,k}$, 
(5.10) says that $|DF_{i,k}(x') \leq C_5\varepsilon$, and the
proof of (5.18) shows that $|DA(x)| \leq (C_5 + C)\,\varepsilon$. 
Then (5.6) holds if we choose
$$
C_2 \geq C_5 + C. 
\leqno (5.26)
$$
We may now return to (5.5). Observe that 
$$\eqalign{
|x_{j,k+1}-z_0| &\leq |x_{j,k+1}-x_{i,k}|+|x_{i,k}-y_0|+|y_0-z_0|
\cr& \leq 2r_k + C_3\varepsilon r_{k} + |A(z_0)|
\leq 2r_k + 2C_3\varepsilon r_{k} + 20 \varepsilon r_k < 3r_k
}\leqno (5.27)
$$
by (5.15) and (5.24), so 
$$
|A(x_{j,k+1})| \leq |A(z_0)| + (C_5+C)\,\varepsilon\, |x_{j,k+1}-z_0|
\leq (C_3 + 3C_5 + C) \varepsilon r_k
\leqno (5.28)
$$
because we just checked that $|DA| \leq (C_5 + C)\,\varepsilon$ on 
$P_{j,k+1}\cap 49B_{j,k+1}$ and by (5.24). So (5.5) holds as soon as
$$
C_1 \geq C_3 + 3C_5 + C.
\leqno (5.29)
$$
This completes our verification of (5.5)-(5.7) for
$k+1$ given (5.8)-(5.11) for $k$.

\msi{\bf Step 3.}
Now we assume (5.5)-(5.7) (for $k$) and show that (5.12) and
(5.13) hold. We start with (5.12).
Let $y\in \Sigma_k$ be given. If $y$ lies out of 
$V_k^{10} = \cup_{j\in J_k} 10B_{j,k}$, (4.5) 
says that $\sigma_0(y)=y$ and (5.12) holds trivially.
So we may assume that $y \in V_k^{10}$.
Choose $j \in J_k$ such that $y \in 10 B_{j,k}$, %v
and then let $A_{j,k}$ and $\Gamma_{A,j,k}$ be as in (5.5)-(5.7).
By (5.7), $y \in \Gamma_{A,j,k}$. That is, $y = x + A_{j,k}(x)$
for some $x \in P_{j,k}$. In addition $x \in 10B_{j,k}$
(because it is the orthogonal projection on $P_{j,k}$ of
$y\in 10B_{j,k}$). By (5.5) and (5.6),
$$
\dist(y,P_{j,k}) = |A_{j,k}(x)|
\leq (C_1+10C_2) \varepsilon r_k.
\leqno (5.30)
$$
Next let $i\in J_k$ be such that $\theta_{i,k}(y) \neq 0$;
then $y\in 10B_{i,k}$ by (3.9), and (2.8) says that
$$
d_{x_{j,k},100r_k}(P_{i,k},P_{j,k}) \leq \varepsilon
\leqno (5.31)
$$
(because $y \in 10B_{j,k} \cap 10B_{i,k}$). Thus
$$
\dist(y,P_{i,k}) \leq \dist(y,P_{j,k}) + 100 \varepsilon r_k
\leq (C_1+10C_2+100) \varepsilon r_k.
\leqno (5.32)
$$
Now the first part of (4.2) yields
$$\eqalign{
|\sigma_k(y) - y|
&\leq \sum_{i \in J_k} \theta_{i,k}(y) |\pi_{i,k}(y)-y|
\cr&= \sum_{i \in J_k} \theta_{i,k}(y) \dist(y,P_{i,k})
\leq (C_1+10C_2+100) \, \varepsilon r_k.
}\leqno (5.33)
$$
This yields (5.12), if
$$
C_6 \geq C_1+10C_2+100.
\leqno (5.34)
$$
For (5.13) let $j\in J_k$ and $y \in \Sigma_k \cap 45B_{j,k}$ be given.
If $y \in \R^n \sm V_k^{10}$, i.e., if $|y-x_{i,k}| \geq 10 r_k$
for all $i \in J_k$, then $\psi_k(y) = 1$ by (3.9) and (3.13),
and $D\sigma_k(y) = I$ by (4.5), so (5.13) holds in this case.
So we may assume that $y \in V_k^{10}$, and choose $i \in J_k$
such that $|y-x_{i,k}| < 10r_k$.
Recall from (4.9) and (4.8) that
$$
|D\sigma_k (y) - L(y)| \leq C \varepsilon, 
\leqno (5.35)
$$
where
$$
L(y) 
= \psi_k(y) D\pi_{i,k}^\perp + D\pi_{i,k} + [y - \pi_{i,k}(y)] 
D\psi_k(y),
\leqno (5.36)
$$
so we want to control
$$\eqalign{
A &= L(y) - D\pi_{j,k} - \psi_k(y) D\pi_{j,k}^\perp 
\cr&= [D\pi_{i,k} - D\pi_{j,k}] 
+ \psi_k(y) [D\pi_{i,k}^\perp - D\pi_{j,k}^\perp]
+ [y - \pi_{i,k}(y)] D\psi_k(y).
}\leqno (5.37)
$$
Recall that $y\in 45B_{j,k} \cap 10B_{i,k}$, so (2.8) says that
$$
d_{x_{i,k},100}(P_{i,k},P_{j,k}) \leq \varepsilon,
\leqno (5.38)
$$
and hence
$$
|D\pi_{i,k} - D\pi_{j,k}| + |D\pi_{i,k}^\perp - D\pi_{j,k}^\perp| 
\leq C \varepsilon.
\leqno (5.39)
$$
Thus we control the first two terms of $A$ (recall that  
$0 \leq \psi_k(y) \leq 1$). Next, 
$$\eqalign{
\big| [y - \pi_{i,k}(y)] D\psi_k(y) \big|
&\leq C r_k^{-1} |y - \pi_{i,k}(y)|
= C r_k^{-1} \dist(y,P_{i,k})
\cr&\leq C r_k^{-1} \dist(y,P_{j,k})+C \varepsilon
}\leqno (5.40)
$$
by (3.15) and (5.38). By (5.7), $y\in \Gamma_{A,j,k}$, so
there is an $x\in 45B_{j,k}$ such that $y=x+A_{j,k}(x)$. Then
$$
\dist(y,P_{j,k}) =|A_{j,k}(x)|
\leq |A_{j,k}(x_{j,k})| + 45 C_2 \varepsilon r_k
\leq (C_1+45 C_2) \varepsilon r_k
\leqno (5.41)
$$
by (5.5) and (5.6). Altogether,
$$
|D\sigma_k (y) - D\pi_{j,k} - \psi_k(y) D\pi_{j,k}^\perp|
\leq |D\sigma_k (y) - L(y)| + |A|
\leq C (C_1+C_2+1) \varepsilon
\leqno (5.42)
$$
by (5.37), (5.35), and (5.39)-(5.41). This proves (5.13)
with the constraint that
$$
C_7 \geq C (C_1+C_2+1).
\leqno (5.43)
$$

\msi{\bf Step 4.}
We now come to the main and final part of the induction argument.
This is also the place where we get small constants %v
that we can use to control the other ones.
We still assume (5.5)-(5.7) (for $k$) and prove (5.8)-(5.11), 
with the help of (5.12) and (5.13) that we just proved.
Let $j\in J_k$ be given, and let $A_{j,k}$ and $\Gamma_{A,j,k}$ 
be as in (5.5)-(5.7). To simplify the notation, set
$$
x = x_{j,k}\,, P = P_{j,k}\,, \pi = \pi_{j,k}\,, 
\pi^\perp = \pi_{j,k}^\perp \,, A=A_{j,k}\,, 
\hbox{ and } \Gamma = \Gamma_{A,j,k} \,.
\leqno (5.44)
$$
Define $h : P \cap B(x,44r_k) \to \R^n$ and 
$h_1: P \cap B(x,44r_k) \to P$ by
$$
h(z) = \sigma_k(z+A(z))
\hbox{ and } h_1(z) = \pi \circ h(z)
\ \hbox{ for } z\in P \cap B(x,44r_k).
\leqno (5.45)
$$
Fix $z\in P \cap B(x,44r_k)$, and set $y = z+A(z)$; note that
$y\in \Sigma_k$ by (5.7), (5.5), and (5.6). Then
$$\eqalign{
|h_1(z)-z| &= |\pi \circ h(z) -z| 
= |\pi(\sigma_k(y)) -z|
\cr& = |\pi(\sigma_k(y)) -\pi(y)|
\leq |\sigma_k(y)-y| \leq  C_6 \, \varepsilon r_k
}\leqno (5.46)
$$ 
by (5.45) and (5.12). Let us now say why (5.46) and a tiny 
bit of degree theory imply that
$$
h_1(P \cap B(x,44r_k))
\hbox{ contains } P \cap B(x,43r_k).
\leqno (5.47)
$$ 

Let $\xi \in P \cap B(x,43r_k)$ be given; we want to define a
few mappings from $\d = \d B(x,44r_k) \cap P$ to the
unit sphere $\d B(0,1)$. First observe that
$$
\xi \hbox{ lies out of the segment $[h_1(z),z]$
for } z\in \partial,
\leqno (5.48)
$$ 
simply because $|z-\xi| \geq r_k$, while 
$|h_1(z)-z| \leq C_6 \varepsilon r_k$ by (5.46). Set
$$
u_\xi(z) = (h_1(z)-\xi)/|h_1(z)-\xi|
\ \hbox{ for $z \in \d$.}
\leqno (5.49)
$$
The denominator does not vanish (by (5.48)), and $u_\xi$ is a 
continuous function from $\partial$ to $\d B(0,1)$.
It is homotopic (among mappings from $\d$ to $\d B(0,1)$)
to $v_\xi$, where
$$
v_\xi(z) = [z-\xi]/|z-\xi|
\ \hbox{ for $z \in \d$,}
\leqno (5.50)
$$
because (5.48) allows us to use the obvious linear path 
from $h_1(z)$ to $z$; that is, we can set
$g(z,t) = [(1-t)h_1(z)+tz-\xi]/\big|(1-t)h_1(z)+tz-\xi\big|$
for $0 \leq t \leq 1$ and $z\in \d$ to connect the two mappings.

In turn, $v_\xi$ is homotopic to $z \to [z-x]/|z-x|$ 
(just move $\xi$ to $x$ continuously along the segment $[\xi,x]$,
and observe that $[\xi,x]$ does not meet $\d$), and this last 
mapping is of degree $1$. Thus $v_\xi$ is not homotopic to a constant, 
and neither is $u_\xi$. See for instance [Du]. %%
Therefore  $\xi \in h_1(P \cap B(x,44r_k))$, 
because otherwise we could define a homotopy from $u_\xi$ 
to a constant, by setting $g(z,t) = 
[h_1\big((1-t)z+tx\big)-\xi ]/\big|h_1\big((1-t)z+tx\big)-\xi \big|$ for 
$z\in \d$ and $0 \leq t \leq 1$. Thus (5.47) holds.

\ms
Next we estimate the derivative of $h$ on $P \cap B(x,44r_k)$. 
Still set $y = z+A(z)$ for $z\in P \cap B(x,44r_k)$; by (5.45),
$$
Dh(z) = D\sigma_k(y) \circ [I' + DA(z)],
\leqno (5.51)
$$ 
where $I'$ denotes the injection from $P'$ to $\R^n$, 
where $P'$ denotes the vector space parallel to $P$.
For notational convenience we forget to write the
composition of $DA(z)$ with the canonical injection from
$P^\perp$ to $\R^n$. That is, we now see $DA(z)$ as going 
from the plane $P$ to $\R^n$ (and not just to $P^\perp$). 

Set $Z = D\sigma_k(y) - D\pi - \psi_k(y) D\pi^\perp$;
then $|Z| \leq C_7 \varepsilon$ by (5.13), which applies
to $y$ because
$$
|x-y| \leq |x-z| + |A(z)| \leq 44r_k + |A(z)|
\leq 44r_k + |A(x)| + 44 C_2 \varepsilon r_k < 45 r_k
\leqno (5.52)
$$ 
by (5.6) and (5.5). Now (5.51) yields
$$\eqalign{
Dh(z) &= D\sigma_k(y) \circ [I' + DA(z)]
= [D\pi+\psi_k(y)D\pi^\perp + Z]\circ [I' + DA(z)]
\cr&
= I' + Z \circ I' 
+ [D\pi+\psi_k(y)D\pi^\perp + Z]\circ DA(z)
}\leqno (5.53)
$$ 
because $D\pi \circ I' = I'$ and $D\pi^\perp \circ I' = 0$. 
Now $|Z \circ I'| \leq |Z| \leq C_7 \varepsilon$
and 
$$\eqalign{
\big|[D\pi+\psi_k(y)D\pi^\perp + Z]\circ DA(z)\big|
&\leq \big|[D\pi+\psi_k(y)D\pi^\perp + Z]\big|\,\big| DA(z)\big|
\cr&\leq (1+C_7\varepsilon) \big| DA(z)\big| \leq 2 C_2 \varepsilon
}\leqno (5.54)
$$
by (5.6), so we get that
$$
|Dh(z) - I'| \leq (2C_2+C_7) \varepsilon
\ \hbox{ for } z\in P \cap B(x,44r_k).
\leqno (5.55)
$$

Return to the study of $h_1 = \pi \circ h$ (see (5.45)).
Since $Dh_1(z) = D\pi \circ Dh(z)$, we see that
$$
|Dh_1(z) - I'| \leq (2C_2+C_7) \varepsilon
\ \hbox{ for } z\in P \cap B(x,44r_k).
\leqno (5.56)
$$
Let us now check that, because of (5.56), 
$$
\hbox{ the restriction of $h_1$ to $P \cap B(x,44r_k)$
is injective.}
\leqno (5.57)
$$
Let $z_1,z_2 \in P \cap B(x,44r_k)$ be given; we 
apply the fundamental theorem of calculus to the function
$h_1(z)-z$ on the segment $[z_1,z_2]$, and get that
$|h_1(z_1) - h_1(z_2) - z_1 -z_2| \leq (2C_2+C_7) \varepsilon |z_1-z_2|$
because $|D[h_1 -I']| \leq (2C_2+C_7) \varepsilon$ by (5.56).
This is impossible if $h_1(z_1) = h_1(z_2)$ but $z_1 \neq z_2$.
So (5.57) holds.

For $w\in P \cap B(x,43r_k)$, (5.47) says that we can find 
$z\in P \cap B(x,44r_k)$ such that $h_1(z)=w$, and (5.57) says 
that it is unique. So we can define 
$$
h_1^{-1} : P \cap B(x,43r_k) \to P \cap B(x,44r_k)
\leqno (5.58)
$$
as follows: $h_1^{-1}(w)$ is the only
$z\in P \cap B(x,44r_k)$ such that $h_1(z)=w$. 
Recall that all our mappings are of class $C^2$; then
(5.56) and the local inversion theorem say that 
$$
h_1^{-1} \hbox{ is of class $C^2$ on $P \cap B(x,43r_k)$, and } 
\big|D[h_1^{-1}]\big| \leq 1+(3C_2+2C_7)\varepsilon.
\leqno (5.59)
$$
Define $F$ on $P \cap B(x,43r_k)$ by
$$
F = \pi^\perp \circ h \circ h_1^{-1}.
\leqno (5.60)
$$
Obviously $F$ is also $C^2$, and
$$
DF(w) = D\pi^\perp \circ Dh(z) \circ D[h_1^{-1}](w)
= D\pi^\perp \circ Dh(z) \circ [Dh_1(z)]^{-1},
\leqno (5.61)
$$
where we set $z = h_1^{-1}(w) \in P \cap B(x,44r_k)$.
Then 
$$\eqalign{
|DF(w)| &\leq (1+(3C_2+2C_7)\varepsilon) |D\pi^\perp \circ Dh(z)|
\leq 2 |D\pi^\perp \circ Dh(z)|
\cr&= 2 |D\pi^\perp \circ [Dh(z)-I']|
\leq 2(2C_2+C_7)\varepsilon
}\leqno (5.62)
$$
by (5.59), because $D\pi^\perp \circ I'=0$, and by (5.55).

We want to show that (5.8)-(5.11) hold for $F=F_{j,k}$,
and we start with (5.11). Denote by $\Gamma_F$ the graph of 
$F$ over $P$; thus $\Gamma_F = \Gamma_{F,j,k}$ with the notation
of (5.11). Let us first check that 
$$
\Gamma_F \cap \pi^{-1}(B(x,40r_k)) \i \Sigma_{k+1}.
\leqno (5.63)
$$
The fact that $\Gamma_F \cap D(x,P,40r_k) \i 
\Sigma_{k+1} \cap D(x,P,40r_k)$ will immediately  
follow, by restricting to $D(x,P,40r_k)$ (see
the definition (5.2)).

Let $\xi$ lie in $\Gamma_F \cap \pi^{-1}(B(x,40r_k))$. 
Then $\xi = w + F(w)$ for some $w\in P \cap B(x,40r_k)$,
and so $z=h_1^{-1}(w)$ is defined. Then $w=h_1(z)$,
and $F(w) = \pi^\perp \circ h(z)$ by (5.60). 
Thus $\xi = w + F(w) = h_1(z) +\pi^\perp \circ h(z) = 
\pi \circ h(z) +\pi^\perp \circ h(z) = h(z) = 
\sigma_k(z+A(z))$ by (5.45). But $z \in P \cap B(x,44r_k)$,
so $z+A(z) \in D(x,P,44r_k)$ by (5.5) and (5.6), and now
(5.7) says that $z+A(z) \in \Sigma_k$. Thus
$\xi \in \Sigma_{k+1}$, and (5.63) holds.

Next we claim that
$$
\Sigma_{k+1} \cap D(x,P,40r_k) \i \Gamma_F \, .
\leqno (5.64)
$$
Let $\xi \in \Sigma_{k+1} \cap D(x,P,40r_k)$
be given, and let $y\in \Sigma_k$ be such that 
$\xi = \sigma_k(y)$. By (5.12), 
$|\xi-y| = |\sigma_k(y) - y| \leq C_6\varepsilon r_k$,
so $y \in \Sigma_k \cap D(x,P,41r_k)$, and by (5.7) it lies on 
the graph of $A$. Thus $y = z + A(z)$ for some 
$z=\pi(y)\in P \cap B(x,41r_k)$. Now $\xi =  \sigma_k(y) = h(z)$ 
by (5.45), so $\pi(\xi) = \pi \circ h (z) = h_1(z)$ by (5.45) and 
(since $\pi(\xi) \in B(x,40r_k)$ because $\xi \in D(x,P,40r_k)$)
we get that $z=h_1^{-1}(\pi(\xi))$. Finally
$\pi^\perp(\xi) = \pi^\perp(h(z)) = F(\pi(\xi))$
by (5.60), which means that $\xi$ lies on the graph of
$F$, as needed for (5.64).

As was just observed, (5.11) follows from (5.63) and (5.64).
We deduce (5.9) from (5.62) as soon as
$$
C_4 \geq 4C_2+2C_7,
\leqno (5.65)
$$
so we are left with (5.8) and (5.10) to check.

First we estimate $|F(h_1(x))|$, where $x$ %v
still denotes the center of $B_{j,k}$ (see (5.44)).
Notice that $x+A(x) \in \overline B(x,C_1 \varepsilon r_k)$ by
(5.5), so it lies in $\Sigma_k$ by (5.7), 
and $h(x) = \sigma_k(x+A(x))$ lies in $\Sigma_{k+1}$ Next, 
$$
|h_1(x) -x| \leq C_6 \varepsilon r_k
\leqno (5.66)
$$
by (5.46), so $h_1^{-1}(h_1(x)) = x$ %v
(see near (5.58)) and
$$\eqalign{
|F(h_1(x))| &= |\pi^\perp \circ h \circ h_1^{-1}(h_1(x))|
=|\pi^\perp \circ h(x)| = |\pi^\perp(\sigma_k(x+A(x)))|
\cr& \leq |\pi^\perp(x+A(x))| + |\sigma_k(x+A(x))-(x+A(x))|
\cr&
\leq |A(x)| + C_6 \varepsilon r_k
\leq (C_1+C_6)\,\varepsilon r_k
}\leqno (5.67)
$$
by (5.60), (5.12), and (5.5). This is not yet good enough
for (5.5), but is a first step.

Let $z\in P \cap 7B_{j,k}$ be given, and set $w = z + F(z)$.
By (5.62), (5.66), and (5.67), $|F(z)| < r_k/10$.
By (5.63), $w \in \Sigma_{k+1}$ and we can find $y \in \Sigma_k$ 
such that $w = \sigma_k(y)$. Then 
$$
|z-y| \leq |F(z)| + |w-y| = |F(z)| + |\sigma_k(y)-y| 
\leq r_k/10 + C_6 \varepsilon r_k < r_k,
\leqno (5.68)
$$
so $y \in 8B_{j,k}$ because $z\in 7B_{j,k}$.
Corollary 4.14 applies with $i(y)=j$, and says that
$$
|\sigma_k(y)-\pi(y)| \leq C \varepsilon r_k
\ \hbox{ and } \
|D\sigma_k(y)-D\pi| \leq C \varepsilon.
\leqno (5.69)
$$
In particular,
$$
|F(z)| = \dist(w,P) \leq |w-\pi(y)|
= |\sigma_k(y)-\pi(y)| \leq C \varepsilon r_k,
\leqno (5.70)
$$
which implies (5.8) by taking $z=x$ and if $C_3$ is large enough.

Return to the general $z\in P \cap 7B_{j,k} \,$; since 
$y \in \Sigma_k \cap 8B_{j,k}$, (5.7) says that $y$ lies on 
the graph of $A$. That is, $y = \pi(y) + A(\pi(y))$ and now
$$
h_1(\pi(y)) = \pi \circ h(\pi(y)) = \pi \circ \sigma_k(y)
= \pi(w) = z
\leqno (5.71)
$$
by (5.45) and the various definitions. Thus $\pi(y) = h_1^{-1}(z)$
(see near (5.58) and recall that $\pi(y) \in P \cap 8B_{j,k}$
because $y \in 8B_{j,k}$). Since
$$
F = \pi^\perp \circ h \circ h_1^{-1} = \pi^\perp \circ \sigma_k
\circ [I'+A] \circ h_1^{-1}
\leqno (5.72)
$$
by (5.60) and (5.45) and where $I'$ is the restriction
of $I$ to $P$, we get that
$$\eqalign{
DF(z) &=
D\pi^\perp \circ D\sigma_k(y) \circ [I'+DA(\pi(y))] 
\circ D[h_1^{-1}](z)\cr
&= D\pi^\perp \circ (D\sigma_k(y)-D\pi) \circ [I'+DA(\pi(y))] 
\circ D[h_1^{-1}](z)}
\leqno (5.73)
$$
because $h_1^{-1}(z) = \pi(y)$, $(I'+A)(\pi(y)) = y$, and
$D\pi^\perp \circ D\pi = 0$. Hence (5.69) yields %v
$$
|DF(z)| \leq C \varepsilon |I'+DA(\pi(y))| \, |D[h_1^{-1}](z)|
\leq 2 C\varepsilon
\leqno (5.74)
$$
by (5.6) and (5.59). This proves (5.10) if $C_5$ is large enough,
and completes our verification of (5.8)-(5.11). 
Our proof of Proposition 5.4 by induction is now complete too.
\qed

\bsi
{\bf 6. Reifenberg-flatness of the image}
\ms

We continue to assume that $(\Sigma_0, \{ B_{j,k} \}, \{ P_{j,k} \})$
is a coherent collection of balls and planes, as in Definition 2.11,
with $\varepsilon > 0$  is small enough (depending only on $n$).

We defined a sequence of functions $f_n : \R^n \to \R^n$ in
Section 4, and by (4.4) the sequence converges uniformly and
we can define $f : \R^n \to \R^n$ by
$$
f(x) = \lim_{k \to +\infty} f_k(x)
\ \hbox{ for } x \in \R^n.
\leqno (6.1)
$$

In this section we want to record some of the geometric properties 
of $\Sigma = f(\Sigma_0)$, and in particular its Reifenberg-flatness.
Of course the main ingredient will be Proposition 5.4.

Let us first check that
$$
\Sigma = f(\Sigma_0) \hbox{ contains } E_\infty,
\leqno (6.2)
$$
where $E_\infty$ is the accumulation set of the
centers $x_{j,k}$, i.e., 
$$\eqalign{
E_\infty &= \big\{ x\in \R^n \, ; \, 
x \hbox{ can be written as }
x = \lim_{m\to +\infty} x_{j(m),k(m)}, 
\hbox{ with } k(m) \in \Bbb N
\cr&\hskip 2.5cm \hbox{ and }
j(m) \in J_{k(m)} \hbox{ for } m \geq 0,  \hbox{ and }
\lim_{m \to  +\infty} k(m) = +\infty \big\}.
}\leqno (6.3)
$$

This is easy. Let $x \in E_\infty$ be given, and write
$x= \lim_{m\to +\infty} x_{j(m),k(m)}$ as above. 
By Proposition 5.4, $B_{j(m),k(m)}$ meets $\Sigma_{k(m)+1}$, 
so we can find $t_m \in \Sigma_0$ such that
$f_{k(m)}(t_m) \in B_{j(m),k(m)}$, i.e., 
$|f_{k(m)}(t_m)-x_{j(m),k(m)} |< r_{k(m)}$. Then
$$
\lim_{k \to +\infty} f_{k(m)}(t_m) = x
\leqno (6.4)
$$
because $r_{k(m)}$ tends to $0$.
Since $|f_k(y)-y| \leq 12$ by (4.4), the sequence $\{ t_m \}$
is bounded, and some subsequence $\{ t_{m(p)} \}$ converges.
Set $t = \lim_{p \to +\infty} t_{m(p)}$. Then
$|f(t)-x| \leq |f(t)-f(t_{m(p)})| + |f(t_{m(p)})-f_{k(m(p))}(t_{m(p)})|
+ |f_{k(m(p))}(t_{m(p)})-x|$; the first term tends to $0$ because $f$
is continuous, the second one because $f$ is the uniform limit
of the $f_{k(m(p))}$, and the third one by (6.4). So $f(t) = x$, and 
(6.2) holds.

\ms
We shall also need the following description of trajectories.

\ms\proclaim Lemma 6.5.
Let $z\in \R^n$ be given, and set $z_k = f_k(z)$
for $k \geq 0$. Then
$$
\hbox{ if $z_k \in \R^n \sm V_k^{10}$ for some $k \geq 0$, then 
$z_l = z_k \in \R^n \sm V_l^{10}$ for $l \geq k$;}
\leqno (6.6)
$$
$$
\hbox{if $z \in \Sigma_0$ and $z_k \in V_{k}^{10}$ for some $k \geq 1$, 
then $z_l \in V_l^{4}$ for $0 \leq l \leq k-1$.} 
\leqno (6.7)
$$

\ms
Recall from (5.12) that 
$$
|\sigma_k(y)-y| \leq C \varepsilon r_k
\ \hbox{ for } k \geq 0 \hbox{ and } y \in \Sigma_k,
\leqno (6.8)
$$
which easily implies by iteration that
$$
|f(x)-f_k(x)| \leq C \varepsilon r_k
\ \hbox{ for } k \geq 0 \hbox{ and } x \in \Sigma_0.
\leqno (6.9)
$$
Also recall from (4.5) that
$$
\sigma_k(y) = y \ \hbox{ when }
y \in \Sigma_k \sm V_k^{10}, 
\leqno (6.10)
$$
where $V_k^{10} = \cup_{j \in J(k)} 10B_{j,k}$ is as in (2.2).

If $z_k \in \R^n \sm V_k^{10}$, (6.10) says that $z_{k+1} = z_k$.
Then $z_{k+1} \in \R^n \sm V_{k+1}^{10}$, because otherwise we 
could find $j \in J_{k+1}$ such that $|z_{k+1}-x_{j,k+1}| \leq
10 r_{k+1}$, and since (2.3) gives $i \in J_k$ such that
$|x_{i,k}-x_{j,k+1}| \leq 2r_k$, we would get that 
$|z_{k}-x_{i,k}| = |z_{k+1}-x_{i,k}| \leq 3 r_k$, a contradiction.
So $z_{k+1} \in \R^n \sm V_{k+1}^{10}$, and we can iterate
and get (6.6).

Now suppose that $z \in \Sigma_0$ and $z_k \in V_k^{10}$ 
for some $k \geq 1$. Thus $|z_k-x_{j,k}| \leq 10r_k$
for some $j \in J_k$, and by (2.3) we can find $i \in J_{k-1}$ 
such that $|x_{i,k-1}-x_{j,k}| \leq 2r_{k-1}$ and hence
$|x_{i,k-1}-z_k| \leq 3r_{k-1}$.

By (6.8), $|z_k - z_{k-1}| = 
|\sigma_k(z_{k-1})-z_{k-1}| \leq C \varepsilon r_k$
because $z_{k-1} \in \Sigma_{k-1}$, and so
$|x_{i,k-1}-z_{k-1}| < 4r_{k-1}$. That is, $z_{k-1} \in V_{k-1}^4$.
The previous values of $l$ are now obtained by induction applied to $z_{k-l}$;
(6.7) and Lemma 6.5 follow.
\qed

\ms
As a simple consequence of Lemma 6.5, let us check that for $k \geq 0$,
$$
\Sigma \sm V_{k}^{11} = \Sigma_{l} \sm V_{k}^{11}
= \Sigma_{k} \sm V_{k}^{11}
\ \hbox{ for all } l \geq k.
\leqno (6.11)
$$
Indeed, if $y \in \Sigma \sm V_{k}^{11}$ and if $z \in \Sigma_0$ 
is such that $f(z) = y$, then $f_k(z)\not \in V_{k}^{10}$ by (6.9), 
and (6.10) says that $f_{l}(z) = f_k(z)$ 
for all $l \geq k$. Hence $y = f(z) = f_k(z) = f_l(z)$ for
$l \geq k$, and $y \in \Sigma_l$. Thus $\Sigma \sm V_{k}^{11}
\i \Sigma_l$ for $l \geq k$.

Conversely, let $y\in \Sigma_l \sm V_{k}^{11}$ for some $l \geq k$,
and let $z \in \Sigma_0$ be such that $f_l(z) = y$. By 
repeated uses of (6.8), $|f_k(z)-f_l(z)| \leq C \varepsilon r_k$,
so $f_k(z) \in \R^n \sm V_{k}^{10}$, and (6.10) says that 
$f_{m}(z) = f_k(z)$ for $m \geq l$. Then $f(z) = f_l(z)$
and $y = f(z) \in \Sigma$. So $\Sigma_l \sm V_{k}^{11}
\i \Sigma$. This completes our proof of (6.11).

\ms
Before proving the Reifenberg flatness of $\Sigma$ we first
complete the Lipschitz description of $\Sigma_k$ given in 
Proposition 5.4.

\ms\proclaim Lemma 6.12.
For $k \geq 0$ and $y\in \Sigma_k$, there is an affine $d$-plane
$P$ through $y$ and a $C\varepsilon$-Lipschitz and $C^2$ function 
$A : P \to P^\perp$ such that 
$$
\Sigma_k \cap B(y,19r_k) = \Gamma \cap B(y,19r_k).
\leqno (6.13)
$$
where $\Gamma$ denotes the graph of $A$ over $P$.

\ms
When $y \in V_k^{30}$, we choose $j\in J_k$ such that
$y\in 30B_{j,k}$ and apply Proposition 5.4. The plane
$P = P_{j,k}$ and the function $A = A_{j,k}$ satisfy the
desired properties, by (5.5)-(5.7).

So we may assume that $y \in \Sigma_k \sm V_k^{30}$.
Let $m \geq 0$ denote the smallest
integer such that $y \in \R^n \sm V_{m}^{30}$; 
we know that $m \leq k$, so $r_k \leq r_m$ and 
$B(y,19r_k) \i \R^n\sm V_{m}^{11}$, so (6.11) says that
$$
\Sigma_k \cap B(y,19r_k) = \Sigma_{m} \cap B(y,19r_k).
\leqno (6.14)
$$
If $m=0$, let $x\in \Sigma_0$ be such that $f_k(x) = y$;
then $|x-y| \leq C \varepsilon$ by (6.8). We take
$P = P_x$ and $A = F_x$, then $A$ is $\varepsilon$-Lipschitz
by (2.4) and (6.13) follows from (2.5) and (6.14).
If $m>0$, $y\in V_{m-1}^{30}$ by minimality of $m$, so we can 
choose $j \in J_{m-1}$ such that $y\in 30B_{j,m-1}$. Then 
$B(y,19r_k) \i 32B_{j,m-1} \i 49B_{j,m-1}$; we take 
$P = P_{j,m-1}$ and $A = A_{j,m-1}$ as in Proposition~5.4,
and (6.13) follows from (5.7) and (6.14).
\qed

\ms
We are finally ready to show that $\Sigma$ is 
$C\varepsilon$-Reifenberg-flat (as defined just after 
the statement of Theorem 2.15).

\ms\proclaim Proposition 6.15.
For $z\in \Sigma$ and $0 < t \leq 1$, there is an affine $d$-plane
$P(z,t)$ through $z$ such that
$$
d_{z,t}(\Sigma,P(z,t)) \leq C \varepsilon.
\leqno (6.16)
$$

\ms
Recall that $d_{z,t}$ is the normalized local Hausdorff distance %v
defined by (1.7). Let $z\in \Sigma$ and $0 < t \leq 1$
be given, and choose $k$ such that $10^{-1}t < r_k \leq t$;
thus $k \geq 0$. Write $z = f(x)$ for some $x\in \Sigma_0$
and set $y = f_k(x)$; note that $|y-z| \leq C \varepsilon r_k$
by (6.9), so $B(z,t) \i B(y,11r_k)$.

Let $P$ be as in Lemma 6.12, and let us check that (6.16) 
holds for $P(z,t) = P$. For $z'\in \Sigma \cap B(z,t)$, 
let $x'\in \Sigma_0$ be such that $f(x')=z'$, and set
$y'=f_k(x')$. Then $|y'-z'| \leq C \varepsilon r_k$ by (6.9),
so $y'\in \Sigma_k \cap B(y,12r_k)$. Lemma 6.12 says that
$\dist(y',P) \leq C \varepsilon r_k$ (recall that $\Gamma$
is a small Lipschitz graph over $P$ that passes through $y$
because $y\in \Sigma_k$). So $\dist(z',P) \leq 
C \varepsilon r_k \leq C \varepsilon t$, as needed.

Conversely, if $w \in P \cap B(z,t)$, Lemma 6.12
gives $y'\in \Sigma_k$ such that $|y'-w| \leq C \varepsilon r_k$;
let $x'\in \Sigma_0$ be such that $f_k(x')=y'$ and set %v
$z'=f(x') \in \Sigma$. Then $|z'-y'| \leq C \varepsilon r_k$ by (6.9),
so $|z'-w| \leq C \varepsilon r_k \leq C \varepsilon t$, as needed 
again. This proves (6.16); Proposition 6.15 follows.
\qed

\msi{\bf Remark 6.17.}
If we accept to use the standard Reifenberg theorem, we are
ready to prove Theorem 1.10 now. Indeed, let $E$ be as in the
statement. Take $\Sigma_0 = P(0,10)$ and, for each $k \geq 0$,
choose a maximal collection $\{ x_{j,k} \}$, $j\in J_k$, in
$E$, subject to the constraint that $|x_{i,k}-x_{j,k}| \geq r_k$,
as in (2.1). As was observed below (2.3), our constraint (2.3)
is satisfied, just because $E \i V_{k}^2$ by maximality.

The Lipschitz properties (2.4)-(2.6) of $\Sigma_0$ are 
satisfied, with $P_x = P(0,10)$ and $F_x = 0$, just because 
$\Sigma_0 = P(0,10)$, and (2.7) (although with the constant 
$C\varepsilon$) follows from (1.6) for $P(0,10)$.

For $j \in J_k$, set $P_{j,k} = P(x_{j,k},10r_k)$;
then (2.8)-(2.10) (with $C\varepsilon$ again)
follow from our coherence conditions (1.8) and (1.9), and 
the triangle inequality. So we get a coherent set of balls
and planes, and we can define $f$ as above.
The reader may be surprised that we only use (1.6) for $P(0,10)$,
but what happens is that the coherence conditions (1.8) and (1.9),
plus the fact that $P(x_{j,k},10r_k)$ contains $x_{j,k}$, force
$E$ to stay close to the $P_{j,k}$ anyway.

By (6.2) and the fact that $E \i V_k^2$ by maximality, we deduce that
$E \i \Sigma$, which by Proposition 6.15 is Reifenberg-flat with a 
constant less than $C \varepsilon$. The existence of a bi-H\"older
mapping as in Theorem 1.10 is the standard Reifenberg Theorem, applied
to $\Sigma$.

The construction of Sections~3-5 is the main part of Reifenberg's
topological disc theorem. In the following sections we get precise distortion 
estimates for the $\sigma_k$'s. While they are not needed for the proof of Theorem 1.10
they yield very useful information.

\bsi 
{\bf 7. Distortion estimates for $D\sigma_k$}
\ms  
We want to see how much our mappings $f_k$ distort
lengths and distances, and since Section 5 gives a good 
local Lipschitz description of the $\Sigma_k$, it will often be
enough to control the derivative $D\sigma_k$. We shall mostly
worry about the effect of $D\sigma_k$ on the vector
space $T\Sigma_k(y)$ parallel to the tangent plane to 
$\Sigma_k$ at $y\in \Sigma_k$, because anyway we shall later
replace $f$ on $\R^n \sm \Sigma_0$ with another function $g$.
We start with a simple estimate that holds everywhere.

\ms\proclaim Lemma 7.1. 
For $k \geq 0$, $\sigma_k$ is a $C^2$-diffeomorphism
from $\Sigma_k$ to $\Sigma_{k+1}$ and, for $y\in \Sigma_k$,
$$
D\sigma_k(y) : T\Sigma_k(y) \to T\Sigma_{k+1}(\sigma_k(y))
\hbox{ is bijective and $(1+C\varepsilon)$-bi-Lipschitz.}
\leqno (7.2)
$$
In addition,
$$
|D\sigma_k(y) \cdot v - v| \leq C \varepsilon |v|
\ \hbox{ for $y\in \Sigma_k$ and } v \in T\Sigma_k(y)
\leqno (7.3)
$$
and 
$$
|\sigma_k(y)-\sigma_k(y') - y+y'| \leq C \varepsilon |y-y'|
\ \hbox{ for } y,y'\in \Sigma_k.
\leqno (7.4)
$$

\ms
We already know from Lemma 6.12 that the $\Sigma_k$ are
$C^2$ submanifolds, and $\sigma_k : \Sigma_k \to \Sigma_{k+1}$ 
is smooth by construction and surjective by definition.
We will
 also know that it is injective as soon as we prove
(7.4), and that it is a diffeomorphism as soon as we check (7.2).

Let us check (7.3) first. Recall from (5.13) that for $j\in J_k$
$$
|D\sigma_k(y)-D\pi_{j,k} - \psi_k(y) D\pi_{j,k}^\perp|
\leq C \varepsilon
\ \hbox{ for } y \in \Sigma_k \cap 45B_{j,k} %v
\leqno (7.5)
$$
but also, from the local description of $\Sigma_k$ in
(5.5)-(5.7), that $T\Sigma_k(y)$ makes an angle less than
$C \varepsilon$ with $P_{j,k}$. Thus, if $v \in T\Sigma_k(y)$, 
$|D\pi_{j,k}\cdot v + \psi_k(y) D\pi_{j,k}^\perp\cdot v - v| 
\leq C \varepsilon |v|$. That is, (7.3) holds for
$y\in V_k^{45} = \cup_{j\in J_k} 45B_{j,k}$. 
On the other hand, it is also trivial on $\Sigma_k \sm V_k^{10}$,
because there $D\sigma_k(y) = I$ by (4.5). So (7.3) holds. 

It immediately follows from (7.3) that 
$D\sigma_k(y) : T\Sigma_k(y) \to T\Sigma_{k+1}(\sigma_k(y))$
is $(1+C\varepsilon)$-bi-Lipschitz, and in particular injective.
It is also surjective, since the two spaces have the same dimension.
So (7.2) holds.

Finally let $y,y' \in \Sigma_k$ be given. If $|y'-y| \geq r_k$,
(7.4) holds because $|\sigma_k(y)-\sigma_k(y') - y+y'|
\leq |\sigma_k(y)-y|+|\sigma_k(y')-y'|\leq C \varepsilon r_k$ by (6.8).
Otherwise, Lemma 6.12 gives a $C^1$ path $\gamma : I \to \Sigma_k$ 
that goes from $y$ to $y'$, and with 
${\rm length}(\gamma) = \int_I |\gamma'(t)| dt 
\leq 2|y'-y|$. Then
$$\eqalign{
|\sigma_k(y)-\sigma_k(y') - y+y'|
&= \Big| \int_I (\sigma_k \circ \gamma)'(t)dt
-\int_I \gamma'(t)dt \Big|
\cr&\hskip -2cm
= \Big| \int_I \big[D\sigma_k(\gamma(t)) \cdot \gamma'(t)dt
-\gamma'(t) \big] \, dt \Big|
\leq C \varepsilon \int_I |\gamma'(t)|
\leq 2C \varepsilon |y-y'|
}\leqno (7.6)
$$
by (7.3); (7.4) and Lemma 7.1 follow. 
\qed

\ms
Lemma 7.1 is essentially all that will be needed
for the bi-H\"older estimates in Reifenberg's classical theorem
or Theorem 1.10, but for our bi-Lipschitz results we need more 
precise estimates.
For the rest of this section, we shall concentrate on what 
happens at $y\in V^8_k$, because then Corollary~4.14 applies 
and 
$$
D\sigma_k(y) 
= \sum_{j \in J_k}\theta_{j,k}(y) D\pi_{j,k}
+ \sum_{j \in J_k} \pi_{j,k}(y) D\theta_{j,k}(y)
\leqno (7.7)
$$
by (4.16). We want to improve slightly over the estimates
in Corollary 4.14, and in particular keep track of the places where
the distances between the $P_{j,k}$ are significantly smaller than 
$\varepsilon$. So we set
$$
\varepsilon_k(y) = \sup\big\{ 
d_{x_{i,k},100r_k}(P_{j,k},P_{i,k}) \, ; \, i,j \in J_k, \,
y \in 10 B_{i,k} \cap 10 B_{j,k} \big\} 
\leqno (7.8) 
$$
for $k \geq 0$ and $y \in V_k^{10}$ (when the supremum
does not concern an empty set), and $\varepsilon_k(y) = 0$
for $y \in \R^n \sm V_k^{10}$. 
Note that $\varepsilon_k(y) \leq \varepsilon$, by (2.8).

\ms\proclaim Lemma 7.9.
For $k \geq 0$ and $y\in \Sigma_k \cap V^8_k$, choose $i\in J_k$
such that $|y-x_{i,k}| \leq 10r_k$. Then
$$
|\sigma_{k}(y)-\pi_{i,k}(y)| \leq C \varepsilon_k(y) \, r_k,
\leqno (7.10)
$$
$$
|D\sigma_{k}(y)-D\pi_{i,k}| \leq C \varepsilon_k(y),
\leqno (7.11)
$$
and 
$$
{\rm Angle}(T\Sigma_{k+1}(\sigma_k(y)),P_{i,k})) 
\leq C \varepsilon_k(y).
\leqno (7.12)
$$

\ms
Indeed recall that 
$$
\sigma_{k}(y)-\pi_{i,k}(y)
=\sum_{j \in J_k}\theta_{j,k}(y) [\pi_{j,k}(y)-\pi_{i,k}(y)]
\leqno (7.13)
$$
by the second part of (4.2) and (3.12). Next
$$
d_{x_{i,k},100r_k}(P_{i,k},P_{j,k}) \leq \varepsilon_k(y)
\ \hbox{ when $\theta_{j,k}(y) \neq 0$,}
\leqno (7.14)
$$
because then $y \in 10 B_{i,k} \cap 10 B_{j,k}$ and
by the definition (7.8). Hence 
$$
|\pi_{i,k}(y)-\pi_{j,k}(y)| \leq C\varepsilon_k(y) \, r_k
\leqno (7.15)
$$
for $j$  as above, and (7.10) follows from (7.13) because 
$\sum_{j \in J_k}\theta_{j,k}(y)=1$. Next,
$$
D\sigma_{k}(y)-D\pi_{i,k}
= \sum_{j \in J_k}\theta_{j,k}(y) [D\pi_{j,k}-D\pi_{i,k}]
+ \sum_{j \in J_k} [\pi_{j,k}(y)-\pi_{i,k}(y)] D\theta_{j,k}(y)
\leqno (7.16)
$$
by (7.7) and again because $\sum_{j \in J_k}\theta_{j,k} = 1$ near $y$.
By (7.14), $|D\pi_{j,k}-D\pi_{i,k}| \leq C \varepsilon_k(y)$
when $\theta_{j,k}(y) \neq 0 \,$; (7.11) follows from this, (7.15), 
and the fact that $|D\theta_{j,k}(y)| \leq C r_k^{-1}$ by (3.10).

Finally we check (7.12). Let $w \in T\Sigma_{k+1}(\sigma_k(y))$
be given, and write $w = D\sigma_k(y) \cdot v$ for some 
$v \in T\Sigma_k(y)$. By (7.2), $|v| \leq (1+C \varepsilon)|w|$.
Denote by $P'_{i,k}$ the vector $d$-plane parallel to $P_{i,k}$;
then
$$
\dist(w,P'_{i,k}) \leq |w-D\pi_{i,k}(v)| = 
|[D\sigma_{k}(y)-D\pi_{i,k}]\cdot v| \leq C \varepsilon_k(y)|v|
\leq 2 C \varepsilon_k(y) |w|
\leqno (7.17)
$$ 
by (7.11). Thus every point of $T\Sigma_{k+1}(\sigma_k(y))$
lies close to $P'_{i,k} \,$; (7.12) follows because
$T\Sigma_{k+1}(\sigma_k(y))$ and $P'_{i,k}$ are both 
$d$-dimensional. This proves Lemma 7.9.
\qed

\ms
We shall obtain better estimates (of order $2$)
in some directions, but in terms of the following
numbers that also take the previous generation 
into account. Set
$$\eqalign{
\varepsilon'_k(y) = \sup\big\{ &
d_{x_{i,l},100r_{l}}(P_{j,k},P_{i,l}) \, ; \,
j\in J_k, \, l \in \{ k-1, k \}, 
\cr& \hskip 3.5cm
i\in J_{l}, \hbox{ and } y \in 10 B_{j,k} \cap 11 B_{i,l} 
\big\}  
}\leqno (7.18) 
$$
for $k\geq 1$ and $y \in  V_k^{10}$ and $\varepsilon'_k(y) = 0$ 
when $y\in \R^n \sm V_k^{10}$ (when there are no pairs $(j,k)$ 
as above). These are the same numbers as in (2.21) and (2.22).
Notice that 
$$
\varepsilon_n(y) \leq \varepsilon'_n(y) \leq C \varepsilon,
\leqno (7.19) 
$$
by (7.8), (2.8), and (2.10). Let us first control some angles and 
distances.

\ms\proclaim Lemma 7.20.
For $k \geq 1$ and $y\in \Sigma_k \cap V^8_k$, 
choose $i\in J_k$ such that $|y-x_{i,k}| \leq 10r_k$,
$l \in J_{k-1}$ such that $|x_{l,k-1}-x_{i,k}| \leq 20 r_k$,
and $z\in \Sigma_{k-1}$ such that $y = \sigma_{k-1}(z)$.
Then
$$
\varepsilon_{k-1}(z) \leq 40 \varepsilon'_k(y),
\leqno (7.21)
$$
$$
{\rm Angle}(T\Sigma_{k}(y),P_{i,k}) %v
\leq C \varepsilon'_k(y),
\leqno (7.22)
$$
and
$$
\big| D\pi_{i,k} \cdot (\pi_{j,k}(y)-y) \big|
\leq C \varepsilon'_k(y)^2 r_k
\hbox{ for all $j \in J_k$ such that } 
y \in 10B_{j,k}. %v
\leqno (7.23)
$$

\ms
Recall that we can find $i$, $l$, and $z$ as above, by 
the definition (2.2) of $V^8_k$, by (2.3),
because $k \geq 1$, and by definition of $\Sigma_k$.
Notice that 
$$
|z-y| = |z-\sigma_{k-1}(z)| \leq C \varepsilon r_k
\ \hbox{ and } z \in \Sigma_{k-1} \cap 4B_{l,k-1}
\leqno (7.24)
$$ 
by (6.8) and because $|z-x_{l,k-1}| \leq |z-y| + |y-x_{i,k}|
+ |x_{i,k}-x_{l,k-1}| \leq C \varepsilon r_k + 30r_k$.

Let $m \in J_{k-1}$ be such that $z \in 10 B_{m,k-1}$. 
Then $y\in 11B_{m,k-1}$ and we claim that
$$
d_{x_{m,k-1},200r_{k-1}}(P_{i,k},P_{m,k-1})
\leq 10 d_{x_{m,k-1},100r_{k-1}}(P_{i,k},P_{m,k-1})
\leq 10 \varepsilon'_k(y).
\leqno (7.25)
$$ %vv  
% I changed more. But feel free to return to your version,
% which I leave under %
% In fact  let $d_{x_{m,k-1},100r_{k-1}}(P_{i,k},P_{m,k-1})=\alpha$ 
% then for $x_{m,k-1}$ there is $y_{m, k-1}\in P_{i,k}$ such that 
% $|x_{m, k-1}-y_{m,k-1}|\le 100 \alpha r_{k-1}$.  
% For $x\in P_{m,k-1}\cap B(x_{m,k-1}, 200 r_{k-1})\backslash  B(x_{m,k-1}, 100 r_{k-1})$ 
% let $w=x_{m,k-1} + 100 r_{k-1}|x-x_{m,k-1}|^{-1} (x-x_{m,k-1})
% \in P_{m,k-1}\cap B(x_{m, k-1}, 100 r_{k-1})$. 
% There exists $w'\in P_{i,k}$ so that $|w-w'|\le 100r_{k-1}\alpha$. 
% A straightforward computation shows that 
% $\tilde w=y_{m,k-1} + |x-x_{m,k-1}||w'-y_{m,k-1}|^{-1} (w'-y_{m,k-1})\in P_{i,k}$ 
% satisfies $|\tilde w- x|\le 2000\alpha r_{k-1}$.
% Thus $\sup\{\dist\,(x,P_{i,k}):x\in P_{m,k-1}\cap B(x_{m,k-1}, 200r_{k-1})\}\le 10\alpha$. 
% A similar computation shows that the same inequality holds when the roles 
% of $P_{i,k}$ and $P_{m, k-1}$ are reversed which ensures by (1.7) 
% that the first inequality in (7.25) holds. 
% The second inequality in (7.25) comes from the definition in (7.18).

Set $\alpha = d_{x_{m,k-1},100r_{k-1}}(P_{i,k},P_{m,k-1})$.
Since $x_{m,k-1}\in P_{m,k-1}$ we can find $y_{m}\in P_{i,k}$ 
such that $|x_{m, k-1}-y_{m}|\leq 100 \alpha r_{k-1}$.  
For $x\in P_{m,k-1}\cap B(x_{m,k-1}, 200 r_{k-1})\sm  B(x_{m,k-1}, 
100 r_{k-1})$, set $w=x_{m,k-1} + 100 r_{k-1}|x-x_{m,k-1}|^{-1} 
(x-x_{m,k-1})$. Then $w\in P_{m,k-1}\cap B(x_{m, k-1}, 100 r_{k-1})$
and there exists $w'\in P_{i,k}$ so that $|w-w'|\le 100r_{k-1}\alpha$. 
A straightforward computation shows that 
$\widetilde w=y_{m} + |x-x_{m,k-1}||w'-y_{m}|^{-1} (w'-y_{m})\in P_{i,k}$ 
satisfies $|\tilde w- x|\leq 2000\alpha r_{k-1}$.
Thus $\,\sup\big\{\dist\,(x,P_{i,k}) \, ; \, 
x\in P_{m,k-1}\cap B(x_{m,k-1}, 200r_{k-1})\big\}\leq 2000\alpha r_{k-1}$. 
A similar computation shows that the same inequality holds when the roles 
of $P_{i,k}$ and $P_{m, k-1}$ are reversed, which ensures by (1.7) 
that the first inequality in (7.25) holds. 
The second inequality in (7.25) comes from the definition in (7.18).

If $m' \in J_{k-1}$ is another index such that 
$z \in 10 B_{m',k-1}$,
$$\eqalign{
d_{x_{m,k-1},100r_{k-1}}(P_{m',k-1},P_{m,k-1})
&\leq 2 d_{x_{m,k-1},200r_{k-1}}(P_{i,k},P_{m,k-1})
\cr& \hskip 0.4cm
+ 2 d_{x_{m',k-1},200r_{k-1}}(P_{i,k},P_{m',k-1})
\leq 40 \varepsilon'_k(y)
}\leqno (7.26)
$$
by the triangle inequality, because 
$B(x_{m,k-1},150r_{k-1}) \i B(x_{m',k-1},200r_{k-1})$,
and by (7.25) and its analogue for $m'$. This proves (7.21)
(compare with (7.8)).

Next recall that $z \in \Sigma_{k-1} \cap 4B_{l,k-1}$ (by (7.24)), 
so we may apply (7.12) and get that 
$$
{\rm Angle}(T\Sigma_{k}(y),P_{l,k-1}) 
\leq C \varepsilon_{k-1}(z) \leq C \varepsilon'_k(y)
\leqno (7.27)
$$
by (7.21). Now (7.22) follows because 
$$
{\rm Angle}(P_{l,k-1},P_{i,k}) 
\leq C d_{x_{l,k-1},100r_k}(P_{l,k-1},P_{i,k}) 
\leq C \varepsilon'_k(y)
\leqno (7.28)
$$
by (7.18) and because $y \in 10B_{i,k} \cap 3B_{l,k-1}$.

Let us now prove (7.23). By (7.24), we can apply (7.10) to 
$z \in \Sigma_{k-1} \cap 4B_{l,k-1}$ and get that 
$$
|y-\pi_{l,k-1}(z)| = |\sigma_{k-1}(z)-\pi_{l,k-1}(z)| 
\leq C \varepsilon_{k-1}(z) \, r_k
\leq C \varepsilon'_{k}(y) \, r_k
\leqno (7.29)
$$
(by (7.10) and (7.21)) and then, if $j \in J_k$ is such that $y \in 10B_{j,k}$
(as in the assumption),
$$\eqalign{
|\pi_{j,k}(y)-y| &= \dist(y,P_{j,k}) \leq \dist(y,P_{l,k-1})
+ 100r_{k-1} d_{x_{l,k-1},100 r_{k-1}}(P_{l,k-1},P_{j,k})
\cr&
\leq |y-\pi_{l,k-1}(z)| + C \varepsilon'_k(y) \, r_k
\leq C \varepsilon'_k(y) \, r_k
}\leqno (7.30)
$$
because $y \in 10B_{j,k} \cap 10 B_{l,k-1}$ and by (7.18).
But $\pi_{j,k}(y)-y$ is orthogonal to $P_{j,k}$, hence nearly 
orthogonal to $P_{i,k}$, so that
$$
\big| D\pi_{i,k} \cdot (\pi_{j,k}(y)-y) \big|
\leq C \varepsilon_k(y) \, |\pi_{j,k}(y)-y|
\leq C \varepsilon_k(y) \varepsilon'_k(y) \, r_k
\leq C \varepsilon'_k(y)^2 r_k
\leqno (7.31)
$$
by (7.8), (7.30), and (7.19); (7.23) and Lemma 7.20 follow.
\qed

\ms \proclaim Lemma 7.32.
For $k \geq 1$ and $y\in \Sigma_k \cap V^8_k$, choose $i\in J_k$
such that $|y-x_{i,k}| \leq 10r_k$. Then
$$
|D\pi_{i,k} \circ D\sigma_{k}(y) \circ D\pi_{i,k} -D\pi_{i,k}| 
\leq C \varepsilon'_k(y)^2,
\leqno (7.33)
$$
and 
$$
\big| |D\sigma_k(y)\cdot v| - 1 \big| \leq C \varepsilon'_k(y)^2 
\ \hbox{ for every unit vector } v \in T\Sigma_k(y).
\leqno (7.34)
$$

\ms
Recall that $T\Sigma_k(y)$ is the vector space parallel
to the tangent plane to $\Sigma_k$ at $y$; thus (7.34)
is a distortion estimate for the restriction of $D\sigma_k$
to $T\Sigma_k(y)$, which will be useful when we need to estimate
$|f(x)-f(y)|$.

Let us first check that for $y \in \R^n$, $k \geq 0$ and 
$j,k \in J_k$ such that $y \in 10 B_{i,k} \cap 10 B_{j,k}$,
$$
\big| D\pi_{i,k} \circ [D\pi_{j,k}-D\pi_{i,k}] \circ D\pi_{i,k} \big| 
\leq C \varepsilon_k(y)^2.  
\leqno (7.35)
$$

Let $u \in \R^n$ be given, and set $v = D\pi_{i,k}\cdot u$
and $w = D\pi_{j,k}\cdot v \,$; then $w-v \in P_{j,k}^\perp \,$, so
it is almost orthogonal to $P_{i,k}$ and
$$\eqalign{
|D\pi_{i,k} \circ &[D\pi_{j,k}-D\pi_{i,k}] \circ D\pi_{i,k}\cdot u|
= |D\pi_{i,k} \circ [D\pi_{j,k}-D\pi_{i,k}]\cdot v|
\cr&= |D\pi_{i,k} \cdot (w-v)| 
\leq C \varepsilon_k(y) |w-v|
\leq C \varepsilon_k(y)^2 |v| \leq C \varepsilon_k(y)^2|u|
}\leqno (7.36)
$$
by (7.8); (7.35) follows. 

Now let $y$ and $i\in J_k$ be as in Lemma 7.32. By (7.7),
$$
D\pi_{i,k} \circ D\sigma_{k}(y) \circ D\pi_{i,k}-D\pi_{i,k}
= A_1 + A_2,
\leqno (7.37)
$$
with
$$
A_1 = -D\pi_{i,k} + \sum_{j \in J_k}\theta_{j,k}(y) 
D\pi_{i,k} \circ D\pi_{j,k}\circ D\pi_{i,k}
\leqno (7.38)
$$
and
$$
A_2 = \sum_{j \in J_k} \big[ D\pi_{i,k} \cdot \pi_{j,k}(y) \big] 
\, \big[ D\theta_{j,k}(y) \circ D\pi_{i,k} \big].
\leqno (7.39)
$$
Recall that $y \in V_k^{8}$, so $\sum_{j \in J_k}\theta_{j,k}(y) = 1$
by (3.12) and (3.13), and 
$$\eqalign{
|A_1| &= \big| \sum_{j \in J_k}\theta_{j,k}(y) 
\big[ D\pi_{i,k} \circ D\pi_{j,k}\circ D\pi_{i,k} - D\pi_{i,k} \big]\big|
\cr&
= \big| \sum_{j \in J_k}\theta_{j,k}(y) 
D\pi_{i,k} \circ [D\pi_{j,k} - D\pi_{i,k}] \circ D\pi_{i,k} \big|
\cr&
\leq \sum_{j \in J_k}\theta_{j,k}(y) 
\big| D\pi_{i,k} \circ [D\pi_{j,k} - D\pi_{i,k}] \circ D\pi_{i,k} \big|
\leq C \varepsilon_k(y)^2
}\leqno (7.40)
$$
by (7.35) and because $y \in 10 B_{i,k} \cap 10 B_{j,k}$
for all $j\in J_k$ such that $\theta_{j,k}(y) \neq 0$.
Similarly, $\sum_{j \in J_k} D\theta_{j,k}(y) = 0$ because 
$y \in V_k^{8}$, so
$$\eqalign{
|A_2| &= \Big|\sum_{j \in J_k} 
\big[ D\pi_{i,k} \cdot (\pi_{j,k}(y)-y) \big] 
\, \big[ D\theta_{j,k}(y) \circ D\pi_{i,k} \big] \Big|
\cr&
\leq C \varepsilon'_k(y)^2 r_k \sum_{j \in J_k} 
|D\theta_{j,k}(y)|
\leq C \varepsilon'_k(y)^2
}\leqno (7.41)
$$
by (7.23) and because $|D\theta_{j,k}(y)| \leq C r_k^{-1}$ 
by (3.10) (recall that the sum has less than $C$ terms as the number of 
balls $B_{j,k}$ such that $B_{j,k}\cap B_{i,k}\not =\emptyset$ is bounded 
by a constant that only depends on $n$).
Now (7.33) follows from (7.41), (7.40), (7.37), and (7.19).

Next let $v$ be a unit vector in $T\Sigma_k(y)$, and write
$D\sigma_k(y)\cdot v = w_1 + w_2 + w_3$,
with $w_1 = D\pi_{i,k}^\perp \circ D\sigma_k(y)\cdot v$,
$w_2 = D\pi_{i,k} \circ D\sigma_k(y) \circ D\pi_{i,k}\cdot v$,
and $w_3 = D\pi_{i,k} \circ D\sigma_k(y) \circ D\pi_{i,k}^\perp\cdot v$.
Observe that 
$$
|w_1| = \big|D\pi_{i,k}^\perp \circ 
[D\sigma_k(y)-D\pi_{i,k}]\cdot v \big|
\leq C \varepsilon_k(y)
\leqno (7.42)
$$
because $D\pi_{i,k}^\perp \circ D\pi_{i,k}=0$ and by (7.11),
$$
|w_2 - D\pi_{i,k}\cdot v| = 
\big|[D\pi_{i,k} \circ D\sigma_k(y) \circ D\pi_{i,k} - D\pi_{i,k}]
\cdot v \big| \leq C \varepsilon'_k(y)^2
\leqno (7.43)
$$
by (7.33), and
$$\eqalign{
|w_3| &=  
\big| D\pi_{i,k} \circ [D\sigma_k(y)-D\pi_{i,k}]
\circ D\pi_{i,k}^\perp\cdot v \big|
\leq C \varepsilon_k(y) \, |D\pi_{i,k}^\perp\cdot v|
\cr&
\leq C \varepsilon_k(y) \, {\rm Angle}(T\Sigma_{k}(y),P_{i,k})) 
\leq C \varepsilon_k(y) \varepsilon'_k(y)
\leq C \varepsilon'_k(y)^2
}\leqno (7.44)
$$
by (7.11), (7.22), and (7.19). Thus 
$$\eqalign{
\big| |D\sigma_k(y)\cdot v|^2 - 1 \big|
&= \big| |w_1|^2 + |w_2+w_3|^2 - 1 \big|
\leq \big| |w_2|^2 - 1 \big| + C\varepsilon'_k(y)^2
\cr&
\leq \big| |D\pi_{i,k}\cdot v|^2 - 1 \big| + C\varepsilon'_k(y)^2
= |D\pi_{i,k}^\perp \cdot v|^2 + C\varepsilon'_k(y)^2
}\leqno (7.45)
$$
because $w_1$ is orthogonal to $w_2+w_3$,
by (7.42), (7.44), (7.43), and because $1 = |v|^2 =  
|D\pi_{i,k}\cdot v|^2 + |D\pi_{i,k}^\perp\cdot v|^2$. Now 
$|D\pi_{i,k}^\perp \cdot v| \leq {\rm Angle}(T\Sigma_{k}(y),P_{i,k})) 
\leq C \varepsilon'_k(y)$ by (7.22), so 
$\big| |D\sigma_k(y)\cdot v|^2 - 1 \big| \leq C\varepsilon'_k(y)^2$,
as needed for (7.34). Lemma 7.32 follows.
\qed

\bsi
{\bf 8. H\"older and Lipschitz properties of $f$ on $\Sigma_0$}
\ms
In this section we use the distortion estimates 
from Section 7 to prove that the restriction of $f$
to $\Sigma_0$ is bi-H\"older in general, and bi-Lipschitz if 
we have a good enough control on the $\varepsilon_k(z)$.
This will be used later to control the function $g$. 

\ms\proclaim Proposition 8.1.
There is a constant $C \geq 0$ such that, with the notation of
the previous sections and if $\varepsilon$ is small enough, 
$$
(1-C\varepsilon) |x-y|^{1+C\varepsilon} \leq |f(x)-f(y)| 
\leq (1+C\varepsilon) |x-y|^{1-C\varepsilon}
\leqno (8.2)
$$
for $x,y\in \Sigma_0$ such that $|x-y| \leq 1$.

\ms
Recall from (6.9) that
$$
|f(x)-f_k(x)| \leq C \varepsilon r_k
\ \hbox{ for } x\in \Sigma_0 \hbox{ and } k \geq 0
\leqno (8.3)
$$
and in particular $|f(x)-x| \leq C \varepsilon$, so we also have 
a good control when $|x-y| \geq 1$, namely,
$\big||f(x)-f(y)| - |x-y|\big| \leq C \varepsilon$.

\ms
Now let $x$, $y \in \Sigma_0$ be given, with $0 < |x-y| \leq 1$,
and set $x_k = f_k(x)$ and $y_k = f_k(y)$ for $k \geq 0$.
Let us check that
$$
|x_m-y_m| \leq (1+C \varepsilon)^{m+1} |x-y|
\ \hbox{ for } m \geq 0.
\leqno (8.4)
$$

Choose a smooth arc $\gamma : I \to\Sigma_0$, that goes from
$x$ to $y$, and such that 
$$
{\rm length}(\gamma) \leq (1+C\varepsilon)|x-y|.
\leqno (8.5)
$$
Such a curve exists, by the description of $\Sigma_0$
as a Lipschitz graph near (2.4), and 
$$
y_m-x_m = \int_{I}D( f_m \circ \gamma)(t) dt
= \int_{I} Df_m(\gamma(t)) \cdot \gamma'(t) dt
\leqno (8.6)
$$
by the fundamental theorem of calculus. Fix
$t\in I$, and set $z_k = f_k(\gamma(t))$.
Also set $v_k = Df_k(\gamma(t)) \cdot \gamma'(t)$
for $0 \leq k \leq m$; thus $v_0 = \gamma'(t)$,
$v_m = Df_m(\gamma(t)) \cdot \gamma'(t)$, and
by (4.1), $v_{k+1} = D\sigma_k(z_k)\cdot v_k$ for $k<m$.
By definition $v_k \in T\Sigma_{k}(z_{k})$ for 
$0 \leq k \leq m$.

For the sake of Proposition 8.1, we just need to 
know that
$$
|v_{k+1}| = |D\sigma_k(z_k)\cdot v_k|
\leq (1+C\varepsilon)\, |v_k|
\leqno (8.7)
$$
by (7.2), which implies that
$$
|Df_m(\gamma(t)) \cdot \gamma'(t)| = |v_m|
\leq (1+C \varepsilon)^{m} \, |v_0|
= (1+C \varepsilon)^{m} \, |\gamma'(t)|.
\leqno (8.8)
$$
Let us record slightly better estimates, that will be 
used later for our bi-Lipschitz results. 
When $z_k \notin V_k^{10}$, then $D\sigma_k(z_k) = I$ by (4.5), 
so $v_{k+1} = v_k$. When $z_k \in V_k^8$, Lemma 7.32 says that
$$
\big| |v_{k+1}| - |v_k| \big|
=\big| |D\sigma_k(z_k)\cdot v_k| - |v_k| \big| 
\leq C \varepsilon'_k(z_k)^2 |v_k|,
\leqno (8.9)
$$
and hence 
$$
|v_{k+1}| \leq \big[1 + C\varepsilon'_k(z_k)^2 \big] 
\, |v_k|.
\leqno (8.10)
$$ 
Now Lemma 6.5 says that if $z_k \in V_k^{10} \sm V_k^8$ 
for some $k$, then $z_l \in V_l^4 \i V_l^8$ for $l < k$, 
and $z_l \in \Sigma_l \sm V_l^{10}$ for $l > k$.
Thus we only need to use (8.7) once, and otherwise 
we can rely on (8.10) or the trivial estimate. Thus
$$
|Df_m(\gamma(t)) \cdot \gamma'(t)| = |v_m|
\leq (1+C \varepsilon) \, |\gamma'(t)| \, 
\prod_{0 \leq k < m \, ; \, z_k \in V_k^8} 
\big[ 1+ C\varepsilon'_k(z_k)^2 \big].
\leqno (8.11)
$$
We return to (8.6) and get that
$$
|x_m-y_m| \leq (1+C \varepsilon)^{m} \,\int_I |\gamma'(t)| dt
\leq (1+C \varepsilon)^{m+1} |x-y|
\leqno (8.12)
$$
by (8.8) and (8.5). This proves (8.4).

\ms
Now we want to check that
$$
|x_m-y_m| \geq (1+C \varepsilon)^{-m-1} |x-y|
\ \hbox{ for $m \geq 0$ such that}
\leqno (8.13)
$$
$$
|x_k-y_k| \leq r_k \ \hbox{ for } 0 \leq k \leq m.
\leqno (8.14)
$$
We may assume that $m \geq 1$, because (8.13)
is trivial for $|x_0-y_0| = |x-y|$. 
Recall from Lemma 6.12 that $\Sigma_m$ coincides with
a small Lipschitz graph in $B(y_m,19r_m)$; then there is
a $C^2$ curve  $\gamma : I \to \Sigma_m$, that goes 
from $x_m$ to $y_m$, such that
$$
{\rm length}(\gamma) \leq (1+C\varepsilon) \, |x_m-y_m|.
\leqno (8.15)
$$

\ms
Recall from Lemma 7.1 that each $\sigma_k : \Sigma_k \to \Sigma_{k+1}$
is a $C^2$ diffeomorphism, so we can define 
$\sigma_k^{-1} : \Sigma_{k+1} \to \Sigma_k$ and 
$f_m^{-1} : \Sigma_m \to \Sigma_{0}$. Now 
$f_m^{-1} \circ \gamma : I \to \Sigma_0$ is a path 
from $x$ to $y$, and
$$
y-x = \int_I  D(f_m^{-1} \circ \gamma)(t) dt
= \int_{I} Df_m^{-1}(\gamma(t)) \cdot \gamma'(t) dt
\leqno (8.16)
$$
by the fundamental theorem of calculus. 
Fix $t\in I$, and set $z_0 = f_m^{-1}(\gamma(t))$ and
$z_k = f_k(z_0) = f_k\circ f_m^{-1}(\gamma(t))$ for 
$1 \leq k \leq m$. Then set 
$v_0 = Df_m^{-1}(\gamma(t)) \cdot \gamma'(t)$.
Observe that 
$$
\gamma'(t) = Df_m(z_0) \cdot v_0
\leqno (8.17)
$$
because $Df_m^{-1}(\gamma(t))$ is the inverse of 
$Df_m(z_0)$. Then set $v_k = Df_k(z_0) \cdot v_0$
for $k \leq m$. In particular, $v_m = \gamma'(t)$ 
by (8.17). The chain rule says that
$$
v_{k+1} = D\sigma_k(z_k)\cdot v_k
\ \hbox{ for } 0 \leq k < m.
\leqno (8.18)
$$
Note also that $v_k \in T\Sigma_{k}(z_{k})$ by the definition.
We now argue as for the upper bound (8.4). First observe that
$$
|v_{k+1}| = |D\sigma_k(z_k)\cdot v_k|
\geq (1+C\varepsilon)^{-1} |v_k|
\leqno (8.19)
$$
by (7.2), so 
$$
|Df_m^{-1}(\gamma(t)) \cdot \gamma'(t)| = |v_0| 
\leq (1+C \varepsilon)^{m} \, |v_m|
= (1+C \varepsilon)^{m} \, |\gamma'(t)|.
\leqno (8.20)
$$
Let us also record here a better estimate.
When $z_k \notin V_k^{10}$, then $D\sigma_k(z_k) = I$ 
by (4.5), so $v_{k+1} = v_k$. When $z_k \in V_k^8$, 
Lemma 7.32 says that
$$
\big| |v_{k+1}| - |v_k| \big|
=\big| |D\sigma_k(z_k)\cdot v_k| - |v_k| \big| 
\leq C \varepsilon'_k(z_k)^2 |v_k|,
\leqno (8.21)
$$
and hence 
$$
|v_{k+1}| \geq \big[1 + C\varepsilon'_k(z_k)^2 \big]^{-1} 
\, |v_k|.
\leqno (8.22)
$$ 
As before, Lemma 6.5 says that we only need to use (8.19) 
once, and otherwise we can rely on (8.22) or the trivial 
estimate. Thus
$$\eqalign{
|Df_m^{-1}(\gamma(t)) \cdot \gamma'(t)| &= |v_0| 
\leq (1+C \varepsilon) \, |v_m| \, 
\prod_{0 \leq k < m \, ; \, z_k \in V_k^8} 
\big[ 1+ C\varepsilon'_k(z_k)^2 \big] 
\cr&
\leq (1+C \varepsilon) \, |\gamma'(t)| \, 
\prod_{0 \leq k < m \, ; \, z_k \in V_k^8} 
\big[ 1+ C\varepsilon'_k(z_k)^2 \big].
}\leqno (8.23)
$$
We return to (8.20), plug it into (8.16), integrate,
and get that
$$
|y-x| \leq \int_{I} |Df_m^{-1}(\gamma(t)) \cdot \gamma'(t)| dt
\leq (1+C \varepsilon)^{m} \, \int_{I} |\gamma'(t)| dt
\leq (1+C \varepsilon)^{m+1} |x_m-y_m|
\leqno (8.24)
$$
by (8.15). This completes our proof of (8.13).

\ms
We are now ready to prove (8.2).
Observe that
$$
10^{-m} = r_m \geq |x_m-y_m| \geq (1+C \varepsilon)^{-m-1} |x-y|
\leqno (8.25)
$$
for every $m$ such that (8.14) holds, by (8.13).
Hence (8.14) cannot hold for $m$ large, because the left-hand 
side tends to $0$ much faster than the right-hand side
(recall that we assume that $x \neq y \,$).

Since $|x-y| \leq 1$, we have (8.14) for $m=0$; hence 
there is a largest $m \geq 0$ such that (8.13) holds,
which for convenience we still call $m$. Thus
$|x_{m+1}-y_{m+1}| > r_{m+1}$, and since
$|x_{m+1}-x_{m}|+|y_{m+1}-y_{m}| \leq C \varepsilon r_m$
by (6.8), we get that
$$
|x_{m}-y_{m}| > {r_{m+1} \over 2}.
\leqno (8.26)
$$
We also have (8.25), which by taking logarithms yields
$m \ln\Big({10 \over 1+C \varepsilon}\Big)
\leq \ln\Big({1 +C \varepsilon\over |x-y|}\Big)$,
hence $\dsp m \leq {1 \over \ln(9)}
\ln\Big({2\over |x-y|}\Big)$. Now (8.4) and (8.13) yield
$$
\Big| \ln\Big({|x_m-y_m|\over |x-y|}\Big) \Big|
\leq (m+1) \ln(1+C\varepsilon)
\leq C\varepsilon \, \ln\Big({2\over |x-y|}\Big).
\leqno (8.27)
$$
At the same time, 
$$
|f(x)-x_m| + |f(y)-y_m| \leq C \varepsilon r_m
\leqno (8.28)
$$
by (8.3) and $|x_{m}-y_{m}| > r_{m+1}/2$ by (8.26), so
$$
\Big| \ln\Big({|f(x)-f(y)|\over |x_m-y_m|}\Big) \Big|
\leq C \varepsilon,
\leqno (8.29)
$$
hence $\dsp\Big| \ln\Big({|f(x)-f(y)|\over |x-y|}\Big) \Big|
\leq C\varepsilon \, \ln\Big({2\over |x-y|}\Big)$ and
$$
\Big({2\over |x-y|}\Big)^{-C\varepsilon} \leq 
{|f(x)-f(y)|\over |x-y|} \leq \Big({2\over |x-y|}\Big)^{C\varepsilon},
\leqno (8.30)
$$
which implies (8.2). Proposition 8.1 follows.
\qed

\msi{\bf Remark 8.31.}
The proof of Proposition 8.1 also yields
$$
(1-C\varepsilon) |x-y|^{1+C\varepsilon} \leq |f_k(x)-f_k(y)| 
\leq (1+C\varepsilon) |x-y|^{1-C\varepsilon}
\leqno (8.32)
$$
for all $k \geq 1$ and $x,y\in \Sigma_0$ such that $|x-y| \leq 1$.
When $k \leq m$ (the largest integer for which (8.14) holds,
we just use (8.4) and (8.13); when $k > m$, we just replace
(8.28) with 
$$
|f_k(x)-x_m| + |f_k(y)-y_m| \leq C \varepsilon r_m
\leqno (8.33)
$$
and conclude as before. Also, (8.32) is nothing more
than (8.2) applied to the mapping $\widetilde f$ that
we would get by replacing all $J_l$, $l \geq k$,
with the empty set.

The same remark applies to Proposition 8.34 below.

\ms
With Proposition 8.1, we have a proof of the weaker
variants of Theorems 1.1 and 1.10 where we only want to define $f$
on the plane $P(0,10)$. Similarly, Proposition 8.34 below will
lead to weaker variants of Theorem 1.13 and Theorem 1.18
once we sort out the relations between the $J_q$ ( see (1.16))
and the $\varepsilon'_k$. [But we shall not do this before
Sections 12 and 13.]

\ms
\proclaim Proposition 8.34. Suppose that for some $M < +\infty$
$$
\sum_{k \geq 0} \varepsilon'_k(f_k(z))^2 \leq M
\ \hbox{ for } z \in \Sigma_0,
\leqno (8.35)
$$
where the $\varepsilon'_k(f_k(z))$ are as in (2.21)-(2.22) 
or (7.18). Then $f : \Sigma_0 \to \Sigma$ is bi-Lipschitz.

\ms 
We do not try to make (8.35) look nicer for the moment; 
we shall make attempts in this direction in Sections 12 and 13. 
We keep the same proof as for Proposition 8.1 above, but
instead of (8.8) and (8.20) we use (8.11) and (8.23).
In both case the product is less than $C(M)$ by (8.35), and we get 
that 
$$
(1+C\varepsilon)^{-1} C(M)^{-1} |x-y| \leq 
|x_m-y_m| \leq (1+C\varepsilon) C(M) |x-y|
\leqno (8.36)
$$
by the proof of (8.12) and (8.24). We apply this to the same 
$m$ as above (the largest one for which (8.13) holds); then
(8.29) still holds. Thus
$$
(2C(M))^{-1} |x-y| \leq |f(x)-f(y)| \leq 2C(M) |x-y|.
\leqno (8.37)
$$
This takes care of the case when $|x-y| \leq 1$, but
the situation when $|x-y| \geq 1$ is even better, because
$|f(x)-x|+|f(y)-y| \leq C \varepsilon$ by (8.3).
Proposition 8.34 follows.
\qed

\bsi
{\bf 9. $C^2$-regularity of the $\Sigma_k$ and
fields of linear isometries defined on $\Sigma_0$}
\ms 

We first focus on the regularity of the approximating
surfaces $\Sigma_k$, and in particular the fact that the tangent
direction $T\Sigma_k(y)$ is a Lipschitz function of $y \in \Sigma_k$.

We shall measure distances between $d$-dimensional
vector subspaces $V$, $V'$ by setting
$$\eqalign{
D(V,V') &= \Max \Big\{ 
\sup\big\{ \dist(v,V') \, ; \, v \in V \cap \overline B(0,1) \big\} ;
\cr&\hskip 4cm
\sup\big\{ \dist(v',V) \, ; \, v' \in V' \cap \overline B(0,1) \big\}
\Big\}.
}\leqno (9.1)
$$
This is equivalent to the notion of angle used so far.

\ms\proclaim Lemma 9.2.
We have that for $k \geq 0$ and $x,x'\in \Sigma_0$ 
such that $|x'-x| \leq 10$,
$$
D(T\Sigma_{k+1}(f_{k+1}(x)),T\Sigma_k(f_k(x)))
\leq C_1 \varepsilon
\leqno (9.3)  
$$
and
$$
D(T\Sigma_{k}(f_{k}(x')),T\Sigma_k(f_k(x)))
\leq C_2 \varepsilon \, r_k^{-1} |f_{k}(x')-f_{k}(x)|.
\leqno (9.4)  
$$

\ms
Let us first prove (9.3). Let $k \geq 0$ and $x\in \Sigma_0$
be given. Set $y = f_k(x)$.
If $y \in \Sigma_n \sm V_k^{11}$, (6.10) says that
$\sigma_k(w) = w$ near $y$, so $T\Sigma_{k+1}(f_{k+1}(x))=
T\Sigma_k(f_k(x))$ and (9.3) is trivial.

If $y \in \Sigma_k \cap V_k^{11}$,  
choose $j \in J_k$ such that $y\in 11B_{j,k}$;
Proposition 5.4 gives a good description in 
$40 B_{j,k}$ of both $\Sigma_k$ and $\Sigma_{k+1}$, 
as a $C \varepsilon$-Lipschitz graphs over $P_{j,k}$. 
Since both $f_k(x) = y$ and $f_{k+1}(x) = \sigma_k(y)$ 
lies well inside $40 B_{j,k}$ (by (6.8)), we get (9.3).

We shall prove (9.4) by induction. When $k=0$, we need to show that 
$D(T\Sigma_{0}(x'),T\Sigma_0(x)) \break\leq C\varepsilon |x'-x|$.
This follows from the local Lipschitz graph description of $\Sigma_0$
in Section 2, and in particular (2.4).

Now suppose that $k \geq 0$, assume that (9.4) holds for $k$, and 
prove it for $k+1$. Set
$$
D = D(T\Sigma_{k+1}(f_{k+1}(x')),T\Sigma_{k+1}(f_{k+1}(x))),
\leqno (9.5) 
$$
$y = f_k(x)$, and $y' = f_k(x')$. Observe that
$$\eqalign{
D& \leq D(T\Sigma_{k}(f_{k}(x')),T\Sigma_k(f_k(x))) + 
2 D(T\Sigma_{k+1}(f_{k+1}(x')),T\Sigma_k(f_k(x')))
\cr&\hskip 5.5cm
+ 2 D(T\Sigma_{k+1}(f_{k+1}(x)),T\Sigma_k(f_k(x)))
\cr&
\leq  C_2 \varepsilon  r_k^{-1} |y'-y| + 2C_1 \varepsilon
}\leqno (9.6)  
$$
by induction assumption and (9.3). If $|y'-y| \geq r_k$,
$|\sigma_k(y')-\sigma_k(y)| \geq |y'-y| -C \varepsilon r_k
\geq |y'-y|/2$, so
$$
C_2 \varepsilon r_k^{-1} |y'-y|
\leq {C_2 \varepsilon \over 5} \, r_{k+1}^{-1} 
|\sigma_k(y')-\sigma_k(y)|
\leqno (9.7)  
$$
because $r_{k+1} = r_k/10$. At the same time,
$|\sigma_k(y')-\sigma_k(y)| \geq |y'-y|/2 \geq r_k/2 = 5r_{k+1}$, so
$$
2C_1 \varepsilon \leq {2 C_1 \varepsilon \over 5}\, r_{k+1}^{-1}
|\sigma_k(y')-\sigma_k(y)|,
\leqno (9.8)  
$$
and (9.4) for $k+1$ follows from (9.5)-(9.8) if we take $C_2 \geq C_1$.
So we may now assume that $|y'-y| < r_k$.

Let $v$ be a unit vector in $T\Sigma_{k+1}(f_{k+1}(x))$, and
use Lemma 7.1 to find $u \in T\Sigma_{k}(y)$ such that
$v = D\sigma_k(y) \cdot u$. Note that $|u| \leq 1+C\varepsilon$
by (7.2). By our induction assumption (9.4), we can find
$u'\in T\Sigma_{k}(y')$ such that
$$
|u'-u| \leq (1+C\varepsilon) C_2 \varepsilon \, r_k^{-1} |y'-y|,
\leqno (9.9)  
$$
and of course $v' = D\sigma_k(y) \cdot u'$
lies in $T\Sigma_{k+1}(f_{k+1}(x'))$, so
$$\eqalign{
\dist(v,T&\Sigma_{k+1}(f_{k+1}(x'))) \leq |v'-v|
= |D\sigma_k(y') \cdot u' - D\sigma_k(y) \cdot u|
\cr&
\leq |D\sigma_k(y') \cdot (u'-u)|+|[D\sigma_k(y')-D\sigma_k(y)] \cdot u|
}\leqno (9.10)  
$$
Let us check that
$$
|D\sigma_k(y')| \leq 1+C\varepsilon.
\leqno (9.11)  
$$
When $y'\in \Sigma_k \sm V_k^{10}$, $D\sigma_k(y')=I$ by 
(4.5). Otherwise, (5.13) says that
$|D\sigma_k(y')-D\pi_{j,k}-\psi_k(y')D\pi_{j,k}^\perp| \leq C 
\varepsilon$ for some $j \in J_k$; (9.11) follows because
$|D\pi_{j,k}+\psi_k(y')D\pi_{j,k}^\perp| \leq 1$. Now
(9.11) and (9.9) yield
$$
|D\sigma_k(y') \cdot (u'-u)|
\leq (1+C\varepsilon)^2 C_2 \varepsilon \, r_k^{-1} |y'-y|
\leq 2 C_2 \varepsilon \, r_k^{-1} |y'-y|.
\leqno (9.12)  
$$

We also want to estimate 
$[D\sigma_k(y')-D\sigma_k(y)] \cdot u$. When we differentiate
the first part of (4.2), we get that
$$
D\sigma_k(y) = I + \sum_{j\in J_k} D\theta_{j,k}(y) [\pi_{j,k}(y)-y]
+ \sum_{j\in J_k} \theta_{j,k}(y) [D\pi_{j,k}-I]. %v
\leqno (9.13)  
$$
Thus $D\sigma_k(y')-D\sigma_k(y) = A+B+D$, where
$$
A = \sum_{j\in J_k} [D\theta_{j,k}(y')-D\theta_{j,k}(y)] 
\,[\pi_{j,k}(y')-y'],
\leqno (9.14)  
$$
$$
B = \sum_{j\in J_k}D\theta_{j,k}(y) \,[\pi_{j,k}(y')-y'-\pi_{j,k}(y)+y]
\leqno (9.15)  
$$
and
$$
D = \sum_{j\in J_k} [\theta_{j,k}(y')-\theta_{j,k}(y)]\, [D\pi_{j,k}-I].
\leqno (9.16)  
$$

Let $j \in J_k$ be such that 
$D\theta_{j,k}(y')-D\theta_{j,k}(y) \neq 0$;
then $y\in 10B_{j,k}$ or $y'\in 10B_{j,k}$, and in both
cases $y' \in \Sigma_k \cap 11B_{j,k}$ because we now 
assume that $|y'-y| \leq r_k$. By (5.5)-(5.7) in Proposition~5.4,
$|\pi_{j,k}(y')-y'| = \dist(y',P_{j,k}) \leq C \varepsilon r_k$.
Thus
$$
|A| \leq C r_k^{-2} |y'-y| |\pi_{j,k}(y')-y'| 
\leq C \varepsilon r_k^{-1} |y'-y|
\leqno (9.17)  
$$
by (3.10).

Next let $j$ be such that $D\theta_{j,k}(y) \neq 0$;
as before, $y\in 10B_{j,k}$ and hence $y'\in 11B_{j,k}$.
By Proposition 5.4, both points lie on a $C\varepsilon$-Lipschitz graph
over $P_{j,k}$, so
$$
|\pi_{j,k}(y')-y'-\pi_{j,k}(y)+y| \leq C \varepsilon |y'-y|
\leqno (9.18)  
$$
and hence
$$
B \leq C \varepsilon r_k^{-1} |y'-y|
\leqno (9.19)  
$$
by (3.10). 

For $D$, observe again that $y\in 11B_{j,k}$ when
$\theta_{j,k}(y')-\theta_{j,k}(y) \neq 0$, then
(5.6) and (5.7) in Proposition~5.4 say that
$|D\pi_{j,k}\cdot u-u| \leq C \varepsilon |u|$ for
$u \in T\Sigma_k(y)$, so that
$$
|D \cdot u| \leq \sum_{j\in J_k} |\theta_{j,k}(y')-\theta_{j,k}(y)|\, 
\, |D\pi_{j,k}\cdot u-u|
\leq C \varepsilon r_k^{-1} |y'-y| |u|
\leqno (9.20)  
$$
by (3.10) again. Finally,
$$
|[D\sigma_k(y')-D\sigma_k(y)] \cdot u|
\leq |A||u|+|B||u| + |D\cdot u| \leq C \varepsilon r_k^{-1} |y'-y|
\leqno (9.21)  
$$
because $|u| \leq 1+C\varepsilon$, and 
$$
\dist(v,T\Sigma_{k+1}(f_{k+1}(x'))) \leq 
(2C_2+C) \varepsilon r_k^{-1} |y'-y|
< 3C_2 \varepsilon r_k^{-1} |y'-y|
\leqno (9.22)
$$
by (9.10), (9.12), (9.21), and if $C_2$ is large enough.
Of course the proof of (9.22) also yields
$\dist(v',T\Sigma_{k+1}(f_{k+1}(x))) \leq 
3C_2 \varepsilon r_{k}^{-1} |y'-y|$ for any unit vector
$v' \in T\Sigma_{k+1}(f_{k+1}(x'))$, so 
$$
D(T\Sigma_{k+1}(f_{k+1}(x')),T\Sigma_{k+1}(f_{k+1}(x)))
\leq 3C_2 \varepsilon \, r_k^{-1} |y'-y|.
\leqno (9.23)  
$$
Recall from (7.4) that
$$\eqalign{
|y'-y| &\leq |\sigma_k(y')-\sigma_k(y)|
+ |\sigma_k(y)-\sigma_k(y') - y+y'|
\cr&
\leq |\sigma_k(y')-\sigma_k(y)|+ C \varepsilon |y-y'|,
}\leqno (9.24)
$$ 
and hence
$$\eqalign{
|y'-y| &\leq 2|\sigma_k(y')-\sigma_k(y)|
= 2 |f_{k+1}(x') - f_{k+1}(x)|.
}\leqno (9.25)
$$ 
Now (9.4) for $k+1$ follows from (9.23) and (9.25),
because $r_{k} = 10 r_{k+1}$.
This completes our proof of Lemma 9.2.
\qed

\ms
Let us encode Lemma 9.2 in terms of orthogonal projections.
For $k \geq 0$ and $y\in \Sigma_k$, denote by $\pi_{y,k}$ the 
orthogonal projection onto $T\Sigma_k(y)$, and by 
$\pi_{y,k}^\perp = I -\pi_{y,k}$ the orthogonal projection 
onto the orthogonal subspace.
Note that $\pi_{y,k}$ and $\pi_{y,k}^\perp$ are $C^1$
functions of $y$, because $\Sigma_k$ is $C^2$ (by Lemma 6.12),
and we could compute $\pi_{y,k}$ locally, with a fixed
orthogonal basis of $\R^n$, a description of $\Sigma_k$ 
in coordinates, and a Gram-Schmidt orthogonalization process
to compute the projections.
In addition, since $\pi_{y,k}$ is a Lipschitz function of 
$T\Sigma_k(y)$, we get the following immediate consequence of 
Lemma~9.2.

\ms\proclaim Lemma 9.26.
For $k \geq 0$ and $x,x'\in \Sigma_0$ 
such that $|x'-x| \leq 10$,
$$
|\pi_{f_{k+1}(x),k+1}-\pi_{f_k(x),k}|
\leq C \varepsilon
\leqno (9.27)  
$$
and
$$
|\pi_{f_{k}(x'),k}-\pi_{f_k(x),k}|
\leq C \varepsilon \, r_k^{-1} |f_{k}(x')-f_{k}(x)|.
\leqno (9.28)  
$$

\ms
Now we want to construct fields of linear isometries defined 
on $\Sigma_0$. This corresponds to the construction of coherent
orthonormal bases in [Mo] and [To], %% Duke et Morrey
except that since we do not want to assume that $\Sigma_0$
is orientable, for instance, we cannot start our construction 
with the knowledge of a smooth choice of orthonormal basis
for $T\Sigma_0(x)$. So instead we shall define smooth mappings
$R_k$ on $\Sigma_0$, with values in the set of linear isometries 
of $\R^n$, and such that for $x\in \Sigma_0$ and $k \geq 0$,
$R_k(x)(T\Sigma_0(x))=T\Sigma_k(f_k(x))$.

\ms\proclaim Proposition 9.29.
Let $\cal R$ denote the set of linear isometries of $\R^n$.
Also set
$$
T_k(x) = T\Sigma_k(f_k(x))
\ \hbox{ for } x \in \Sigma_0
\hbox{ and } k \geq 0.
\leqno (9.30)
$$
There exist $C^1$ mappings $R_k : \Sigma_0 \to {\cal R}$, with the
following properties:
$$
R_0(x) = I \hbox{ for } x \in \Sigma_0,
\leqno (9.31)
$$
$$
R_k(x)(T_0(x))=T_k(x)
\ \hbox{ for } x \in \Sigma_0
\hbox{ and } k \geq 0,
\leqno (9.32)
$$
$$
|R_{k+1}(x)-R_{k}(x)| \leq C \varepsilon
\ \hbox{ for } x \in \Sigma_0
\hbox{ and } k \geq 0,
\leqno (9.33)
$$
and, if we set 
$$
\widetilde R_k(y) = R_k \circ f_k^{-1}(y)
\hbox{ for } y \in \Sigma_k
\leqno (9.34)
$$ 
and denote by $D_y\widetilde R_k$ the differential 
of $\widetilde R_k(y)$ with respect to 
$y\in \Sigma_k$,
$$
|D_y\widetilde R_{k+1}(y)| 
\leq C_1 r_{k+1}^{-1} \varepsilon
\ \hbox{ for $k \geq 0$ and } y \in \Sigma_{k+1}.
\leqno (9.35)
$$

\ms
Note that $f_k^{-1} : \Sigma_k \to \Sigma_0$  is well defined and $C^2$, 
by Lemma 7.1. We again give a special name to $C_1$
in (9.35) to clarify role of constants in the proof by induction.

As strongly suggested by (9.31), we take $R_0(x) = I$ for $x \in \Sigma_0$.
Now let $k \geq 0$ be given, assume that we already constructed 
$R_k: \Sigma_0 \to {\cal R}$ with the desired properties, and 
construct $R_{k+1}$.
As a first attempt, we set $y = f_{k+1}(x)$ for $x\in \Sigma_0$ and try
$$
S_k(x) = \pi_{y,k+1} \circ R_k(x) \circ \pi_{x,0}
+\pi_{y,k+1}^\perp \circ R_k(x) \circ \pi_{x,0}^\perp.
\leqno (9.36)
$$
This is probably not an isometry, but at least
$$
S_k(x)(T_0(x)) \i T_{k+1}(x) \hbox{ and }
S_k(x)(T_0(x)^\perp) \i T_{k+1}(x)^\perp.
\leqno (9.37)
$$
Next, $S_k(x)$ is a $C^1$ function of $x$,
because $R_k(x)$ is $C^1$ (by induction assumption),
the $f_k$ are $C^1$ on $\Sigma_0$, $\Sigma_0$ is $C^2$, and 
$\pi_{y,k+1}$ is a $C^1$ function of $y \in \Sigma_{k+1}$.
Set
$$
\widetilde S_{k}(y) = S_k \circ f_{k+1}^{-1}(y)
\hbox{ for } y \in \Sigma_{k+1}.
\leqno (9.38)
$$
$\widetilde S_k$ is $C^1$ because $S_k$ is $C^1$ and $f_{k+1}$ is a $C^2$ 
diffeomorphism from $\Sigma_0$ to $\Sigma_{k+1}$ (see Lemma~7.1).
Since $R_k(x) = \widetilde R_{k}(f_k(x))
= \widetilde R_{k}(\sigma_k^{-1}(y))$ when $y = f_{k+1}(x)$,
we get that
$$
\widetilde S_{k}(y) 
= \pi_{y,k+1} \circ \widetilde R_k(\sigma_k^{-1}(y)) \circ \pi_{x,0}
+\pi_{y,k+1}^\perp \circ \widetilde R_k(\sigma_k^{-1}(y)) 
\circ \pi_{x,0}^\perp.
\leqno (9.39)
$$
We want to differentiate this with respect to $y\in \Sigma_{k+1}$.
First denote by $A_1$ the part of $D_y \widetilde S_k(y)$ that comes
from differentiating $\pi_{y,k+1}$ and $\pi_{y,k+1}^\perp=I-\pi_{y,k+1}$.
Note that $|D_y\pi_{y,k+1}| \leq C \varepsilon r_{k+1}^{-1}$
by (9.28), so $A_1 \leq 2C \varepsilon r_{k+1}^{-1}$ (because 
$|\widetilde R_k(\sigma_k^{-1}(y))| = 1$).

For the part of $D_y \widetilde S_k(y)$ that comes
from differentiating $\widetilde R_k(\sigma_k^{-1}(y))$,
we use the fact that $D\sigma_k^{-1}(y)$ is 
$(1+C\varepsilon)$-bi-Lipschitz on the tangent plane,
by Lemma 7.1. If $k \geq 1$, we get that 
$$
|A_2| \leq (1+C\varepsilon) |D\widetilde R_k|
\leq C_1 (1+C\varepsilon) r_k^{-1} \varepsilon
\leq 2 C_1 r_k^{-1} \varepsilon
\leqno (9.40)
$$
by induction assumption. When $k=0$, $D\widetilde R_k = 0$
by (9.31), so $A_2 = 0$.

Finally let $A_3$ be the part of $D_y \widetilde S_k(y)$ that comes
from differentiating $\pi_{x,0}$ and $\pi_{x,0}^\perp$. We know that 
$|D_x\pi_{x,0}| \leq C \varepsilon$ by (9.28) (or directly (2.4)),
and the differential of $x = f_{k+1}^{-1}(y)$ with respect to
$y$ has a norm less than $(1+C\varepsilon)^{k+1}$, by (7.2) and the
chain rule. So
$$
|A_3| \leq C (1+C\varepsilon)^{k+1} \varepsilon \leq C r_k^{-1} \varepsilon
\leqno (9.41)
$$
because $r_k = 10^{-k}$. Altogether,
$$
|D_y \widetilde S_k(y)| \leq |A_1|+|A_2|+|A_3|
\leq (2C_1 + C) r_k^{-1} \varepsilon,
\leqno (9.42)
$$
which gives some hope for (9.35). 

Observe that since $R_k$ maps $T_0(x)$ to $T_k(x)$
(by (9.32)) and hence $T_0(x)^\perp$ to $T_k(x)^\perp$
(because it is a linear isometry), 
$R_k(x) = \pi_{f_k(x),k} \circ R_k(x) \circ \pi_{x,0}
+ \pi_{f_k(x),k}^\perp \circ R_k(x) \circ \pi_{x,0}^\perp$ 
and hence
$$
|R_k(x) - S_k(x)| \leq |\pi_{f_k(x),k}-\pi_{y,k+1}| 
+ |\pi_{f_k(x),k}^\perp-\pi_{y,k+1}^\perp|
= 2 |\pi_{f_k(x),k}-\pi_{y,k+1}|
\leq C \varepsilon
\leqno (9.43)
$$
by (9.36), because $|R_k(x)| \leq 1$, by (9.27), and because
$y = f_{k+1}(x)$. As a consequence, $S_k(x)$ is nearly an isometry,
and even
$$
S_k(x) \in U = \big\{ S \in {\cal L}(\R^n, \R^n) \, ; \, 
|S S^\ast - I| \leq \eta \big\},
\leqno (9.44)
$$
where the small $\eta\geq \varepsilon^{1/4}$ will be chosen soon. 
%vv  again, I will modify a lot here and keep your version below
For convenience we denote by $S^\ast$ the transpose of $S$.
We set 
$$
R_{k+1}(x) = H(S_k(x)),
\leqno (9.45)
$$
where $H : U \to {\cal R}$ is a nonlinear projection on the 
set of linear isometries that we define by
$$
H(S) = (S S^\ast)^{-1/2} S
\ \hbox{ for } S \in U.
\leqno (9.46)
$$
Here the simplest way to define $(S S^\ast)^{-1/2}$ is 
to take $(S S^\ast)^{-1/2} = \sum_{n \geq 0} a_n (S S^\ast-I)^n$,
where $\sum_{n \geq 0} a_n x^n$ is the expansion
of $(1+x)^{-1/2}$ near $0$. The series converges 
as soon as $|S S^\ast-I| < 1$, which is the case for $S \in U$. 
% Since $U$ is a compact thin tubular neighborhood of 
% ${\cal R}$ and $H$ is $C^2$ in $U$, to estimate the Lipschitz constant
% of $H$ in $U$ it is enough to estimate its derivative on ${\cal R}$ and 
% note that the second derivative of $H$ in $U$ is bounded. 
% If $T=R(I+A)$ with $R\in {\cal R}$ and $|A|$ small then 
% $H(T)= R+ R{{A-A^\ast}\over{2}} + {R\over 4}(A^2+ (A^\ast)^2 + A^\ast A + 3 AA^\ast) +RG(A)$
% where $|G(A)|\le C|A|^3$ for $|A|$ small and 
% $|G(A_1)-G(A_2)|\le C(|A_1|^2+|A_2|^2)|A_1-A_2|$. 
% Thus the first derivative of $H$ at $R$ has norm 1, 
% and the second derivative is bounded in $U$. 
% For $S\in U$ there exists $R\in{\cal R}$ so that $|S-R|\le C\eta$ and 
% $||DH(R)|-|DH(S)||\le |DH(R)-DH(S)|\le \sup_{U}|D^2H||R-S|\le C\eta$. 
% Hence $|DH(S)|\le 1 +C\eta\le 1+ 10^{-2}$ provided $\eta$ is small enough, and
% $$
% \hbox{$H$ is $(1+10^{-2})$-Lipschitz on $U$.}
% \leqno (9.47)
% $$
% As promised, $H(S) \in \cal R$ for $S \in U$, because
% $H(S)H(S)^\ast = (S S^\ast)^{-1/2} S S^\ast (S S^\ast)^{-1/2}
% = I$ because $(S S^\ast)^{-1/2}$ commutes with $SS^\ast$ and 
% its square is $(S S^\ast)^{-1}$ (say that we could manipulate 
% power series). Also observe that
% $$
% H(S) = S \hbox{ for } S \in {\cal R},
% \leqno (9.48)
% $$
% just because $S S^\ast = I$.
The use of $H$ is our substitute for the 
Gram-Schmidt orthogonalization process used by Morrey [Mo] %%
and Toro [To] to define fields of orthonormal bases. %%
As promised, 
$$
H(S) \in \cal R \ \hbox{ for } S \in U, 
\leqno (9.47)
$$
since
$H(S)H(S)^\ast = (S S^\ast)^{-1/2} S S^\ast (S S^\ast)^{-1/2}
= I$ because $(S S^\ast)^{-1/2}$ commutes with $SS^\ast$ and 
its square is $(S S^\ast)^{-1}$ (say that we could manipulate 
power series). Also,
$$
H(S) = S \ \hbox{ for } S \in {\cal R},
\leqno (9.48)
$$
just because $S S^\ast = I$.
Next we want to show that if $\eta$ is small enough, 
$$
\hbox{$H$ is $(1+10^{-2})$-Lipschitz on $U$.}
\leqno (9.49)
$$
First observe that $|S| = |SS^\ast|^{1/2} \leq (1+\eta)^{1/2}$
for $S\in U$; then 
$$
\dist(S,{\cal R}) \leq |S-H(S)|
= |S|\,|(S S^\ast)^{-1/2}-I| \leq \eta
\ \hbox{ for } S \in U
\leqno (9.50)
$$
(use (9.48) and the power series expansion). 
If  $S, S' \in U$ are such that $|S'-S| \geq 200\eta$,
then
$$\eqalign{
|H(S)-H(S')| &\leq |H(S)-S|+|S-S'|+|S'-H(S')|
\cr&\leq |S-S'| + 2 \eta \leq (1+ 10^{-2}) |S-S'|,
}\leqno (9.51)
$$
by (9.50) and as needed, so for the proof of (9.49) we may assume 
that $S$ and $S'$ lie in a same ball $B$ of radius $202\eta$ 
centered on $\cal R$.
Note that $H$ is defined on $B$, and its second derivative on $B$ 
is bounded by $100$, trivially by the power series expansion of 
$(SS^\ast)^{-1/2}$. Denote by $D_H(S)$ the differential of
$H$ at $S \in B$; we shall check soon that
$$
||D_H(R)|| \leq 1 
\ \hbox{ for } R \in {\cal R}
\leqno (9.52)
$$
and this immediately implies that $||D_H(S)|| \leq 1 + C \eta$
for $S \in B$, which will complete the proof of (9.49) because 
$B$ is convex.

We still need to check (9.52). Let $R\in {\cal R}$ be given, and
let us check that
$$
D_H(R) \cdot A = {1 \over 2} (A-RA^\ast R)
\ \hbox{ for } A \in {\cal L}(\R^n).
\leqno (9.53)
$$
Let us expand $H(S)$ when $S = R+A$, with $A$ small: observe that
$$
SS^\ast - I = RR^\ast + RA^\ast + AR^\ast + AA^\ast -I
= RA^\ast + AR^\ast + O(|A|^2), 
\leqno (9.54)
$$
then
$(SS^\ast)^{-1/2} = I -{1 \over 2}(RA^\ast + AR^\ast) +  O(|A|^2)$
and 
$$\eqalign{
H(S)& = (SS^\ast)^{-1/2}(R+A)
= R + A -{1 \over 2}(RA^\ast R + AR^\ast R) +  O(|A|^2)
\cr&= R + {1 \over 2}(A - RA^\ast R) +  O(|A|^2).
}\leqno (9.55)
$$
So (9.53), (9.52) follows easily, and we can choose $\eta$
so that (9.49) holds.
%v  Et maintenant il faut augmenter les numeros

From (9.45), (9.48), and (9.43) we deduce that
$$\eqalign{
|R_{k+1}(x)-R_k(x)| &= |H(S_k(x)) - R_k(x)|
=|H(S_k(x)) - H(R_k(x))| 
\cr&\leq (1+10^{-2}) |S_k(x)-R_k(x)|
\leq C \varepsilon,
}\leqno (9.56)
$$
so (9.33) holds. Also, 
$$
\widetilde R_{k+1}(y) = R_{k+1}(f_{k+1}^{-1}(y))
= H(S_k(f_{k+1}^{-1}(y)))
= H(\widetilde S_{k}(y))
\leqno (9.57)
$$
by the definition (9.34) and (9.38), so the chain rule
gives
$$
|D_y\widetilde R_{k+1}(y)| \leq (1+10^{-2}) \, |D_y \widetilde S_k(y)| 
\leq (3C_1 + 2C) r_k^{-1} \varepsilon < C_1 r_{k+1}^{-1}
\leqno (9.58)
$$
by (9.49) and (9.42), if $C_1$ is large enough, and because
$r_{k} = 10 r_{k+1}$. This proves (9.35).

We still need to prove (9.32), and since we do not
 understand square 
roots of operators, we shall take orthonormal bases.
Denote by $\pi$ the orthogonal projection onto $T_k(x)$, and set
$A = \pi_{y,k+1} \circ \pi \circ \pi_{y,k+1}$. This is a self-adjoint
operator on $\R^n$, and it maps $T_{k+1}(x)$ to itself (recall that
$\pi_{y,k+1}$ is the orthogonal projection onto $T_{k+1}(x)$), so
its restriction to $T_{k+1}(x)$ is self-adjoint. Thus there 
is an orthonormal basis $e_1, \ldots, e_d$ of $T_{k+1}(x)$ such that
$A(e_l) = \lambda_l e_l$ for $1 \leq l \leq d$ and some real
numbers $\lambda_l$. Note also that $A$ vanishes on $T_{k+1}(x)^\perp$.

Similarly, $A' = \pi_{y,k+1}^\perp \circ \pi^\perp \circ \pi_{y,k+1}^\perp$ 
is self-adjoint, so there is an orthonormal basis $e_{d+1}, \ldots, e_n$ 
of $T_{k+1}(x)^\perp$ such that $A'(e_l) = \lambda_l e_l$ for 
$d+1 \leq l \leq n$. Then the matrix of $A+A'$ in the basis
$e_1,\ldots, e_n$ is diagonal, with entries $\lambda_l$.

Observe that $R_k(x)$ sends $T_0(x)$ to $T_k(x)$, by the
induction assumption (9.32), and sends $T_0(x)^\perp$ to $T_k(x)^\perp$,
because it is an isometry. Hence %v
$R_k(x) \circ \pi_{x,0} =
\pi \circ R_k(x)$ (recall that $\pi$ is the orthogonal projection
on $T_k(x)$) and similarly
$R_k(x) \circ \pi_{x,0}^\perp = \pi^\perp \circ R_k(x)$. Then
$$\eqalign{
S_k(x) &=  \pi_{y,k+1} \circ R_k(x) \circ \pi_{x,0}
+\pi_{y,k+1}^\perp \circ R_k(x) \circ \pi_{x,0}^\perp
\cr&
= \pi_{y,k+1} \circ \pi \circ R_k(x) 
+\pi_{y,k+1}^\perp \circ \pi^\perp \circ R_k(x)
\cr&
= (\pi_{y,k+1} \circ \pi +  \pi_{y,k+1}^\perp \circ \pi^\perp) 
\circ R_k(x)
}\leqno (9.59) 
$$
by (9.36). Next
$$\leqalignno{
S_k(x)S_k(x)^\ast &=
(\pi_{y,k+1} \circ \pi +  \pi_{y,k+1}^\perp \circ \pi^\perp) 
\circ R_k(x) \circ R_k(x)^\ast \circ 
(\pi\circ \pi_{y,k+1} +  \pi^\perp\circ \pi_{y,k+1}^\perp) 
\cr& = (\pi_{y,k+1} \circ \pi +  \pi_{y,k+1}^\perp \circ \pi^\perp) 
\circ (\pi\circ \pi_{y,k+1} +  \pi^\perp\circ \pi_{y,k+1}^\perp) 
\cr&
= \pi_{y,k+1} \circ \pi \circ \pi_{y,k+1}
+ \pi_{y,k+1}^\perp \circ \pi^\perp \circ \pi_{y,k+1}^\perp
= A+A'
& (9.60) 
}%\leqno (9.60) 
$$
so the matrix of $S_k(x)S_k(x)^\ast$ in the basis $e_1,\ldots, e_n$ 
is diagonal with entries $\lambda_1, \ldots, \lambda_n$. %v
Thus $R_{k+1}(x) = H(S_k(x)) = D S_k(x)$, where $D$ is diagonal
with entries $\lambda_l^{-1/2}$. In particular, %v
$D$ preserves the spaces $T_{k+1}(x)$ and $T_{k+1}(x)^\perp$, and
$$
R_{k+1}(x)(T_0(x)) = D[S_k(x)(T_0(x)]
\i D[T_{k+1}(x)] \i T_{k+1}(x)
\leqno (9.61) 
$$
by (9.37). The inclusion is an identity because 
$R_{k+1}$ is an isometry, and this proves (9.32);
Proposition 9.29 follows.
\qed

\msi{\bf Remark 9.62.}
It would not be too difficult to prove that with our easy-to-get
additional regularity assumption (2.6) on $\Sigma_0$, the surfaces 
$\Sigma_k$ are of class $C^{m_0}$, with bounds like 
$$
|D^mA| \leq C_m M_m 2^m + C_m \varepsilon r_k^{1-m}
\leqno (9.63) 
$$ 
in the small Lipschitz representation of Lemma 6.12,
or similar estimates for the $A_{j,k}$ and $F_{j,k}$
of Proposition 5.4.

With more work, we could try to improve the estimates
on the the restriction of $D\sigma_k(y)$ to $T\Sigma_k(y)$,
get better estimates on $D_y \widetilde R_{k+1}$
in Proposition 9.29, or improve (9.63), for instance when 
the square summability condition (8.35) holds.
We do not do these computations.

\bsi
{\bf 10. The definition of $g$ on the whole $\R^n$}
\ms

We continue with the notations and assumptions of the previous
sections.
We shall soon be ready to define the mapping $g$ promised
in the various statements of Sections 1 and 2. We still need a nearest point
projection on $\Sigma_0$, defined in a tubular neighborhood
of $\Sigma_0$. If $\Sigma_0$ is a plane, the mapping defined
in the next lemma is simply the orthogonal projection onto $\Sigma_0$.

\ms\proclaim Lemma 10.1.
Set $V = \big\{ z \in \R^n \, ; \, \dist(z,\Sigma_0) < 40 \big\}$.
For each $z\in V$, there is a unique point $p(z)\in \Sigma_0$
such that $|p(z)-z| \leq 50$ and $p(z)-z$ is orthogonal to
$T\Sigma_0(p(z))$. 
In addition, the mapping $p : V \to \Sigma_0$ is of class 
$C^1$, and 
$$
|p(z')-p(z)| \leq (1+C\varepsilon) |z'-z|
\ \hbox{ for } z,z' \in V \hbox{ such that } |z'-z|\leq 1.
\leqno (10.2)
$$
Similarly, if we set $q(z) = z - p(z)$ for $z\in V$,
$$
|q(z')-q(z)| \leq (1+C\varepsilon) |z'-z|
\ \hbox{ for } z,z' \in V \hbox{ such that } |z'-z|\leq 1.
\leqno (10.3)
$$

\ms
We shall first define $p$ locally. Let $x \in \Sigma_0$ be given, and
let $P_x$ and $F_x$ be as in the local Lipschitz description of 
$\Sigma_0$ near (2.4) and (2.5). Denote by $\Gamma$ the graph
of $F_x$ over $P_x$ and, for $x'\in \Gamma$, denote by
$T(x')$ the vector space parallel to the tangent plane 
to $\Gamma$ at $x'$. 
Denote by $\pi_T(x')$ the orthogonal projection onto 
$T(x')$, and set $\pi_T^\perp(x') = I - \pi_T(x')$. 
Set $H = P_x \times (P_x^\perp \cap B(0,100))$ and define a function 
$\Phi : H \to \R^n$ by
$$
\Phi(u,v) = u + F_x(u) + \pi_T^\perp(u + F_x(u)) \cdot v
\leqno (10.4)
$$
for $(u,v)\in H$. Recall from (2.4) that $T(u + F_x(u))$ makes 
a small angle with $P_x$, so 
$$
|\pi_T^\perp(u + F_x(u))-\pi_x^\perp| \leq 2 \varepsilon, 
\leqno (10.5)
$$
where $\pi_x^\perp$ still denotes the orthogonal 
projection onto $P_x^\perp$. In addition, the fact that
$|D^2F_x| \leq \varepsilon$ on $P_x$ (by (2.4)) implies that
$$
|D_u \pi_T^\perp(u + F_x(u))| \leq C \varepsilon,
\leqno (10.6)
$$
where we denote by $D_u$ the differential with respect to $u \in P_x$.
Indeed, we could compute $\pi_T^\perp(u + F_x(u))$
from $DF_x(u)$ by a painful but explicit Gram-Schmidt 
orthogonalization process.

The mapping $\Phi$ is of class $C^1$ (more if we assume (2.6)), 
and its differential is given by
$$
D_u\Phi(u,v) = I_u + DF_x(u) + D_u \pi_T^\perp(u + F_x(u))\cdot v
\ \hbox{ and } \ 
D_v\Phi(u,v) = \pi_T^\perp(u + F_x(u)) \circ I_v
\leqno (10.7)
$$
(where $I_u$ and $I_v$ simply denote the canonical injections 
from the planes parallel to $P_x$ and $P_x^\perp$ into $\R^n$). Thus
$$
|D\Phi(u,v) - I| \leq C \varepsilon, 
\leqno (10.8)
$$
by (2.4), (10.5), and (10.6). Because of the simple shape of the
domain of definition of $\Phi$, we deduce from (10.8) and the 
fundamental theorem of calculus that $\Phi$ is a $C^1$ 
diffeomorphism from $H$ to $\Phi(H)$. Since
$$
|\Phi(u,v)-(u+v)| \leq C \varepsilon
\leqno (10.9)
$$ 
by (2.4) and (10.5), $\Phi(H)$ contains 
$P_x + (P_x^\perp \cap B(0,99)) = \big\{ z\in \R^n \, ; \, 
\pi_x^\perp(z) \in B(0,99) \big\}$.

\ms
We now return to the lemma itself. Let $z \in B(x,45)$ be given.
Then we can find $(u,v) \in H$ such that $z = \Phi(u,v)$.
Set $p(z) = u + F_x(u)$. By (10.9), $|u-\pi_x(z)| 
= |\pi_x(u+v-z)| =|\pi_x(u+v-\Phi(u,v))|\leq C \varepsilon$
and similarly $|v-\pi_x^\perp(z)| \leq C \varepsilon$.
Thus $u \in B(x,46)$, $p(z) = u + F_x(u) \in B(x,47)$,
and hence $p(z) \in \Sigma_0$ by (2.5). Also,
$$\eqalign{
|z-p(z)| &\leq |z-u| + |F_x(u)| \leq |z-u| + \varepsilon 
\leq |z-\pi_x(z)| + C \varepsilon 
\cr&
= |\pi_x^\perp(z)| + C\varepsilon
= |\pi_x^\perp(z-x)| + C\varepsilon
\leq |z-x| + C\varepsilon < 46
}\leqno (10.10)
$$
by (2.4) and because $\pi_x^\perp(x)=0$. Finally,
$z-p(z) = \Phi(u,v) - p(z) = \Phi(u,v) - u - F_x(u) =
\pi_T^\perp(u + F_x(u)) \cdot v = \pi_T^\perp(p(z)) \cdot v$ 
is orthogonal to $T\Sigma_0(p(z))$, as required for the lemma. 

Next we check the uniqueness. %v
Let $\xi\in \Sigma_0 \cap B(z,50)$ be such that 
$z-\xi \in T\Sigma_0(\xi)^\perp$. By (2.5) and because 
$z \in B(x,45)$, $\xi \in \Gamma$, so $\xi = u + F_x(u)$ for some 
$u\in P_x$. Also, $z-\xi =\pi_T^\perp(\xi)\cdot v$ %v
for some $v \in P_x^\perp$ (recall from (10.5) that 
$\pi_T^\perp(\xi) : P_x^\perp \to T\Sigma_0(\xi)^\perp$ is bijective). 
In addition, $|v|\leq (1+C\varepsilon)|z-\xi|\leq (1+C\varepsilon) 50$
so $(u,v) \in H$. But then $z = \Phi(u,v)$ and since $\Phi$ is 
injective, $(u,v)$ is the same as in the definition of $p(z)$. Thus
$\xi = u + F_x(u) = p(z)$, as needed for the uniqueness.

This gives a good definition of $p(z)$ for $z\in B(x,45)$.
There $p$ and $q=I-p$ are differentiable, and
$$
Dp(z) = [I + DF_x(u)] \circ D\pi_x \circ D\Phi^{-1}(z),
\leqno (10.11)
$$
with $u = \pi_x(\Phi^{-1}(z))$. Thus 
$|Dp(z)-D\pi_x| \leq C \varepsilon$ by (10.8).

We get the desired function $p$ on $V$ by covering $V$ by balls
$B(x,45)$, $x\in \Sigma_0$. There is no difficulty with potentially
different definitions in overlapping domains, because of the uniqueness. 
Finally, for (10.2) and (10.3), we can choose $x\in \Sigma_0$ such that
$|x-z| \leq 40$, and then compute in $B(z,2) \i B(x,45)$ and
integrate $Dp$ and $Dq=I-Dp$.
This completes the proof of Lemma 10.1.
\qed

\ms
Let us also check that 
$$
|q(z)| = |z-p(z)| = \dist(z,\Sigma_0)
\ \hbox{ for } z\in V.
\leqno (10.12)
$$
Let $p\in \Sigma_0$ be such that $|z-p| = \dist(z,\Sigma_0)$;
such a $p$ exists by compactness, and $p \in \overline B(z,40)$
by definition of $V$. The differential of $|z-y|^2$
along $\Sigma_0$ vanishes at $y=p$ (because $|z-p|$ is minimal),
so $z-p$ is orthogonal to $T\Sigma_0(p)$. By Lemma 10.1, 
$p=p(z)$ and (10.12) holds.

\ms
We are now ready to define $g$. We set
$$
g(z) = z \ \hbox{ for } z \in \R^n \sm V,
\leqno (10.13)
$$ 
$$
g(z) =  f(z) \ \hbox{ for } z \in \Sigma_0,
\leqno (10.14)
$$
and 
$$
g(z) =  \sum_{k \geq 0} \rho_k(q(z)) \,
\Big\{ f_k(p(z)) + R_k(p(z)) \cdot q(z) \Big\}
\ \hbox{ for } z \in V \sm \Sigma_0,
\leqno (10.15)
$$
where $f_k$ is as in Section 4, $R_k$ was constructed in Section 9,
and the cut-off functions $\rho_k$ are defined below.
For a given $z$, the sum in (10.15) will have at most three
terms, by (10.18) below.

Choose $h : \R_+ \to [0,1]$ smooth, nondecreasing, and such that
$h(t) = 1$ for $t \geq 2$ and $h(t) = 0$ for $t \leq 1$. Then
set 
$$
\rho_0(y) = h(|y|)
\ \hbox{ and }
\rho_k(y) = h(r_k^{-1}|y|) - h(r_{k-1}^{-1}|y|)\ge 0\ 
\hbox{ for } k \geq 1.
\leqno (10.16)
$$
Notice that
$$
\sum_{k\ge 0} \rho_k(y) = 1 \ \hbox{ for } y \in \R^n \sm \{ 0 \},
\leqno (10.17)
$$
$\rho_0(y) =0$ when $|y| \leq 1$, and, for $k \geq 1$,
$$
\rho_k(y) = 0 
\hbox{ unless } r_k < |y| < 20 r_k.
\leqno (10.18)
$$

\ms
In what follows, it will be convenient to set systematically 
$$
x = p(z) \in \Sigma_0 \ \hbox{ and }
y = q(z) = z-p(z) \in T\Sigma_0(x)^\perp
\leqno (10.19)
$$
for $z \in V$. With these notations, (10.15) becomes 
the nicer-looking
$$
g(z) =  \sum_{k \geq 0} \rho_k(y) 
\big\{ f_k(x) + R_k(x) \cdot y \big\}
\ \hbox{ for } z = x+y \in V \sm \Sigma_0.
\leqno (10.20)
$$

\ms \proclaim Lemma 10.21.
The function $g$ is continuous on $\R^n$,
$$
g(z) = z 
\ \hbox{ on } \big\{ z\in V \, ; \, 
\dist(z,\Sigma_0) \geq 2 \big\},
\leqno (10.22)
$$
and
$$
|g(z)-z| \leq C \varepsilon
\ \hbox{ for } z \in \R^n.
\leqno (10.23)
$$

\ms
We first check (10.22). If $z\in V$ is such that
$\dist(z,\Sigma_0) \geq 2$, then $|y| = |z-p(z)| \geq 2$, 
and so $\rho_0(y) = 1$ and $\rho_k(y) = 0$ for $k \geq 1$. 
Hence $g(z) = f_0(x) + R_0(x) \cdot y
= x + y = z$ by (4.1) and (9.31). 

The continuity of $g$ across $\partial V$ follows from
(10.22) and (10.13). Let us now check the continuity across
$\Sigma_0$. For $z\in V \sm \Sigma_0$ (and with $x=p(z)$),
$$\eqalign{
|g(z)-f(x)| &\leq \sum_{k \geq 0} \rho_k(y) 
\big| f_k(x)-f(x) + R_k(x) \cdot y \big|
\cr&\leq \sum_{k \geq 0} \rho_k(y) 
\big\{ |f_k(x)-f(x)|+|y| \big\}
}\leqno (10.24)
$$
by (10.20), (10.17), and because $R_k(x)$ is an isometry.
In addition, $r_k < |y|$ when $\rho_k(y) \neq 0$, and then
$|f_k(x)-f(x)| \leq C \varepsilon r_k \leq C \varepsilon |y|$
by (6.9), so 
$$
|g(z)-f(x)| \leq (1+C\varepsilon) |y|
= (1+C\varepsilon) \dist(z,\Sigma_0),
\leqno (10.25)
$$
by (10.12). Then $g$ is continuous across $\Sigma_0$,
because $f$ is continuous on $\Sigma_0$.

Finally we check (10.23). By (10.13) and (10.22), we may
assume that $z\in V$ and even $\dist(z,\Sigma_0) \leq 2$.
Then
$$\eqalign{
|g(z)-z| &\leq \sum_{k \geq 0} \rho_k(y) 
\big| f_k(x)-x + (R_k(x)-I) \cdot y \big|
\cr&\leq \sum_{k \geq 0} \rho_k(y) 
\big\{ |f_k(x)-x|+|R_k(x)-I| \, |y| \big\}
\cr&
\leq C \varepsilon + |y| \sum_{k \geq 0} \rho_k(y) |R_k(x)-I|
\leq C \varepsilon + C \varepsilon |y| \sum_{k \geq 0} k \rho_k(y)
}\leqno (10.26)
$$
because $z=x+y$, and by (10.20), (10.17), (6.9),
and (9.33).
In addition, when $\rho_k(y) \neq 0$, (10.18) says that $r_k<|y| < 20 
r_k$, so $\sum_{k \geq 0} k \rho_k(y) \leq C \ln(20/|y|)$, 
and (10.23) follows from (10.26).
\qed

\ms
It will be useful in Section 14 to know that
$$
(1-C\varepsilon)\dist(z,\Sigma_0) \leq \dist(g(z),\Sigma) 
\leq (1+C\varepsilon) \dist(z,\Sigma_0)
\ \hbox{ for } z\in \R^n.
\leqno (10.27)
$$
This is clear when $z\in \Sigma_0$, because then 
$g(z)=f(z) \in \Sigma$, so we may assume that 
$z\in \R^n \sm \Sigma_0$.
The second inequality follows from (10.25). The first one 
holds trivially when $\dist(z,\Sigma_0) \geq 1$, because
$\Sigma$ is $C\varepsilon$-close to $\Sigma_0$ by (6.9),
and $|g(z)-z| \leq C \varepsilon$ by (10.23).
So we may assume that $z \in V$ and $\dist(z,\Sigma_0) \leq 1$. 
Let $m$ be the smallest integer such that $\rho_m(y) \neq 0$.
Thus $m > 0$ because $|y| = \dist(z,\Sigma_0) \leq 1$ 
(by (10.12)), and so $r_m < |y| < 20r_m$ by (10.18).

Apply Lemma 6.12 to $\xi =f_{m-1}(x)$ and the integer $m-1$;
there is a plane $P$ through $\xi$ such that $\Sigma_{m-1}$ 
coincides in $B(\xi,190r_m)$ with a $C\varepsilon$-Lipschitz 
graph over $P$. Hence
$$
{\rm Angle}(P,T_{m-1}(x)) =
{\rm Angle}(P,T\Sigma_{m-1}(\xi))
\leq C \varepsilon
\leqno (10.28)
$$
(recall that $T_{m-1}(x) = T\Sigma_{m-1}(f_{m-1}(x))$
by (9.30)), and also
$$
\dist(w,P) \leq C \varepsilon r_m
\ \hbox{ for } w \in \Sigma_{m-1}\cap B(\xi,190r_m)
\leqno (10.29)
$$
(because $P$ and the Lipschitz graph go through $\xi$).

By (10.18), $\rho_k(y) = 0$ unless $k = m, m+1$, or $m+2$.
For such $k$, $|f_k(x)-\xi| = |f_k(x)-f_{m-1}(x)| 
\leq C \varepsilon r_m$ by (6.8), and 
$|R_k(x)-R_{m-1}(x)| \leq C \varepsilon$ by (9.33).
Thus (10.20) yields
$$\eqalign{
|g(z) - \xi - R_{m-1}(x)\cdot y| 
&= \Big| \sum_{k \geq 0} \rho_k(y) 
\big\{ f_k(x)-\xi + [R_k(x)- R_{m-1}(x)] \cdot y \big\}\Big|
\cr&\leq C \varepsilon r_m
+ C \varepsilon |y| \leq C' \varepsilon |y|
}\leqno (10.30)
$$
because $\sum_{k} \rho_k(y) = 1$ and $r_m < |y|$. 
Set $w = \xi + R_{m-1}(x)\cdot y$; thus 
$|g(z)-w| \leq C \varepsilon |y|$. Note that
$|w-\xi| = |y| \leq 20r_m$, so
$$\eqalign{
\dist(w,\Sigma_{m-1}) &= 
\dist(w,\Sigma_{m-1}\cap B(\xi,50r_m)) 
\geq \dist(w,P) - C \varepsilon r_m
%\cr&
%\geq (1-C\varepsilon)|y| - C \varepsilon r_m
%\geq (1-C'\varepsilon)|y|
}\leqno (10.31)
$$
because $\xi\in \Sigma_{m-1}$ and $|w-\xi| \leq 20r_m$,
and by (10.29). %v
Now $w-\xi = R_{m-1}(x)\cdot y$ is 
orthogonal to $T_{m-1}(x)$, by (9.32), because $R_{m-1}(x)$
is an isometry, and because $y=q(z)$ is orthogonal to 
$T_0(x) = T\Sigma_0(x)$. By (10.28), $w-\xi$
is also nearly orthogonal to $P$, and so
$$
\dist(w,P) \geq (1-C \varepsilon) \, |w-\xi|
= (1-C \varepsilon) \, |y|. 
\leqno (10.32)
$$
Thus
$$\eqalign{
\dist(g(z),\Sigma) 
&\geq \dist(g(z),\Sigma_{m-1}) - C \varepsilon r_{m}
\geq \dist(w,\Sigma_{m-1}) - C \varepsilon (r_m +|y|)
\cr&\geq \dist(w,P) - C \varepsilon (r_m +|y|)
\geq (1-C \varepsilon) \, |y|
= (1-C \varepsilon) \, \dist(z,\Sigma_0)
}\leqno (10.33)
$$
because $\Sigma$ is very close to $\Sigma_{m-1}$
(by (6.9)), by (10.30), (10.31), and (10.32), %v
because $r_m < |y|$, and by (10.12).
This is the remaining inequality in (10.27).

\msi{ \bf Remark 10.34.}
We decided to use the coordinates $p(z) \in \Sigma_0$
and $q(z) = z-p(z)$ to define $g$ from the restriction of
the $f_k$ to $\Sigma_0$; another option would
have been to construct $g$, as we did for $f$, as a limit
of compositions $g_k$, where $g_0 = I$ and 
$$
g_{k+1} = s_k \circ g_k
\ \hbox{ for } k \geq 0.
\leqno (10.35)
$$
This is the scheme that was followed in [DDT], for instance. %%
We want $s_k$ to coincide with $\sigma_k$ on $\Sigma_k$, so as to 
obtain $g_k = f_k$ and $g=f$ on $\Sigma_0$, and the simplest seems to
use the Lipschitz graph description of $\Sigma_k$ that we get
from Proposition 5.4. That is, if we were just to define $s_k(z)$
near some $B_{j,k}$, $j\in J_k$, a first attempt would be to 
use 
$$
X_{j,k}(z) = \pi_{j,k}(z) + A_{j,k}(\pi_{j,k}(z)) \in \Sigma_k
\leqno (10.36)
$$ 
(a vertical projection of $z$ on $\Sigma_k$, constructed
with the Lipschitz function $A_{j,k}$ of Proposition~5.4)
and take $s_k(z) = \sigma_k(X_{j,k}(z)) + z-X_{j,k}(z)$.

This would not be very efficient, because for the bi-Lipschitz results
we want $\sigma_k$ to be as close as possible to an isometry 
(we want to compose lots of different $\sigma_k$), in particular in 
the regions where the $P_{j,k}$ vary very slowly. With the formula above,
if the $P_{j,k}$ turn of about $\alpha$ near $B_{j,k}$, we can expect
$s_k$ to be $C \alpha$-close to an isometry, and we would prefer
$C \alpha^2$, so that we can sum the distortions as in Proposition 8.34.
For this, a better attempt would be to try
$$
\eta_{j,k}(z) = \sigma_k(X_{j,k}(z)) + (I-R_{j,k}(X_{j,k}(z)))\cdot 
(z-X_{j,k}(z)),
\leqno (10.37)
$$
where the role of the small perturbation $R_{j,k}$ is to
correct some linear terms in the expansion of the derivative
$D\eta_k(z)$, to make it closer to an isometry. 
We cannot arrange this precisely everywhere, so we focus on $\Sigma_k$, 
because this is the place where we may need to iterate many mappings
(far from $\Sigma_k$, we shall take $s_k(z) = z$).
This is why we like to evaluate things on $\Sigma_k$, and
hope that the estimates will not
  deteriorate too fast when we leave
$\Sigma_k$.

Computations (that would need to be checked) seem to lead to the
choice of
$$
R_{j,k}(x) = 
D\pi_{j,k} \circ DA_{j,k}^\ast(\pi_{j,k}(x)) \circ D\pi_{j,k}^\perp
- D\sigma_k(x)^\ast \circ D\pi_{j,k}^\perp
\leqno (10.38)
$$
for $x\in \Sigma_k \cap 20B_{j,k}$, and where $DA_{j,k}^\ast$
and $D\sigma_k(x)^\ast$ denote adjoints of linear mappings.
%The precomposition with $D\pi_{j,k}^\perp$ is not shocking, 
%because $DA_{j,k}$ really maps to $P_{j,k}^\perp$. 
Now we suggest to take
$$
s_k(z) = \psi_k(z) \, z + \sum_{j \in J_k} \theta_{j,k}(z) \,\eta_{j,k}(z)
= z + \sum_{j \in J_k} \theta_{j,k}(z) [\eta_{j,k}(z)-z]
\leqno (10.39)
$$
for $z \in \R^n$ and a little like in (4.2). 
Note that $s_k(z) = \sigma_k(z)$ for $z\in \Sigma_k$, because 
then $X_{j,k}(z) = z$, so we will get that $g=f$ on $\Sigma_0$.
Now one should estimate the $Ds_k$ on $\R^n$ as we estimated
the $D\sigma_k$ on $\Sigma_k$, and mimic the proof of Section 8;
this involves slightly ugly computations and in particular we need better
estimates on the second derivatives $D^2\sigma_k$ and $D^{2}A_{j,k}$,
but at the end we seem to get a bi-Lipschitz mapping $g$ when (8.35) 
holds.

Altogether using (10.39) seems to require more computations
(even if we save the construction of isometry fields),
but we mention it because it could be useful in some contexts 
(suppose we want to cut $g$ into small bi-Lipschitz pieces).

\bsi
{\bf 11. H\"older and Lipschitz properties of $g$ on $\R^n$}
\ms

In this section we complete the proof of Theorems 2.15 and 2.23,
and in particular establish the desired bi-H\"older and 
bi-Lipschitz estimates for $g$.

\ms\proclaim Proposition 11.1.
There is a constant $C \geq 0$ such that, with the notation of
the previous sections and if $\varepsilon$ is small enough, 
$$
{1 \over 4} |z'-z|^{1+C\varepsilon} \leq |g(z')-g(z)| 
\leq 3 |z'-z|^{1-C\varepsilon}  
\leqno (11.2)
$$
for $z,z'\in \R^n$ such that $|z'-z| \leq 1$.

\ms
We shall make no attempt here to check that the constants
${1 \over 4}$ and $3$ in (11.2) could be replaced 
$1 \pm C\varepsilon$, even though we would be ready to claim
that this is possible (but by being more meticulous about 
orthogonality in the decompositions).

Since by (10.13) and (10.22) $g(z)=z$ out of $\big\{ z\in V \, ; \, 
\dist(z,\Sigma_0) \geq 2 \big\}$, it is enough to prove
(11.2) when $z$ and $z'$ lie in $\big\{ z\in V \, ; \, 
\dist(z,\Sigma_0) \leq 3 \big\}$. Set
$$
x = p(z), y = q(z), x'=p(x'), \hbox{ and } y'=q(z')
\leqno (11.3)
$$
as above. We may as well assume that 
$z' \neq z$, that $|y'| \leq |y|$ (by symmetry), and
that $y,y' \neq 0$ (we can always let $z$ and $z'$
tend to $\Sigma_0$ once we have (11.2) away from $\Sigma_0$).

Let us first settle the case when $|f(x')-f(x)| \geq 10|y|$.
In this case 
$$
|g(z)-f(x)| + |g(z')-f(x')| \leq (1+C \varepsilon) (|y|+|y'|)
\leq {2 (1+C \varepsilon) \over 10} \, |f(x')-f(x)|
\leqno (11.4)
$$
by (10.25). The second inequality in (11.2) holds because
$$\eqalign{
|g(z')-g(z)| &
\leq |f(x')-f(x)| + |g(z)-f(x)| + |g(z')-f(x')|
\cr&
\leq {13 \over 10} \, |f(x')-f(x)|
\leq {13 \over 10} \, (1+C \varepsilon) |x'-x|^{1-C\varepsilon}
\cr&\leq {13 \over 10} \, (1+C \varepsilon)^2 |z'-z|^{1-C\varepsilon}
\leq 2|z'-z|^{1-C\varepsilon}
}\leqno (11.5)
$$
by (8.2) and because $p$ is locally $(1+C\varepsilon)$-Lipschitz
(by Lemma 10.1). Similarly
$$\eqalign{
|g(z')-&g(z)| 
 \geq |f(x')-f(x)| - |g(z)-f(x)| - |g(z')-f(x')|
\cr&\geq {8-C\varepsilon \over 10}\, |f(x')-f(x)|
\geq {8 \over 10}\,(1-C\varepsilon) |x'-x|^{1+C\varepsilon}
\geq {7 \over 10}\, |x'-x|^{1+C\varepsilon}
}\leqno (11.6)
$$
by (8.2), and at the same time
$$
|g(z')-g(z)| \geq {8-C\varepsilon \over 10}\, |f(x')-f(x)|
\geq (8-C\varepsilon) |y| \geq 7 |y| \geq {7 \over 2}\,|y-y'|
\leqno (11.7)
$$
because $|y'|\leq |y|$. Also observe that
$|z'-z| \leq |x'-x|+|y-y'|$ brutally because $z=x+y$ and $z'=x'+y'$,
so 
$$
|z'-z| \leq {2 \over 7} \, |g(z')-g(z)|
+ \Big({10 \over 7}\Big)^{1/(1+C\varepsilon)}
|g(z')-g(z)|^{1/(1+C\varepsilon)}
\leqno (11.8)
$$
by (11.7) and (11.6).
The first half of (11.2) is trivial if $|g(z')-g(z)| \geq 1$
(because $|z'-z| \leq 1$); otherwise
$|g(z')-g(z)| \leq |g(z')-g(z)|^{1/(1+C\varepsilon)}$ and
(11.8) says that
$$
|z'-z| \leq \Big[{2 \over 7} + 
\Big({10 \over 7}\Big)^{1/(1+C\varepsilon)}\Big]
|g(z')-g(z)|^{1/(1+C\varepsilon)},
\leqno (11.9)
$$
which implies the first half of (11.2).

\ms
So we may assume that $|f(x')-f(x)| < 10|y|$. 
Let $m \geq 0$ denote the smallest integer such that
$\rho_m(y) \neq 0$.%vv
If $m> 0$, (10.18) says that $r_m < |y| < 20r_m$;
otherwise $|y| > 1$ (because $\rho_0(y) = 0$ when 
$|y| \leq 1$), and $r_m < |y| < 20r_m$ as well.
% because $\rho_0(y) = 0$ when $|y| \leq 1$, 
% and then (10.18) says that $r_m < |y| < 20r_m$.
Since $|f_m(x)-f(x)| + |f_m(x')-f(x')| \leq C \varepsilon r_m$ by
(6.9), we also get that
$$
|f_m(x')-f_m(x)| < 11|y| \leq 220 r_m
\leqno (11.10)
$$
We want to estimate $g(z)-g(z')=\Delta_1+\Delta_2+\Delta_3+\Delta_4$,
where by (10.20)
$$
\Delta_1 = \sum_{k \geq 0} \rho_k(y) [f_k(x)-f_k(x')],
\leqno (11.11)
$$
$$
\Delta_2 = \sum_{k \geq 0} \rho_k(y) R_k(x)\cdot(y-y'),
\leqno (11.12)
$$
$$
\Delta_3 = \sum_{k \geq 0} \rho_k(y) [R_k(x)-R_k(x')]\cdot y',
\leqno (11.13)
$$
and 
$$
\Delta_4 = \sum_{k \geq 0} [\rho_k(y)-\rho_k(y')] 
\big\{ f_k(x') + R_k(x')\cdot y'\big\}.
\leqno (11.14)
$$
Let us replace some sums with single terms. First check that
$$
|\Delta_1 - [f_m(x)-f_m(x')]| \leq C \varepsilon |f_m(x)-f_m(x')|.
\leqno (11.15)
$$
By (10.18), the only integers $k$ such that $\rho_k(y) \neq 0$
can only be $m$, $m+1$, and $m+2$. Since 
$\sum_{k \geq 0} \rho_k(y)=1$, we just need to check
that
$$
|[f_k(x)-f_k(x')] - [f_m(x)-f_m(x')]| \leq C \varepsilon |f_m(x)-f_m(x')|.
\leqno (11.16)
$$
for these $k$'s. Set $\xi = f_m(x) \in \Sigma_m$ and $\xi' = f_m(x')$;
then (7.4), applied to $\xi$ and $\xi'$, says that
$$
|\sigma_m(\xi)-\sigma_m(\xi')-\xi+\xi'| 
\leq C \varepsilon |\xi-\xi'|,
\leqno (11.17)
$$
which is (11.16) for $m+1$. Applying again (7.4) to
$\sigma_m(\xi)$ and $\sigma_m(\xi') \in \Sigma_{m+1}$, yields
$$\eqalign{
|\sigma_{m+1}(\sigma_m(\xi))-\sigma_{m+1}(\sigma_m(\xi'))
-\sigma_m(\xi)+\sigma_m(\xi')| 
&\leq C \varepsilon |\sigma_m(\xi)-\sigma_m(\xi')|
\cr&\leq 2 C \varepsilon |\xi-\xi'|,
}\leqno (11.18)
$$
where the second inequality comes from (11.17). We add (11.17)
and (11.18) and get that
$$
|\sigma_{m+1}(\sigma_m(\xi))-\sigma_{m+1}(\sigma_m(\xi'))-\xi+\xi'| 
\leq 3C \varepsilon |\xi-\xi'|,
\leqno (11.19)
$$
which is (11.16) for $m+2$. So (11.16) and (11.15) hold.
Similarly,
$$
|\Delta_2 - R_m(x)\cdot(y-y')| 
= \big|\sum_k \rho_k(y) [R_k(x)-R_m(x)]\cdot (y-y')\big|
\leq C \varepsilon |y-y'|
\leqno (11.20)
$$
because $\sum_{k \geq 0} \rho_k(y)=1$ and
$|R_k - R_m| \leq C \varepsilon$ for $k = m+1$ or $m+2$,
by (9.33). Next,
$$
|\Delta_3| \leq |y'| \sum_k \rho_k(y) |R_k(x)-R_k(x')|.
% = |y'|\sum_k\rho_k(y)|\widetilde R_k(f_k(x))-\widetilde R_k(f_k(x'))|
\leqno (11.21)
$$
We need to estimate $|R_k(x)-R_k(x')|$, but first
let us check that for $m \leq k \leq m+2$
there is a path $\gamma$ in $\Sigma_k$, that goes from 
$f_k(x)$ to $f_k(x')$, and such that
$$
{\rm length}(\gamma) \leq 2 |f_k(x)-f_k(x')|.
\leqno (11.22)
$$
Indeed, $|f_k(x)-f_k(x')| \leq |f_m(x)-f_m(x')|+C\varepsilon r_m \leq 221 r_m$ 
by (6.8) and (11.10), and we know from Lemma 6.12 that $\Sigma_k$ 
coincides with a small Lipschitz graph in every ball of radius 
$19 r_k$ centered on $\Sigma_k$. If $|f_k(x)-f_k(x')| \leq 18 r_k$,
we get $\gamma$ immediately from Lemma~6.12; otherwise, we will
 first
need to connect $f_k(x)$ to $f_k(x')$ by a small chain of points.
This will be easy, but let us do the argument anyway.

First suppose that $m \geq 2$. 
Lemma 6.12 gives a good Lipschitz control of $\Sigma_{m-2}$
in $B(f_{m-2}(x),19r_{m-2})$, which we can use to connect
$f_{m-2}(x)$ to $f_{m-2}(x')$ by a path $\gamma' \i \Sigma_{m-2}$ 
of length $L$, with
$$\eqalign{
L &\leq (1+C\varepsilon) |f_{m-2}(x)-f_{m-2}(x')|
\leq (1+C\varepsilon) (|f_{k}(x)-f_{k}(x')|+C\varepsilon r_{m-2})
\cr&\leq |f_k(x)-f_k(x')| + C \varepsilon r_k
}\leqno (11.23)
$$
(by (6.8) and because $|f_k(x)-f_k(x')| \leq 221 r_m$).
We choose less than $10^4$ points $z_l$ in $\gamma'$,
with consecutive distances less than $17 r_k$, so that the
first one is $f_{m-2}(x)$ and the last one is $f_{m-2}(x')$.
Each $z_l$ is of the form $f_{m-2}(x_l)$ for some
$x_l \in \Sigma_0$, so $z'_l = f_{k}(x_l)$ lies in $\Sigma_k$
and $|z'_l-z_l| \leq C \varepsilon r_k$ by (6.8). Now the
$z'_l$ are a string of points of $\Sigma_k$, whose consecutive 
distances are less than $18r_k$, and the total length of the
string is
$$\eqalign{
L' &= \sum |z'_{l+1}-z'_l| 
\leq \sum \big[C \varepsilon r_k + |z_{l+1}-z_l| \big]
\leq 10^4 C \varepsilon r_k + \sum |z_{l+1}-z_l|
\cr&
\leq 10^4 C \varepsilon r_k + L
\leq |f_k(x)-f_k(x')| + C \varepsilon r_k
}\leqno (11.24)
$$
by (11.23). We now get the desired curve $\gamma$ by
applying Lemma 6.12 to find a curve in $\Sigma_k$ of
length $(1+C\varepsilon)|z'_{l+1}-z'_l|$, that goes from
$z'_l$ to $z'_{l+1}$, and then putting all these curves
together. Notice that then
$$
{\rm length}(\gamma) \leq (1+C\varepsilon) L'
\leq (1+C\varepsilon) |f_k(x)-f_k(x')| + C \varepsilon r_k
\leq 2 |f_k(x)-f_k(x')|
\leqno (11.25)
$$
because $|f_k(x)-f_k(x')| > 18 r_k$.

We are left with the case when $m < 2$. But in this case,
we can use $\Sigma_0$, $x$, and $x'$ instead of $\Sigma_{m-2}$,
$f_{m-2}(x)$, and $f_{m-2}(x')$ above, because %vv
$|x'-x| \leq (1+C\varepsilon)|z'-z| \leq 1+C\varepsilon$
by assumption, and so we have more than enough control on 
$\Sigma_0 \cap B(x,2)$ to find $\gamma' \i \Sigma_0$).

Let us now use the curve $\gamma$ to estimate the right-hand side
of (11.21). Note that
$$\eqalign{
|R_k(x)-R_k(x')| &=
|\widetilde R_k(f_k(x))-\widetilde R_k(f_k(x'))|
\leq C \varepsilon r_k^{-1} {\rm length}(\gamma)
\cr&\leq C \varepsilon r_k^{-1} |f_k(x)-f_k(x')|
\leq C \varepsilon r_k^{-1} |f_m(x)-f_m(x')| 
}\leqno (11.26)
$$
by (9.34), (9.35),
(11.22), and (11.16), so (11.21) yields
$$\eqalign{
|\Delta_3| &\leq 
|y'|\sum_k\rho_k(y) \, |R_k(x)-R_k(x')|
\cr&
\leq C \varepsilon r_k^{-1} |f_m(x)-f_m(x')|\, |y'| %vv
%+ C \varepsilon |y'|
%\cr&
\leq C \varepsilon |f_m(x)-f_m(x')| %v
%+ C \varepsilon |y'|
}\leqno (11.27)
$$
because $|y'|\le |y| \leq 20 r_m$ by definition of $m$. Finally,
$$\eqalign{
\big|\Delta_4 \big| &= \Big|\sum_{k \geq 0} [\rho_k(y)-\rho_k(y')] 
\big\{ f_k(x') + R_k(x')\cdot y'\big\}\Big|
\cr&
= \Big|\sum_{k \geq 0} [\rho_k(y)-\rho_k(y')] 
\big\{ f_k(x')-f_m(x') + [R_k(x')-R_m(x')]\cdot y'\big\}\Big|
\cr&
\leq C \sum_{k = m}^{m+2} r_m^{-1}  |y-y'| 
\big[ |f_k(x')-f_m(x')| + |y'| \, |R_k(x')-R_m(x')| \big]
\cr&
\leq C \varepsilon |y-y'| + C \varepsilon r_m^{-1} |y-y'|\,|y'| 
\leq C \varepsilon |y-y'|
}\leqno (11.28)
$$
because $\sum_{k \geq 0} [\rho_k(y)-\rho_k(y')] = 0$,
and by (6.8) and (9.33). Altogether,
$$
|g(z)-g(z')-[f_m(x)-f_m(x')]-R_m(x)\cdot(y-y')| 
\leq C \varepsilon |f_m(x)-f_m(x')| +C \varepsilon |y-y'|
\leqno (11.29)
$$
by (11.15), (11.20), (11.27), and (11.28).
A first consequence of this is that
$$
|g(z)-g(z')| \leq (1+C\varepsilon) \big\{|f_m(x)-f_m(x')| + |y-y'| \big\}
\leqno (11.30)
$$
Next, the proof of (8.2) also yields that
$|f_m(x)-f_m(x')| \leq (1+C\varepsilon) |x-x'|^{1-C\varepsilon}$,
by Remark~8.31, so 
$$
|g(z)-g(z')| 
\leq (1+C\varepsilon) \big\{|x-x'|^{1-C\varepsilon} + |y-y'| \big\}.
\leqno (11.31)
$$
Recall from Lemma 10.1 that $|x-x'| \leq (1+C\varepsilon)|z-z'|$
and $|y-y'| \leq (1+C\varepsilon)|z-z'|$ (see (11.3) too).
Set $v(t) = \big\{|x-x'| + t \big\}^{1-C\varepsilon}$
for $0 \leq t \leq (1+C\varepsilon)|z-z'|$. Then
$v'(t) = (1-C\varepsilon) \big\{|x-x'| + t \big\}^{-C\varepsilon}
\geq (1-C\varepsilon)\, 3^{-C \varepsilon} \geq (1-C'\varepsilon)$
because $|x-x'| + t \leq 3 |z-z'| \leq 3$, so
$$\eqalign{
\big\{|x-x'| + |y-y'| \big\}^{1-C\varepsilon}
&= v(|y-y'|) \geq v(0) + (1-C'\varepsilon) |y-y'|
\cr&
= |x-x'|^{1-C\varepsilon} + (1-C'\varepsilon) |y-y'|
\cr&
\geq (1-C'\varepsilon)\big\{|x-x'|^{1-C\varepsilon} + |y-y'| \big\}
}\leqno (11.32)
$$
and hence, by (11.31),
$$\eqalign{
|g(z)-g(z')| 
& \leq (1+C\varepsilon) \big\{|x-x'| + |y-y'| \big\}^{1-C\varepsilon}
\leq 3 |z-z'|^{1-C\varepsilon}
}\leqno (11.33)
$$
because $|x-x'| + |y-y'| \leq 2 (1+C\varepsilon)|z-z'|$
by Lemma 10.1. %v
This gives the second inequality in (11.2).

\ms
We now look for lower bounds. Let us first project (11.29) on 
the tangent direction $T_m(x)$ to $\Sigma_m$ at $f_m(x)$. %v
Call $\pi$ the orthogonal projection
onto $T_m(x)$ (it was also called $\pi_{f_m(x),m}$ before), %v
and similarly set $\pi' = \pi_{f_m(x'),m}$. Observe that
$$
|\pi - \pi'| = |\pi_{f_m(x),m}-\pi_{f_m(x'),m}|
\leq C \varepsilon r_m^{-1}|f_m(x)-f_m(x')|
\leqno (11.34)
$$
by (9.28). Since $y$ is orthogonal to $T_0(x)$, (9.32) says that
$R_m(x)\cdot y$ is orthogonal to $T_m(x)$, and similarly
$R_m(x')\cdot y'$ is orthogonal to $T_m(x')$. Now
$$\eqalign{
|\pi(R_m(x)\cdot(y-y'))| &= |\pi(R_m(x)\cdot y')|
\leq |\pi'(R_m(x)\cdot y')|+ |\pi - \pi'||R_m(x)\cdot y'|
\cr& = |\pi'([R_m(x)-R_m(x')]\cdot y')|+ |\pi-\pi'||y'|
\cr& \leq \big[|R_m(x)-R_m(x')|+|\pi-\pi'| \big] \, |y'|
\cr&\leq C \varepsilon r_m^{-1}|f_m(x)-f_m(x')| \, |y'|
+|\pi-\pi'| \, |y'|
\cr&\leq C \varepsilon r_m^{-1}|f_m(x)-f_m(x')| \, |y'|
\leq C \varepsilon |f_m(x)-f_m(x')| 
%&(11.35)
}\leqno (11.35)
$$
because $R_m(x)$ is an isometry and then $\pi'(R_m(x')\cdot y')=0$ 
by orthogonality, by the first inequalities of (11.26) for $k=m$, 
(11.34), and the fact that $|y'|\le |y| \leq 20 r_m$.

Next we care about $\pi(f_m(x)-f_m(x'))$. Recall from
(11.10) that $|f_m(x)-f_m(x')| \leq 220 r_m$. 
If $|f_m(x)-f_m(x')| > 18 r_m$, choose $k = m-1$
or $k = m-2$, as large as possible, so that
$|f_m(x)-f_m(x')| \leq 18 r_k$; otherwise, choose
$k=m$. Note that $k \geq 0$, because 
$|x'-x| \leq (1+\varepsilon)|z'-z| \leq 1+\varepsilon$
and hence $|f_m(x)-f_m(x')| \leq 2$.

Apply Lemma 6.12 to the integer $k$ chosen above and the point
$\xi=f_k(x)$. We get that $\Sigma_k$ coincides with a 
$C\varepsilon$-Lipschitz graph $\Gamma$ over some plane 
$P$ through $\xi$. Note that
$$
{\rm Angle}(P',T_{k}(x)) =
{\rm Angle}(P',T\Sigma_{k}(f_k(x)))
\leq C \varepsilon,
\leqno (11.36)
$$
where we denote by $P'$ the vector space
parallel to $P$, and we remind the reader that 
$T_{k}(x) = T\Sigma_{k}(f_{k}(x))$ by (9.30).

If $k=m$, we immediately get that
$$
{\rm Angle}(f_m(x)-f_m(x'),P') \leq C \varepsilon
\leqno (11.37)
$$
(or $f_m(x)=f_m(x')$, but then (11.39) below is trivial), 
because both $f_m(x)$ and $f_m(x')$ 
lie on $\Gamma$. If $k < m$, we only get that 
${\rm Angle}(f_k(x)-f_k(x'),P') \leq C \varepsilon$
(because $f_k(x), f_k(x') \in \Gamma$), but since
$|f_m(x)-f_k(x)|+|f_m(x')-f_k(x')|\leq C \varepsilon r_k$ by
(6.8) and $|f_m(x)-f_m(x')| > 18 r_m$ because $k<m$,
we also get (11.37). Now
$$\eqalign{
{\rm Angle}(f_m(x)-f_m(x'),T_{m}(x)) 
&\leq C \varepsilon + {\rm Angle}(P',T_{m}(x)) 
\cr&\leq C' \varepsilon + {\rm Angle}(P',T_{k}(x)) 
\leq C'' \varepsilon
}\leqno (11.38)
$$
by (11.37), (9.3) and the definition (9.30), and (11.36).
Hence
$$
|\pi(f_m(x)-f_m(x'))| \geq {9 \over 10}|f_m(x)-f_m(x')|.
\leqno (11.39)
$$
Altogether,
$$\leqalignno{
|g(z)-g(z')| &\geq |\pi(g(z)-g(z'))| 
\cr& \geq |\pi(f_m(x)-f_m(x')+R_m(x)\cdot(y-y'))| 
- C \varepsilon |f_m(x)-f_m(x')| - C \varepsilon |y-y'|
\cr&
\geq |\pi(f_m(x)-f_m(x'))| 
- C \varepsilon |f_m(x)-f_m(x')| - C \varepsilon |y-y'|
& (11.40)
\cr& \geq {9-C\varepsilon \over 10} \, |f_m(x)-f_m(x')| 
-C \varepsilon |y-y'|
}%\leqno (11.40)
$$
by (11.29), (11.35), and (11.39). If $|y-y'|\leq 2|f_m(x)-f_m(x')|$,
we get that
$$\eqalign{
|g(z)-g(z')| &\geq {8 \over 10} \, |f_m(x)-f_m(x')|
\geq {8 \over 30}\, \big\{ |f_m(x)-f_m(x')| + |y-y'| \big\}
\cr&
\geq {8 \over 30} \, 
\big\{ (1-C\varepsilon) |x-x'|^{1+C\varepsilon} + |y-y'| \big\}
}\leqno (11.41)
$$
by (8.2). On the other hand, recall that
$|x-x'| \leq (1+C\varepsilon) |z'-z| \leq 1+C\varepsilon$
by Lemma 10.1 and (11.3), and similarly $|y-y'| \leq 1+C\varepsilon$.
Set $v(t) = (|x-x'|+t)^{1+C\varepsilon}$
for $0 \leq t \leq 1+C\varepsilon$; then
$v'(t) = (1+C\varepsilon) (|x-x'|+t)^{C\varepsilon}
\leq 1+C'\varepsilon$, so the fundamental theorem of calculus yields
$$\eqalign{
|z-z'|^{1+C\varepsilon} &\leq (|x-x'|+|y-y'|)^{1+C\varepsilon}
= v(|y-y'|) 
\leq v(0) + (1+C'\varepsilon) |y-y'|
\cr&=|x-x'|^{1+C\varepsilon} + (1+C'\varepsilon)|y-y'|
\leq 4 |g(z)-g(z')|,
}\leqno (11.42)
$$
by (11.41). So the first part of (11.2) holds in this case.

We may finally assume that $|y-y'|> 2|f_m(x)-f_m(x')|$. Then 
(11.29) implies that
$$\leqalignno{
|g(z)-g(z')| &\geq |R_m(x)\cdot(y-y')|-|f_m(x)-f_m(x')|
- C \varepsilon |f_m(x)-f_m(x')| - C \varepsilon |y-y'|
\cr& \geq |y-y'|-|f_m(x)-f_m(x')|-C \varepsilon |y-y'|
\geq \big({1 \over 2} - C \varepsilon \big) |y-y'|
&(11.43)
}%\leqno (11.43)
$$
because $R_m(x)$ is an isometry, and also
$$\eqalign{
|g(z)-g(z')| &\geq |y-y'|-|f_m(x)-f_m(x')|-C \varepsilon |y-y'|
\cr&\geq (1-C\varepsilon) |f_m(x)-f_m(x')|
\geq (1-C'\varepsilon) |x-x'|^{1+C\varepsilon}
}\leqno (11.44)
$$
by the beginning of (11.43) and (8.2). To end the estimate, 
we multiply (11.43) by $9/15$, multiply (11.44) by $6/15$, add the
two, and get an estimate better than (11.41), which as we 
already know implies the first part of (11.2). 
This completes our proof of Proposition~11.1.
\qed
\ms
The bi-Lipschitz version of Proposition 11.1 will be easier.

\ms
\proclaim Proposition 11.45. Suppose that (8.35) holds
for some $M < +\infty$. Then $g : \R^n \to \R^n$ is bi-Lipschitz. 

\ms 
Recall that the $\varepsilon'_k(f_k(z))$ are as in (7.18),
and $\varepsilon'_k(f_k(z)) = 0$ when $f_k(z) \in \R^n \sm V_k^{10}$.
The condition is the same as in Proposition 8.34, so we know
that $f : \Sigma_0 \to \Sigma$ is bi-Lipschitz.
The estimates used for Proposition 11.1 are still valid
now; we just need to conclude differently.

As before, we may assume that $0 < |y'| \leq |y|$, and we start
with the case when $|f(x')-f(x)| \geq 10 |y|$. Then the second line 
of (11.5) yields
$$
|g(z')-g(z)| \leq {13 \over 10} \, |f(x')-f(x)|
\leq C |x'-x| \leq 2C|z'-z|
\leqno (11.46)
$$
by Lemma 10.1. Similarly, the beginning of (11.6) yields
$$
|g(z')-g(z)| \geq {8-C\varepsilon \over 10} \, |f(x')-f(x)|
\geq C^{-1} |x'-x|
\leqno (11.47)
$$
(again by Proposition 8.34), and at the same time 
$$
|g(z')-g(z)| \geq {8-C\varepsilon \over 10} \, |f(x')-f(x)|
\geq (8-C\varepsilon) \, |y| 
\geq {8-C\varepsilon \over 2} \, |y'-y|
\leqno (11.48)
$$
because  $|f(x')-f(x)| \geq 10 |y|$ and $|y'| \leq |y|$. Then
$|g(z')-g(z)| \geq C^{-1}|z'-z|$ because obviously 
$|z'-z| \leq |x'-x|+|y'-y|$.

When instead $|f(x')-f(x)| < 10 |y|$, (11.30) immediately yields
$$\eqalign{
|g(z)-g(z')| 
&\leq (1+C\varepsilon) \big\{|f_m(x)-f_m(x')| + |y-y'| \big\}
\cr&\leq C \big\{ |x-x'| + |y-y'| \big\} \leq 3C |z-z'|
}\leqno (11.49)
$$
by Lemma 10.1. We are left with the lower bound for this second case.
If $|y-y'|\leq 2|f_m(x)-f_m(x')|$, the first line of (11.41) yields
$$\eqalign{
|g(z)-g(z')| 
&\geq {8 \over 30}\, \big\{ |f_m(x)-f_m(x')| + |y-y'| \big\}
\cr&
\geq C^{-1} \, \big\{ |x-x'| + |y-y'| \big\}
\geq C^{-1} \, |z-z'|,
}\leqno (11.50)
$$
as needed. Finally, if $|y-y'| > 2|f_m(x)-f_m(x')|$,
(11.43) implies that $|g(z)-g(z')| \geq {1 \over 3}\, |y-y'|$
and then also, as in the first part of (11.44),
$$
|g(z)-g(z')|\geq {2 \over 3}\, |f_m(x)-f_m(x')|
\geq C^{-1} |x-x'|,
\leqno (11.51)
$$
so we may conclude as in the previous case.
\qed

\msi{\bf Proof of Theorems 2.15 and 2.23.} 
Let us just observe here that we completed the proof
of these two theorems. For Theorem 2.15, the hypotheses are 
the same as throughout Sections 3-11; (2.16) follows from
(10.13) and (10.22), (2.17) follows from (10.23), 
(2.18) is the same as (11.2), and $\Sigma$ contains $E_\infty$
by (6.2). For Theorem 2.23, we added the assumption (2.24),
which is the same as (8.35), and required that $g$ be bi-Lipschz;
this is proved in Proposition 11.45.
\qed

\bsi
{\bf 12. Variants of the Reifenberg theorem}
\ms

In this section we want to state and prove a few variants of
Reifenberg's topological disk theorem. We tried to arrange things so 
that the statements will be easy to read independently from the 
previous sections; of course the proofs will not.

For all the statements, we are given a smooth $d$-dimensional
manifold $\Sigma_0 \i \R^n$, and we assume (exactly as in Section 2) 
that 
$$\eqalign{
&\hbox{for every $x\in \Sigma_0$, there is an affine $d$-plane
$P_x$ through $x$ and}
\cr&\hbox{a $C^2$ function $F_x : P_x \to P_x^\perp,$ such that
(2.4) and (2.5) hold.}
}\leqno (12.1)
$$
As usual, $P_x^\perp$ is the vector space of dimension $n-d$
which is perpendicular to $P_x$. Recall that (2.4) and (2.5) 
say that in $B(x,200)$, $\Sigma_0$ coincides with an 
$\varepsilon$-Lipschitz graph, with a similar estimate on
the graphed function and its second derivative. 
Thus (12.1) is a quantitative way
to require that $\Sigma_0$ be quite flat at the unit scale.
The constant $\varepsilon > 0$ will need to be small enough,
depending on $n$ and $d$.

The most standard example of set $\Sigma_0$ is undoubtedly
a $d$-plane, but it does not cost us much to allow more
complicated manifolds $\Sigma_0$. Note however that all the 
complication occurs at large scales, and that our construction 
is local, so the apparent generality is not too shocking.

Next, we are given a set $E \i \R^n$ that we want to study, 
and on which we shall make various flatness assumptions. 
Finally, we shall use a set $U \i \R^n$ to localize the statements. 
A typical choice of $U$ would be a large ball. 
We shall not put specific conditions on $U$, but since the 
conclusions will occur on $U$ and the hypotheses will be made on
$$
U^+ = \big\{ x\in \R^n \, ; \, \dist(x,U) \leq 2 \big\},
\leqno (12.2)
$$
it is not in our interest to take a complicated $U$.
Our first statement is a generalization of Theorem 1.1.

\ms\proclaim Theorem 12.3. Let $\varepsilon > 0$ be small
enough, depending on $n$ and $d$.
Let $E, U, \Sigma_0 \in \R^n$ be given, and assume that 
(12.1) holds. Also assume that 
(12.1) holds. Also assume that 
$$
\dist(x,\Sigma_0) \leq \varepsilon
\hbox{ for } x\in E \cap U^+
\ \hbox{ and } \ 
\dist(x,E) \leq 1/2 
\hbox{ for } x\in \Sigma_0 \cap U^+,
\leqno (12.4)
$$
and that for $x\in E \cap U^+$ and $r\in (0,1]$, there is
an affine $d$-plane $P = P(x,r)$ through $x$ such that
$$\eqalign{
&\dist(y,P) \leq \varepsilon r \hbox{ for } y\in E\cap B(x,110r)
\cr& \hskip 3cm 
\hbox{ and }
\dist(y,E) \leq \varepsilon r \hbox{ for } y\in P\cap B(x,110r).
}\leqno (12.5)
$$
Then there is a bijective mapping $g : \R^n \to \R^n$
such that
$$
g(x) = x \ \hbox{ when } \dist(x,U) \geq 13 
\leqno (12.6)
$$
$$
|g(x)-x| \leq C\varepsilon \ \hbox{ for } x\in \R^n,
\leqno (12.7)
$$
$$
{1 \over 4} |x'-x|^{1+C\varepsilon} \leq |g(x')-g(x)| 
\leq 3 |x'-x|^{1-C\varepsilon}  
\leqno (12.8)
$$
for $x,x'\in \R^n$ such that $|x'-x| \leq 1$, and 
$$
\overline E \cap U = g(\Sigma_0) \cap U.
\leqno (12.9)
$$
The constant $C$ depends only on $n$ and $d$.

\ms
Let us make a few comments before we prove this theorem.
We decided not to require $E$ to be closed, but replacing
$E$ with its closure essentially does not change the hypotheses 
or the conclusion.

When $U=\R^n$, (12.9) just says that $\overline E = g(\Sigma_0)$,
so we have a good parameterization of $E$ by $\Sigma_0$, which
extends to a bi-H\"older homeomorphism of $\R^n$.

For the proof we do not really need $\Sigma_0$ to be a manifold
everywhere, because we only need to know (2.4) and (2.5) at points
$x\in \Sigma_0$ such that $\dist(x,E) \leq 1$.

The constants ${1\over 4}$ and $3$ in (12.8) are not optimal, and
can probably be replaced with constants that are arbitrarily
close to $1$ (even with the function $g$ constructed
above). See the remark below Proposition 11.1 for 
a hint on how to start a proof. Our constants $100$ and $110$
look annoying, but we could easily make them smaller by appling
a dilation to $E$, $\Sigma_0$, and $U$. %vv

\ms 
We now prove Theorem 12.3. 
We want to construct a CCBP (see Definition 2.11).
We already have the set $\Sigma_0$, with (2.4) and (2.5).
Next we choose the $x_{j,k}$. Set
$$
E_0 = \big\{ x\in E \, ; \, \dist(x,U) \leq 1 \big\}
\leqno (12.10)
$$
and, for $k \geq 0$, let $\{ x_{j,k} \}$, $j\in J_k \,$, be a maximal
subset of $E_0$ with the constraint that $|x_{i,k}-x_{j,k}| \geq r_k$.
[Recall that $r_k = 10^{-k}$.] By maximality,
$$
E_0 \i \bigcup_{j\in J_k} \overline B(x_{j,k},r_k)
\leqno (12.11)
$$
for each $k \geq 0$, and in particular (2.3) holds. 
Also, (2.7) follows from (12.4).

For $j\in J_k$, we choose a $d$-plane $P_{j,k}$ such that
(12.5) holds with $P=P_{j,k}$, $x=x_{j,k}$, and $r=r_k$;
such a plane exists precisely by assumption. Now we need
to check (2.8)-(2.10). We shall use the following lemma,
whose standard elementary proof is left to the reader.

\ms\proclaim Lemma 12.12.
Let $P_1$ and $P_2$ be affine $d$-planes. Let $z \in P_1$
and $r >0$, and suppose that for some $\tau < 1$,
$$
\dist(y,P_2) \leq \tau r
\ \hbox{ for } y\in P_1 \cap B(z,r).
\leqno (12.13)
$$
Then $d_{z,200r}(P_1,P_2) \leq C \tau$.

\ms
As usual, $C$ may depend on $n$ and $d$, but not on $\tau$
$r$ or $z$. See (1.7) for the definition of $d_{z,100r}(P_1,P_2)$
and the proof of (7.25) for a hint on how to start. %v

First we prove (2.8). Let $i,j \in J_k$ be such that 
$|x_{i,k}-x_{j,k}| \leq 100r_k$, and let us try to apply
Lemma 12.12. By (12.5), we can find $z\in P_{j,k}$ such that
$|z-x_{j,k}| \leq \varepsilon r_k$. Then, for each
$y\in P_{j,k} \cap B(z,r_k)$, (12.5) gives $y'\in E$ such that
$|y'-y| \leq \varepsilon r_k$ and, since $y' \in 110B_{i,k}$
because $|x_{i,k}-x_{j,k}| \leq 100r_k$, a new application of
(12.5) gives $y''\in P_{i,k}$ such that $|y''-y'| \leq \varepsilon r_k$.
Thus Lemma 12.12 applies with $r = r_k$ and $\tau = 2 \varepsilon$;
we get that 
$$
d_{x_{j,k},100 r_k}(P_{i,k},P_{j,k})
\leq 2 d_{z,200 r_k}(P_{i,k},P_{j,k}) \leq C \varepsilon
\leqno (12.14)
$$
because $B(x_{j,k},100 r_k) \i B(z,200 r_k)$.

For (2.9), let $i\in J_0$ and $x\in \Sigma_0$ be such
that $|x_{i,0}-x| \leq 2$. We want to apply
Lemma 12.12 to control $d_{x_{i,0},100}(P_{i,0},P_x)$.
First use (12.5) to choose $z\in P_{i,0}$ such that 
$|z-x_{i,0}| \leq \varepsilon$. Then, for 
$y\in P_{i,0} \cap B(z,2/3)$, (12.5) gives $y'\in E$
such that $|y'-y| \leq \varepsilon$. Note that $y'\in U^+$,
because $\dist(x_{i,0},U) \leq 1$ (by (12.10)). So
(12.4) applies, and gives $y'' \in \Sigma_0$ such that 
$|y''-y'| \leq \varepsilon$. Finally, by (2.4) and (2.5),
we can find $w\in P_x$ such that $|w-y''| \leq \varepsilon$.
Altogether, Lemma 12.12 applies to $P_{i,0}$, $P_x$, $z$,
and $r = 2/3$. Thus
$$
d_{x_{i,0},100}(P_{i,0},P_x) 
\leq d_{z,400/3}(P_{i,0},P_x) \leq C \varepsilon
\leqno (12.15)
$$
because $B(x_{i,0},100) \i  B(z,400/3)$.

Finally we prove (2.10) the same way. Let $i\in J_k$ and 
$j\in J_{k+1}$ be such that $|x_{i,k}-x_{j,k+1}| \leq 2 r_k$.
Choose $z\in P_{i,k}$ such that 
$|z-x_{j,k+1}| \leq \varepsilon r_k$.
For $y\in P_{j,k+1} \cap B(z,r_{k})$, (12.5) gives
$y\in E$ such that $|y'-y| \leq \varepsilon r_{k}$
and, since $y'\in E \cap 3B_{i,k}$, we also get
$y'' \in P_{i,k}$ such that $|y''-y'| \leq \varepsilon r_k$.
So Lemma 12.12 applies, and
$$
d_{x_{i,k},20r_k}(P_{i,k},P_{j,k+1})
\leq d_{z,200r_k}(P_{i,k},P_{j,k+1}) \leq C \varepsilon,
\leqno (12.16)
$$
as needed. This completes the verification of the CCBP
conditions (see Definition 2.11). The fact that we only
obtained (2.8)-(2.10) with the constant $C\varepsilon$
does not matter.

We may now apply Theorem 2.15 to the CCBP at hand, and
we get a mapping $g$ for which we now check (12.6)-(12.9).

Set $\Sigma' = \big\{ x\in \Sigma_0 \, ; \, \dist(x,E_0) \leq 10 
\big\}$. Let us check that
$$
g(z)=z \ \hbox{ when } \dist(z,\Sigma') \geq 2;
\leqno (12.17)
$$
obviously (12.6) will follow, because $\dist(x,U) \leq 11$
for $x\in \Sigma'$ (by (12.10)).

For $x\in \Sigma_0 \sm \Sigma'$, (4.5) says that $\sigma_k(x)=x$ and 
$D\sigma_k(x) = I$ for $k \geq 0$, so $f_k(x) = x$ and the successive 
tangent directions $T_k(x)$ are all equal to $T_0(x)$. 
The construction of $R_k(x)$ yields $R_k(x) = I$ for all $k$ 
(notice in particular that if $R_k(x) = I$, (9.36) yields 
$S_k(x) = I$, which in turn yields $R_{k+1}(x)=I$ by (9.45)).

Now let $z\in \R^n$ be such that $\dist(z,\Sigma') \geq 2$,
and let us check that $g(z)=z$. If $z\in \Sigma_0$, 
$g(z) = f(z)$ by (10.14), and $f(z) = z$ because 
$z\in \Sigma_0 \sm \Sigma'$. We may thus assume that
$z \in V$ and $\dist(z,\Sigma_0) \leq 2$, because otherwise
$g(z)=z$ by (10.13) or (10.22). Now $|p(z)-z| \leq 2$, by (10.12),
so $p(z) \in \Sigma_0 \sm \Sigma'$, and by the discussion above
$f_k(p(z)) = p(z)$ and $R_k(p(z)) = I$, so that 
$g(z) = \sum_k \rho_k(q(z))\big\{p(z)+q(z) \big\} = z$
by (10.15) and (10.17). Thus (12.17) and (12.6) hold.

Next (12.7) and (12.8) are the same as (2.17) and (2.18)
(or (10.23) and (11.2)), so we are left with (12.9) to check. 
If $x\in \overline E \cap U$,
(12.11) says that for each $k \geq 0$, we can find $j\in J_k$
such that $|x-x_{j,k}| \leq r_k$. Then $x\in E_\infty$,
the limit set defined by (2.19), and Theorem 2.15 says that
$x\in \Sigma = g(\Sigma_0)$.

Conversely, let $w\in g(\Sigma_0) \cap U$ be given,
and set $d = \dist(w,E)$. Thus we want to show that $d=0$.
Let $z\in \Sigma_0$ be such that $w=g(z)=f(z)$. Observe that
$|z-w| \leq C \varepsilon$ by (12.7) or (6.9), so 
$d \leq 1/2+C\varepsilon < 2/3$ by (12.4). Suppose that $d >0$, 
and let $k \geq 0$ be such that $r_{k+1} \leq d \leq r_k$.

By (6.9), $|w-f_k(z)| \leq C \varepsilon r_k$. By
definition of $d$, we can find $\xi \in E$ such that
$|\xi-w| \leq 3d/2$. Notice that then $\xi \in E_0$,
because $w\in U$ and $d < 2/3$.
Then $\xi \in \overline B_{j,k}$ for some $j \in J_k$
(by (12.11)), and $f_k(z) \in \Sigma_k \cap 3B_{j,k}$
because $|x_{j,k}-f_k(z)| \leq |x_{j,k}-\xi|+|\xi-w|+|w-f_k(z)|
\leq r_k + 3d/2 + C \varepsilon r_k < 3r_k$ because
$d \leq r_k$. Thus Proposition 5.4 says
that $\dist(f_k(z), P_{j,k}) \leq C \varepsilon r_k$.
Choose $y\in P_{j,k}$ such that $|y-f_k(z)| \leq C \varepsilon r_k$;
obviously $y\in 4B_{j,k}$, so by (12.5) and our choice of
$P_{j,k}$, we can find $y' \in E$ such that 
$|y'-y| \leq \varepsilon r_k$. Finally,
$d \leq |y'-w| \leq|y'-y|+|y-f_k(z)|+|f_k(z)-w|
\leq C \varepsilon r_k$, which contradicts the
definition of $r_k$ and proves that $d=0$.

This completes our proof of Theorem 12.3.
\qed

\ms 
Next we generalize Theorem 1.10.

\ms\proclaim Theorem 12.18. Let $\varepsilon > 0$ be small
enough, depending on $n$ and $d$.
Let $E$ and $\Sigma_0 \in \R^n$ be given, and assume that 
(12.1) holds. Also assume that 
$$
\dist(x,\Sigma_0) \leq \varepsilon
\hbox{ for } x\in E ,
\leqno (12.19)
$$
that for $x\in E$ and $k \geq 0$, we are given
an affine $d$-plane $P_k(x)$ through $x$ such that,
with the notation (1.7) for local Hausdorff distances,
$$
d_{x,100r_k}(P_k(x),P_k(x')) \leq \varepsilon
\ \hbox{ for $k \geq 0$ and } x,x' \in E
\hbox{ such that } |x'-x| \leq 100 r_k,
\leqno (12.20)
$$
$$
d_{x,r_k}(P_k(x),P_{k+1}(x)) \leq \varepsilon
\ \hbox{ for $k \geq 0$ and } x \in E,
\leqno (12.21)
$$
and
$$
d_{x,100}(P_0(x),P_y) \leq \varepsilon
\ \hbox{ for $x \in E$ and $y\in \Sigma_0$ such that } 
|x-y| \leq 2,
\leqno (12.22)
$$
where $P_y$ is as in the description of $\Sigma_0$
in (12.1), (2.4), and (2.5). Then there is a bijective 
mapping $g : \R^n \to \R^n$ such that
$$
g(x) = x \ \hbox{ when } \dist(x,E) \geq 12,
\leqno (12.23)
$$
$$
|g(x)-x| \leq C\varepsilon \ \hbox{ for } x\in \R^n,
\leqno (12.24)
$$
$$
{1 \over 4} |x'-x|^{1+C\varepsilon} \leq |g(x')-g(x)| 
\leq 3 |x'-x|^{1-C\varepsilon}  
\leqno (12.25)
$$
for $x,x'\in \R^n$ such that $|x'-x| \leq 1$, and 
$$
E \i g(\Sigma_0).
\leqno (12.26)
$$
In addition, $\Sigma = g(\Sigma_0)$ is Reifenberg-flat, in 
the sense that for $x \in \Sigma$ and $r\in (0,1]$, there is an 
affine $d$-plane $Q(x,r)$ through $x$ such that 
$d_{x,r}(\Sigma,Q(x,r)) \leq C \varepsilon$.
The constant $C$ in (12.24) and (12.25) depends only on $n$ and $d$.

\ms

Note that if in each ball centered on $E$ there are 
$d+1$ ``sufficiently" affinely independent points of $E$, then 
conditions (12.20) and (12.21) are automatically satisfied.
But in general, something like (12.20)-(12.22) is needed;
see Counterexample 12.28.
\ms

Recall that $r_k = 10^{-k}$ for $k \geq 0$.
Here we did not say that points of $E \cap B(x,r_k)$
lie close to $P_k(x)$, but this is implied by (12.20),
because $P_k(x')$ contains $x'$. Also, we do not need
to localize this theorem, as we could just have restricted
our attention to $E\cap U$.

Theorem 12.18 is stronger than Theorem 1.10. If $E$ is as in 
Theorem 1.10, take $\Sigma_0 = P(0,10)$; then (12.1) is obvious,
and (12.19) follows from (1.6) for $P(0,10)$,
(12.20) is the same as (1.9), (12.21) is the same as (1.8)
for $k \geq 1$, and (12.22) holds by (1.8) for $k-1$.
The conclusions of Theorem 12.18 are stronger; in particular
we can also take $Q(x,r) = P(0,10)$ for $r \geq 1$,
by (12.24).

Now we prove Theorem 12.18. As before, we already have $\Sigma_0$,
and for each $k \geq 0$, we choose a maximal collection
$\{ x_{j,k} \}, j\in J_k$ of points of $E$, with the constraint 
that $|x_{i,k}-x_{j,k}| \geq r_k$.
Then we set $P_{j,k} = P_k(x_{j,k})$ for $k \geq 0$ and
$j \in J_k$; (2.3) comes from the maximality of the collection
$\{ x_{j,k} \}, j\in J_k$, (2.8) follows from (12.20), (2.9)
comes from (12.22), and for (2.10) we observe that if
$i\in J_k$ and $j \in J_{k+1}$ are such that 
$|x_{i,k}-x_{j,k+1}| \leq 2r_k$, then
$$\leqalignno{
d_{x_{i,k},20r_k}(&P_{i,k},P_{j,k+1})
= d_{x_{i,k},20r_k}(P_k(x_{i,k}),P_{k+1}(x_{j,k+1}))
\cr&
\leq 2d_{x_{i,k},40r_k}(P_k(x_{i,k}),P_{k}(x_{j,k+1}))
+2d_{x_{i,k},40r_k}(P_{k}(x_{j,k+1}),P_{k+1}(x_{j,k+1}))
\cr&
\leq 5 d_{x_{i,k},100r_k}(P_k(x_{i,k}),P_{k}(x_{j,k+1}))
+ 3d_{x_{j,k+1},50r_k}(P_{k}(x_{j,k+1}),P_{k+1}(x_{j,k+1}))
&(12.27)
\cr& \leq 5 \varepsilon + 
C d_{x_{j,k+1},r_k}(P_{k}(x_{j,k+1}),P_{k+1}(x_{j,k+1}))
\leq C\varepsilon
}%\leqno (12.27)
$$
by the definition (1.7) of $d$ and the triangle inequality, 
because $B(x_{i,k},40r_k) \i B(x_{j,k+1},50r_k)$ and, 
for the last line, (12.20), simple geometry using the fact 
that we are computing distances between $d$-planes, and (12.21).

So we have a CCBP (as in Definition 2.11), and 
Theorem 2.15 gives a mapping $g$. As before, (12.24) and (12.25)
are the same as (2.17) and (2.18). Concerning (12.26), observe that
for $x\in E$ and $k \geq 0$, there is an $x_{j,k}$ such that
$|x-x_{j,k}| \leq r_k$ (by maximality of the family 
$\{ x_{j,k} \}, j\in J_k$), so $x\in E_\infty$
(the limit set from (2.19)), and Theorem 2.15 says that %v
$x\in \Sigma = g(\Sigma_0)$.

Next (12.23) is proved as (12.6) above: first one checks that
$f_k(x)=x$ and $T_k(x) = T_0(x)$ for $x\in \Sigma_0$
such that $\dist(x,E) \geq 10$ (and all $k\geq 0$), and then
one gets that $g(z) = z$ unless $z\in V$, $\dist(z,\Sigma_0) \leq 2$,
and $p(z) \in \Sigma' = \big\{ x\in \Sigma_0 \, ; \, 
\dist(x,E) \leq 10 \big\}$. See the proof of (12.17).
Finally $\Sigma$ is Reifenberg-flat by Proposition 6.15.
Theorem 12.18 follows.
\qed

\msi{\bf Counterexample 12.28.}
The coherence conditions (12.20)-(12.22) are really needed
in the statement of Theorem 12.18. Let us construct a 
two-dimensional set $E \i \R^3$
such that for every $x\in E$ and $r > 0$, there is a plane
$P(x,r)$ such that
$$
\dist(y,P(x,r)) \leq \varepsilon r
\ \hbox{ for } y \in E \cap B(x,r),
\leqno (12.29)
$$
but $E$ is not contained in a any Reifenberg-flat set $\Sigma$.

Let $P_0$ be the horizontal plane through the origin,
set $S = P_0 \cap \partial B(0,1)$, and let $E$ be a M\"obius strip
of very small width $\tau > 0$ whose central curve is $S$.
Choose $E$ so that if $T(x)$ denotes the direction of 
the tangent plane $T(x)$ at $x\in E$, $|DT(x)| \leq 10$, say.

Here we take $\Sigma_0 = P_0$; note that
$$
\dist(x,P_0) \leq \tau
\ \hbox{ for } x\in E,
\leqno (12.30)
$$
so there is no difficulty with the initial condition (2.19). 
The approximation by planes is fine too. For
$r \geq \varepsilon^{-1}\tau$, we simply choose
$P(x,r) = P_0$ and use (12.30), while for $r < \varepsilon^{-1}\tau$
(and if $\tau < c \varepsilon^2$) we can choose the tangent
plane to $E$ at $x$ and use the slow variation of $T(x)$.

Now $E$ is not contained in a Reifenberg-flat set $\Sigma = g(P_0)$,
simply because it is not orientable. The reason why it does not 
satisfy the assumption of Theorem 12.18 is similar: there is no
nearly continuous choice of $P(x,r)$, $x\in E$ and $0 < r < 1$,
that coincides with the choices above for $t$ small and large.

%% second choice for the picture
The reader may wonder whether things got wrong here because we did
not choose the right model $\Sigma_0$, but this is not so.
We can construct a different counterexample as follows.
See Figure 1.  
Start from $P_0$ as above, choose a tiny square $Q \i P_0$
of sidelength $l$, choose two opposite sides of $\partial Q$,
and let $I \i Q$ denote the interval that connects the middles
of these two sides. Let $H_0 \i Q$ denote the very thin stripe
of width $\tau l$ centered along $I$, and let $H$ be obtained
from $H_0$ by twisting it one half turn around $I$ (and fairly
regularly). Finally set $E = (P_0 \sm Q) \cup H$.
As before, if $\tau < c \varepsilon^2$, we can find planes 
$P(x,r)$ such that (12.29) holds, and yet $\Sigma$ is not
contained in a Reifenberg-flat set because $E$ is not orientable.
\vglue 2.5cm %%% verifier le display page 63
\hskip 3.0cm
\epsffile{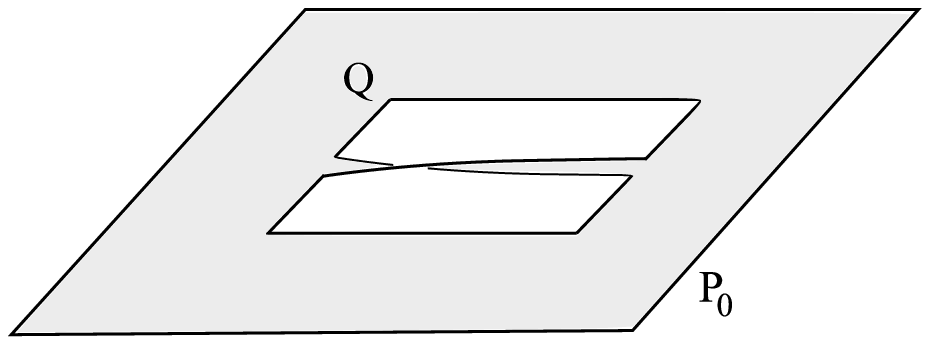}
\medskip
\centerline{{\bf Figure 1.} The set $E$}
\medskip

\ms
Next we want to mention sufficient conditions for
$g$ above to be bi-Lipschitz. Let us use notation
that fits both Theorems~12.3 and 12.18. When 
$E$ and $U$ are as in Theorem~12.3,
$x \in E_0 = \big\{ x\in E \, ; \, \dist(x,U) \leq 1 \big\}$
(as in (12.10)) and $k\geq 0$, we choose a $d$-plane
$P_k(x)$ so that (12.5) holds for $P =P_k(x)$ and $r = r_k$.
In the context of Theorem 12.18, set $E_0 = E$; we already chose 
$P_k(x)$ for $x\in E$ and $k \geq 0$ in the statement.
In both cases, set
$$\eqalign{
\alpha_k(x) &= d_{x,r_k}(P_{k+1}(x),P_k(x)) 
%\cr& \hskip 2cm
+ \sup_{y\in E_0\cap B(x,35r_k)} d_{x,r_{k}}(P_{k}(x),P_k(y)) 
}\leqno (12.31)
$$
for $x\in E_0$ and $k \geq 0$, and then
$$
J(x) = \sum_{k \geq 0} \alpha_k(x)^2.
\leqno (12.32)
$$
(essentially as in (1.19) and (1.20)).

\ms\proclaim Corollary 12.33.
If, in addition to the hypotheses of Theorem 12.3 or 12.18,
we have that $J(x) \leq M$ for $x\in E_0$ (and for some choice of 
planes $P_k(x)$ as above), then $g$ is bi-Lipschitz. More precisely, 
there is a constant $K$, that depends only on $n$, $d$, 
and $M$, such that
$$
K^{-1} |x-y| \leq |g(x)-g(y)| \leq K |x-y|
\hbox{ for } x,y \in \R^n.
\leqno (12.34)
$$

\ms
Notice that Theorem 1.21 is a special case of this.
Because of Theorem 2.23, Corollary~12.33 will follow as soon as
we check that the sufficient condition (2.24) is satisfied when $J$ 
is bounded on $E_0$. That is, it will be enough to check that
$$
\sum_{k \geq 1} \varepsilon'_k(f_k(z))^2 \leq C M
\leqno (12.35)
$$
for $z \in \Sigma_0$, where the $\varepsilon'_k$ 
are defined by (2.21) and (2.22). We dropped $k=0$
from (12.35) because $\varepsilon'_0(z) \leq 1$, so it 
will not alter the boundedness of the sum.

So we let $z\in \Sigma_0$ be given, and set $z_k = f_k(z)$ 
for $k \geq 1$. Recall from (2.21) and (2.22) that
$$\eqalign{
\varepsilon'_k(z_k) = \sup\big\{ &
d_{x_{i,l},100r_{l}}(P_{j,k},P_{i,l}) \, ; \,
j\in J_k, \, l \in \{ k-1, k \}, 
\cr& \hskip 3.5cm
i\in J_{l}, \hbox{ and } 
z_k \in 10 B_{j,k} \cap 11 B_{i,l} \big\}  
}\leqno (12.36) 
$$
for $k \geq 1$, with the convention that
$\varepsilon'_k(z_k) = 0$ if $z_k\in \R^n \sm V_k^{10}$.
Recall also that we set $P_{j,k} = P_k(x_{j,k})$ 
for $k\geq 0$ and $j\in J_k$.

Choose $x\in E_0$ such that
$$
|x-f(z)| \leq 2\dist(f(z),E_0).
\leqno (12.37) 
$$
We claim that 
$$
d_{x_{i,l},100r_{l}}(P_{j,k},P_{i,l}) 
\leq C (\alpha_k(x)+\alpha_{k-1}(x))  
\leqno (12.38) 
$$
when $k \geq 1$, $j\in J_k$, and $l \in \{ k-1, k \}$
are such that $z_k \in 10 B_{j,k} \cap 11 B_{i,l}$.
First observe that 
$$\leqalignno{
|x_{i,l}-x| &\leq |x_{i,l}-z_k|+ |z_k-f(z)|+|f(z)-x|
\leq 11r_l + C \varepsilon r_l + 2\dist(f(z),E_0)
\cr&\leq 11r_l + C \varepsilon r_l + 2\dist(z_k,E_0)
\leq 11r_l + C \varepsilon r_l + 2|x_{i,l}-z_k|
\leq 34 r_l
&(12.39) 
}%\leqno (12.39) 
$$
by (6.9), because $z_k \in \cap 11 B_{i,l}$,
and because $x_{i,l} \in E_0$ by construction.
Similarly, $|x_{i,k}-x| \leq 34 r_k$. Notice that
$$
d_{x_{i,l},100r_{l}}(P_{j,k},P_{i,l}) \leq C (d_1+d_2+d_3),
\leqno (12.40) 
$$
with
$$
d_1 = d_{x,150r_{l}}(P_{j,k},P_{k}(x)),
\ d_2 = d_{x,150r_{l}}(P_{k}(x),P_{l}(x)),
\ d_3 = d_{x,150r_{l}}(P_{l}(x),P_{i,l}),
\leqno (12.41) 
$$
and where we may even drop the middle term $d_2$
when $l=k$. Now 
$$
d_1 \leq C d_{x,r_{k}}(P_{j,k},P_{k}(x)) \leq C \alpha_l(x)
\leqno (12.42) 
$$
by elementary geometry because we are dealing with $d$-planes, 
and by (12.31). Similarly, $d_3 \leq  C \alpha_k(x)$
and (if $l=k-1$) $d_2 \leq C \alpha_{k-1}(x)$. 

So (12.38) holds, and hence $\varepsilon'_k(z_k) \leq
C (\alpha_k(x)+\alpha_{k-1}(x))$ because of (12.36).
We sum over $k$ and get that
$$
\sum_{k \geq 1} \varepsilon'_k(z_k)^2 
\leq C \sum_{k \geq 0} \alpha_k(x)^2 \leq C M,
\leqno (12.43)
$$
by (12.32) and as needed for (12.35). Corollary 12.33 follows.
\qed

\ms
A slightly unpleasant feature of $J$ in (12.32) and 
Corollary 12.33 is that they depend on a choice of
planes $P_k(x)$. In the context of Theorem 12.18 (when
we study Reifenberg-flat sets with holes),
Counterexample 12.28 shows that this is probably part of the
inherent difficulties of the problem. 

In the context of Theorem 1.1 (when we have bilateral 
approximation by planes), this is less of an issue. 
The next result shows that if we are ready to use a slightly 
stronger necessary condition, then any reasonable
choice of planes will work.

We now give a sufficient condition for the boundedness of $J$
in terms of $\beta_\infty$-numbers. Let $\beta_\infty(x,r)$
and $J_\infty$ be as in (1.11) and (1.12).

\ms\proclaim Corollary 12.44. Let $\varepsilon$, $E$, $U$,
and $\Sigma_0$ satisfy the hypotheses of Theorem 12.3.
Suppose in addition that for some $M \geq 1$,
$$
J_\infty(x) =: \sum_{k \geq 0} \beta_\infty(x,r_k)^2 \leq M  
\ \hbox{ for }
x \in E_0 = \big\{ x\in E \, ; \, \dist(x,U) \leq 1 \big\},
\leqno (12.45)
$$
where the $\beta_\infty(x,r)$ are defined in (1.11).
Then (we can choose the $P_k(x)$ above so that) the function 
$g$ of Theorem 12.3 is bi-Lipschitz: we can find
$K = K(n,d,M)$ such that (12.34) holds.

\ms
So we can replace $J$ from (12.32) with the more explicit
$J_\infty(x)$.

Our condition (12.45) is not necessary and sufficient, but 
the square exponent is right, and (12.45) is not so far off.
See the discussion for Ahlfors-regular sets, where it will 
appear that a BMO-like condition on $L^q$ variants of the
$\beta_\infty(x,r_k)$ is needed. Also recall that functions
like $J_\infty$ were introduced by 
P. Jones and C Bishop %%% and Bishop?
and have been widely used in various contexts involving 
parameterizations of sets. %%% add some references? No more?

Notice that we are only using a one-sided $\beta_\infty$ function
here, but we shall rely on the fact that the two-sided version
$\overline \beta_\infty(x,r_k) = \inf_{P} d_{x,r}(E,P)$ stays
small (by (12.5)) to control the variations of the $P_k(x)$.

As for the previous results, Corollary 12.44 is an extension of
Theorem 1.13.

\ms
Let us now deduce Corollary 12.44 from Corollary 12.33.
We first need to choose $d$-planes $P_k(x)$
for $x\in E_0$ and $k \geq 0$. When
$$
k \geq 3 \hbox{ and } \beta_\infty(x,120r_{k}) \leq \varepsilon,
\leqno (12.46)
$$
we choose $P_k(x)$ to be any $d$-plane through $x$ so that
$$
\dist(y,P_k(x)) \leq 10^3 \beta_\infty(x,120r_{k}) r_{k}
\leq 10^3 \varepsilon r_{k}
\ \hbox{ for } y\in E\cap B(x,120r_k).
\leqno (12.47)
$$
Otherwise, if (12.46) fails, we simply choose
$P_k(x) = P(x,r_k)$, where $P(x,r_k)$ comes from
our assumption (12.5). Let us check that even when (12.46)
holds, we have that 
$$
d_{x,200r_{k}}(P_k(x),P(x,r_k))  \leq C \varepsilon.
\leqno (12.48)
$$
Apply Lemma 12.12 to $P_1 = P(x,r_k)$, $P_2 = P_k(x)$,
$z=x$, $r= r_k$, and $\tau = C\varepsilon$. 
The assumption (12.13) is satisfied because if
$y\in P(x,r_k) \cap B(x,r_k)$, (12.5) says that
we can find $\xi\in E$ such that $|\xi-y| \leq \varepsilon r_k$,
and then $\dist(\xi,P_k(x)) \leq 10^3 \varepsilon r_{k}$
by (12.47). The conclusion of Lemma 12.12 is exactly (12.48).

Because of (12.48) and (12.5), we also have that
$$\eqalign{
&\dist(y,P_k(x)) \leq C\varepsilon r_k 
\hbox{ for } y\in E\cap B(x,110r_k)
\cr& \hskip 3cm 
\hbox{ and }
\dist(y,E) \leq C\varepsilon r_k 
\hbox{ for } y\in P_k(x)\cap B(x,110r_k).
}\leqno (12.49)
$$
That is, the planes $P_k(x)$ also satisfy 
the property (12.5) (although with the slightly larger
constant $C\varepsilon$), which means that we could choose
them in the statement of Corollary 12.33.

All we have to do now is check that for this choice 
of planes, $J(x)$ is bounded on $E_0$, and then Corollary 12.33
will give the result. We shall prove that
$$
\alpha_k(x) \leq C \beta_\infty(x,r_{k-3}) 
\leqno (12.50)
$$
for $x\in E_0$ and $k \geq 3$; 
this and the definition (12.32) will then imply that 
$$
J(x) \leq C J_\infty (x) + \sum_{k < 3} \alpha_k(x) 
\leq C J_\infty (x) + 4
\leqno (12.51)
$$
and Corollary 12.44 will follow because we assume that
$J_\infty$ is bounded on $E_0$.

Let $x\in E_0$ and $k \geq 3$ be given; in view of
the definition (12.31) of $\alpha_k(x)$, we need to check 
that 
$$
d_{x,r_k}(P_{k+1}(x),P_k(x)) \leq C \beta_\infty(x,r_{k-3}) 
\leqno (12.52)
$$
and
$$
d_{x,r_k}(P_{k}(x),P_k(y)) \leq C \beta_\infty(x,r_{k-3}) 
\leqno (12.53)
$$
for $y\in E_0 \cap B(x,35r_k)$, and then (12.50) will
follow. We shall only prove (12.53), as (12.52) is simpler.

Choose an orthonormal basis $e_1, \cdots, e_d$ of
the vector space parallel to $P_k(x)$, and set
$p_0 = x$ and $p_l = p_0 + r_k e_l$ for $1 \leq l \leq d$.
By (12.49), we can choose points $x_l \in E$ such that 
$$
|x_l-p_l| \leq C \varepsilon r_k
\ \hbox{ for } 0  \leq l \leq d
\leqno (12.54)
$$
(and we may even take $x_0=x$). We shall first assume
that both $P_k(x)$ and $P_k(y)$ were chosen according
to (12.47). Observe that for 
$0 \leq l \leq d$, 
$$
|x_l-y| \leq |x_l-x| + 35r_k \leq 37r_k,
\leqno (12.55)
$$
so $x_l \in B(y,120r_k)$ and
$$
\dist(x_l,P_k(y)) 
\leq 10^3 \beta_\infty(y,120r_{k}) \, r_k
\leqno (12.56)
$$
by (12.47). Similarly, 
$$
\dist(x_l,P_k(x)) 
\leq 10^3 \beta_\infty(x,120r_{k}) \, r_k
\leqno (12.57)
$$
by the last part of (12.55) and (12.47), so we can find
$\xi_l \in P_k(x)$ and $\zeta_l \in P_k(y)$ such that
$$\eqalign{
|\xi_l-x_l|+|\zeta_l-x_l| 
&\leq \dist(x_l,P_k(y)) + \dist(x_l,P_k(x))
\cr&
\leq 10^3 \beta_\infty(y,120r_{k}) \, r_k
+ 10^3 \beta_\infty(x,120r_{k}) \, r_k.
}\leqno (12.58)
$$
Let us also check that
$$
\beta_\infty(y,120r_{k}) \leq 20 \beta_\infty(x,r_{k-3}).
\leqno (12.59)
$$
Let $P$ be a plane through $x$ such that 
$\dist(w,P) \leq \beta_\infty(x,r_{k-3}) r_{k-3}$
for $w \in E \cap B(x,r_{k-3})$ then in particular 
$\dist(y,P) \leq \beta_\infty(x,r_{k-3}) r_{k-3}$.
Let $P'$ be the translation of $P$ that goes through
$y$; then
$\dist(w,P') \leq 2 \beta_\infty(x,r_{k-3}) r_{k-3}$
for $w \in E \cap B(y,120r_{k}) \i B(x,r_{k-3})$,
and (12.59) follows. The same proof, without any need for
the translation, shows that $\beta_\infty(x,120r_{k}) 
\leq 10 \beta_\infty(x,r_{k-3})$, and so (12.58) yields
$$
|\xi_l-x_l|+|\zeta_l-x_l| 
\leq C  \beta_\infty(x,r_{k-3}) \, r_k.
\leqno (12.60)
$$
This holds when $P_k(x)$ and $P_k(y)$ were chosen according
to (12.47). Otherwise, we simply use (12.49) (for $x$ and 
for $y$) to choose $\xi_l \in P_k(x)$ and $\zeta_l \in P_k(y)$ 
such that
$$
|\xi_l-x_l|+|\zeta_l-x_l| \leq C \varepsilon r_k.
\leqno (12.61)
$$
If (12.46) fails for $y$, then 
$\beta_\infty(x,r_{k-3}) \geq \beta_\infty(y,120r_{k})/20
\geq \varepsilon/20$ by (12.59), and (12.60) holds too.
If (12.46) fails for $x$, we even get that
$\beta_\infty(x,r_{k-3}) \geq \beta_\infty(x,120r_{k})/10
\geq \varepsilon/10$. So (12.60) holds in all cases.

We shall conclude with the following lemma.

\ms\proclaim Lemma 12.62.
Let $z \in \R^n$, $r > 0$, $\tau \in (0,10^{-1})$,
two affine $d$-planes $P_1$ and $P_2$, and $d$ 
mutually orthogonal unit vectors $e_1, \cdots, e_d$ be given.
Suppose that, for $0 \leq l \leq d$, we are given 
points $\xi_l \in P_1$ and $\zeta_l \in P_2$,
so that $\xi_0 \in B(x,r)$, 
$$
|\xi_l-\zeta_l| \leq \tau r
\ \hbox{ for } 0 \leq l \leq d,
\leqno (12.63)
$$
and 
$$
|\xi_l-\xi_0-re_l| \leq r/10
\ \hbox{ for } 1 \leq l \leq d.
\leqno (12.64)
$$
Then
$$
d_{z,\rho}(P_1,P_2) \leq C \tau 
\ \hbox{ for } r \leq \rho \leq 10^4r.
\leqno (12.65) 
$$

We leave the proof to the reader, but claim that
since we can immediately reduce to the case when
$P_1 = \R^d$ and use coordinates, it would be easy 
to verify.
\qed

\ms
Let us apply Lemma 12.62 with $r=r_k$, $z=p_0=x$, $P_1 = P_k(x)$,
$P_2 = P_k(y)$, and $\tau = C\beta_\infty(x,r_{k-3})$; (12.63)
follows from (12.60) and (12.64) holds by (12.54) and
(12.58) (or (12.61) if the right-hand side
of (12.58) is larger than $C\varepsilon$).
Now (12.65)  says that $d_{x,r_k}(P_k(x),P_k(y)) 
\leq C \beta_\infty(x,r_{k-3})$, which is (12.53). 

As we said earlier, (12.52) is easier, (12.50) follows
from (12.52) and (12.53), and Corollary 12.44 follows 
from (12.50).
\qed

\bsi
{\bf 13. Local lower-Ahlfors regularity and a better sufficient
bi-Lipschitz condition}
\ms

The next sections will be devoted to locally 
Reifenberg-flat Ahlfors-regular sets. In most of this one,
and to the authors slight surprise, we do not need to assume 
that $E$ is locally Ahlfors-regular yet. We need and 
prove the lower bound (see (13.2)), and the results would probably 
be hard to apply when $H^d$ restricted to $E$ is locally too large.

The main result of this section is that we can replace the 
Jones function $J_\infty$ in Corollary 12.44 with the 
often smaller $J_1$ based on $L^1$ norms. 
See Corollary 13.4. 

When $n=d+1$, we shall use Corollary 13.4
to give another sufficient condition for the existence of 
a bi-Lipschitz parameterization of $E$, in terms of the unit 
normal to $E$. This condition is
reminiscent of conditions given by Semmes in the context
of Chord-Arc Surfaces with Small Constants [Se1,2,3]. %% CASSC
See Corollary 13.46.

Let $E \i \R^n$ be given, and set
$$
\beta_1(x,r) = \inf_{P} \Big\{
\, {1 \over r^d} \int_{y \in E \cap B(x,r)} 
{\dist(y,P) \over r} \, dH^d(y) \Big\},
\leqno (13.1)
$$
for $x\in \R^n$ and $r >0$ (as in (1.15)),
where the infimum is taken over all $d$-planes $P$ 
through $B(x,r)$ (there is no point in taking $P$ further away)
and $H^d$ denotes the $d$-dimensional Hausdorff measure
(see [Ma] or [Fe], and %% les 2 livres
recall that $H^d$ coincides with the surface measure on
smooth $d$-dimensional submanifolds).

Even though we do not need to assume this to define 
$\beta_1(x,r)$ and prove the result below, it is often %v
easier to use $\beta_1$ when
$$
C_0^{-1} r^d \leq H^d(E\cap B(x,r)) \leq C_0 r^d
\leqno (13.2)
$$
(as we shall assume in the next sections). Note that
if (the second half of) (13.2) holds, then we can
deduce from H\"older's inequality and the definition 
(1.15) that
$$
C_0^{q-p\over pq} \beta_q(x,r) \leq  \, \beta_p(x,r) 
\leq C_0^{1 \over p} \,\beta_\infty(x,r)
\ \hbox{ for } 1 \leq q < p < +\infty,
\leqno (13.3)
$$ 
and $\beta_1(x,r)$ is easier to control than the
other $\beta_q$. If instead $H^d(E\cap B(x,r))$ is too large
(and in particular if it is infinite), we shall probably
not be able to estimate $\beta_1(x,r)$, and we may as well 
use $\beta_\infty$ as in the previous section.
Of course we may try to replace $H^d$ in (13.2) with a different 
measure, or normalize differently, but the choice of (13.1) 
and (13.2) seems very reasonable in the present context.

\ms\proclaim Corollary 13.4. 
Let $\varepsilon$, $E$, $U$, and $\Sigma_0$ satisfy 
the hypotheses of Theorem 12.3.
Suppose in addition that for some $M \geq 1$,
$$
J_1(x) =: \sum_{k \geq 3} \beta_1(x,r_k)^2 \leq M  
\ \hbox{ for } %vv \geq 3 to make 13.49 simpler
x \in E_0 = \big\{ x\in E \, ; \, \dist(x,U) \leq 1 \big\}.
\leqno (13.5)
$$
Then we can choose the $d$-planes $P_{j,k}$ in the
proof of Theorem 12.3 so that, in addition to the
properties stated in Theorem 12.3, $g$ is $K$-bi-Lipschitz
on $\R^n$ (i.e., (12.34) holds), with $K = K(n,d,M)$.

\ms
There is not much point in using Corollary 13.4 when
we have no control on $H^d(E)$, but when
(13.2) holds for $x \in E \cap U^+$ and $r \leq 1$,
(13.3) says that Corollary 13.4 is better than
its analogue with $J_\infty$ (Corollary 12.44).
In fact, (13.5) not far from being optimal:
as we shall see in Section 15, if $E$ is a bi-Lipschitz 
image of $\Sigma_0$, then we have BMO-type estimates for the 
$J_q$ for $1 \leq q < {2d \over d-2}$. And $J_\infty$
could fail to be integrable, essentially because
we have no Sobolev embedding for large exponents.
See Remark 15.38. 

\ms
We shall try to use the same sort of proof as
for Corollary 12.44, but first we shall establish lower bounds
on $H^d(E\cap B(x,r))$, which are obviously needed if we
want the $\beta_1(x,r)$ to give some control on the geometry.

\ms\proclaim Lemma 13.6. 
Let $\varepsilon$, $E$, $U$, and $\Sigma_0$ 
be as in Theorem 12.3. Then 
$$
H^d(\overline E\cap B(x,r_0)) 
\geq (1-C \varepsilon)\, \omega_d \, r_0^d
%\hbox{ for $x_0\in E_0$ and } 0 < r_0 < 10^{-1}, 
\leqno (13.7)
$$
for $x_0 \in E_1 = \big\{ x\in E \, ; \, 
\dist(x,U) \leq 3/2 \big\}$ and $0 < r_0 \leq 10^{-1}$, 
and where $\omega_d = H^d(\R^d \cap B(0,1))$ denotes
the measure of the unit ball in $\R^d$.

\ms
Let $x_0\in E_1$ and $0 < r_0 \leq 10^{-1}$ be given.
We shall only need to know that for 
$x\in E \cap B(x_0,5r_0)$ and $0 < r \leq r_0$,
we can find a $d$-plane $P = P(x,r)$ such that
(12.5) holds. This follows from the assumptions
of Theorem 12.3, because $x\in U^+ = \big\{ x\in E \, ; \, 
\dist(x,U) \leq 2 \big\}$. 

It will be more convenient to renormalize
and work with $B(x_0,r_0)$ replaced by $B(0,10)$. So we set
$$
F = 10 r_0^{-1} (E-x_0).
\leqno (13.8)
$$ 
By what was just said, for each $x\in F \cap B(0,50)$ 
and $0 < r \leq 10$, there is a $d$-plane $P(x,r)$ such that 
$$
d_{x,110r}(P(x,r),F) \leq \varepsilon/110,
\leqno (13.9)
$$
as in (12.5).

Then $F$ satisfies the assumptions of Theorem 12.3,
with $\Sigma_0 = P(0,1)$ and $U = B(0,40)$.
Indeed, $\Sigma_0 = P(0,10)$ satisfies (12.1) trivially,
(12.4) follows from (13.9), and so does (12.5).
We do not
 even need to multiply $\varepsilon$ by a constant.

So Theorem 12.3 gives a bi-H\"{o}lder mapping $g: \R^n \to \R^n$ 
such that in particular
$$
|g(z)-z| \leq C \varepsilon 
\ \hbox{ for } z \in \R^n
\leqno (13.10)
$$
and 
$$
\overline F \cap B(x,40) = g(P(0,1)) \cap B(x,40).
\leqno (13.11)
$$
as in (12.7) and (12.9).

Denote by $\pi$ the orthogonal projection onto $P=P(0,1)$,
and set $h = \pi \circ g$. Note that for $z\in P$,
$$
|h(z)-z| \leq |h(z)-g(z)|+|g(z)-z|
\leq \dist(g(z),P) + |g(z)-z| \leq 2 |g(z)-z|
\leq 2C \varepsilon 
\leqno (13.12)
$$
by (13.10) and because $z\in P$. From this and a
little bit of degree theory, we deduce that
$h(P)$ contains $P \cap B(0,10)$. The proof is the same
as for (5.47); we could also say that a continuous
mapping $h : P \to P$ such that $h(x)=x$ for $x$ large
(which is the case here by construction) is surjective.

Set $D = P \cap B(0,10-3C\varepsilon)$, where $C$ is as in
(13.10) and (13.12). Let $w\in D$, and let $z\in P$ 
be such that $h(z)=w$. Thus $w \in B(0,10-C\varepsilon)$
by (13.12), and $g(w) \in B(0,10)$ by (13.10). 
In addition, $g(w) \in \overline F$ by (13.11),
and $w = h(z) = \pi(g(w))$.
So $D \i \pi(\overline F \cap B(0,10))$, and
$$
H^d(\overline F \cap B(0,10))
\geq H^d(\pi(\overline F \cap B(0,10)))
\geq H^d(D) = \omega_d (10-3C\varepsilon)^d
\leqno (13.13)
$$
because $\pi$ is $1$-Lipschitz. By (13.8),
$$
H^d(\overline E \cap B(x_0,r_0)) 
= (r_0/10)^d H^d(\overline F \cap B(0,10))
\geq \omega_d (1-3C\varepsilon/10)^d r_0^d,
\leqno (13.14)
$$
as needed for (13.7). Lemma 13.6 follows.
\qed

\ms
Now we want to follow the proof of Corollary 12.44. 
We shall have to be slightly more careful
about the choice of points $x_{j,k}$. 
Here we will
 pay the price for deciding that the $d$-plane 
$P_{j,k}$ should go through $x_{j,k}$, because it does
not make sense to force the planes $P$ to go through
the center of $B(x,r)$ when we define $\beta_1(x,r)$.
That is, the best approximating planes $P$ could well
pass some distance away from $x$.

For each $k \geq 0$, we start with a collection 
$\{ \widetilde x_{j,k} \}$, $j\in J_k$, of points of $E_0$,
which is maximal under the constraint that
$|\widetilde x_{i,k}-\widetilde x_{j,k}| \geq 4r_k/3$
when $i \neq j$, and we promise to choose 
$$
x_{j,k} \in E \cap B(\widetilde x_{j,k},r_k/3).
\leqno (13.15)
$$
Then
$$
E_0 \i \bigcup_{j\in J_k} \overline B(\widetilde x_{j,k},4r_k/3)
\i \bigcup_{j\in J_k} B(x_{j,k},5r_k/3)
\leqno (13.16)
$$
and (2.3) follows as before because $\dist(x_{i,k+1},E_0) \leq
r_{k+1}/3$ for $i \in J_{k+1}$. We still get that the limit
set $E_\infty$ of (2.19) is the closure of $E_0$.

Also notice that $x_{j,k} \in E \cap U^+$ because
$|x_{j,k}-\widetilde x_{j,k}| < r_k/3$ and 
$\widetilde x_{j,k} \in E_0$  (see (12.10) and (12.2))
so $\dist(x_{j,0},\Sigma_0) \leq \varepsilon$
for $j \in J_0$, by (12.4).
That is, (2.7) holds.

Next we want to choose the $d$-planes $P_{j,k}$, $j\in J_k$. 
We start in the most interesting case when
$$
k \geq 2 \hbox{ and } \beta_1(\widetilde x_{j,k},120r_{k}) 
\leq \varepsilon.
\leqno (13.17)
$$
We choose a first $d$-plane $P'_{j,k}$ such that
$$
(120r_{k})^{-d} \int_{y \in E \cap B(\widetilde x_{j,k}, 120r_{k})} 
{\dist(y,P'_{j,k}) \over 120r_{k}}  \, dH^d(y)
\leq 2 \beta_1(\widetilde x_{j,k},120r_{k}) \leq 2 \varepsilon
\leqno (13.18)
$$
(compare with the definition (13.1)),
and then use Chebyshev's inequality to choose 
$x_{j,k} \in E\cap B(\widetilde x_{j,k},r_k/3)$ so that
$$\eqalign{
\dist(x_{j,k},P'_{j,k}) 
&\leq H^d(E\cap B(\widetilde x_{j,k},r_k/3))^{-1}
\int_{E\cap B(\widetilde x_{j,k},r_k/3)} \dist(y,P'_{j,k}) \, dH^d(y)
\cr&
\leq C r_{k}^{-d}
\int_{E \cap B(\widetilde x_{j,k},120r_k)} \dist(y,P'_{j,k}) \, dH^d(y)
\cr&
\leq C \beta_1(\widetilde x_{j,k},120r_{k}) r_k
\leq C \varepsilon r_k
}\leqno (13.19)
$$
by (13.7), (13.18), and (13.17). Now we let $P_{j,k}$ be the
$d$-plane parallel to $P'_{j,k}$ that contains $x_{j,k}$.

In the other case, when (13.17) fails, we simply
take $x_{j,k} = \widetilde x_{j,k}$ and 
$P_{j,k} = P'_{j,k} = P(x_{j,k},r_k)$ (the $d$-plane
given by (12.5)).

We shall prove directly that we have the
compatibility conditions (2.8)-(2.10) and the
summability condition (2.24) (the adaptations 
that would be needed to follow the proof of Corollary 12.44
would be even more painful).

Fix $k \geq 0$ and $j \in J_k$.
We shall now choose $d+1$ points $z_l \in E \cap B(x,r_k)$,
$0 \leq l \leq d$, that will control the position of
$P'_{j,k}$. First choose an orthonormal basis 
$\{ e_1, \cdots e_d \}$ of the vector space parallel to 
$P(x_{j,k},r_k)$. Set
$$
p_0 = \widetilde x_{j,k}
\hbox{ and } p_l = p_0 + {1 \over 2} \, r_k e_l  
\hbox{ for } 1 \leq l \leq d.
\leqno (13.20)
$$
By definition of $P(x_{j,k},r_k)$ (i.e., (12.5)),
we can find $w_l \in E \cap B(p_l,C\varepsilon r_k)$.
Observe that
$$
|w_l - x_{j,k}|
< |w_l - \widetilde x_{j,k}|+ {r_k \over 3} 
\leq |p_l - \widetilde x_{j,k}|+C\varepsilon r_k + {r_k \over 3} 
\leq {5 r_k \over 6}  + C \varepsilon r_k,
\leqno (13.21)
$$
by (13.15) and (13.20), so 
$$
B(w_l,r_{k+2}) \i B(x_{j,k},r_k).
\leqno (13.22)
$$

First assume that (13.17) holds. Observe that %v
$|w_l-\widetilde x_{j,k}| \leq 10^{-1}$ 
by the end of (13.21) and because $k \geq 2$, so
$w_l \in E_1$ because $\widetilde x_{j,k} \i E_0 \,$, and so
and  $H^d(E\cap B(w_l,r_{k+2})) \geq C^{-1} r_{k}^d$
by Lemma~13.6. We use Chebyshev's inequality 
to find $z_l \in E \cap B(w_l,r_{k+2})$ such that
$$\leqalignno{
\dist(z_l,P'_{j,k}) 
&\leq H^d(E\cap B(w_l,r_{k+2}))^{-1}
\int_{E\cap B(w_l,r_{k+2})} \dist(y,P'_{j,k}) \, dH^d(y)
\cr&
\leq C r_{k}^{-d}
\int_{E \cap B(\widetilde x_{j,k},120r_k)} \dist(y,P'_{j,k}) \, dH^d(y)
%\cr&
\leq C \beta_1(\widetilde x_{j,k},120r_{k}) r_k
\leq C \varepsilon r_k
& (13.23)
}%\leqno (13.23)
$$
because $B(w_l,r_{k+2}) \i B(\widetilde x_{j,k},120r_k)$ (by (13.22)), 
the choice of $x_{j,k}$, and by (13.18). %v
Let us record that
$$
|z_l - z_0 - {r_k \over 2}\, e_l|
= |z_l - z_0 - p_l + p_0| 
\leq 2 r_{k+2} + |w_l - w_0 - p_l + p_0| < 3r_{k+2}
\leqno (13.24)
$$
by (13.20) and because $w_l \in B(p_l,C\varepsilon r_k)$.
We shall later need to know that
$$
d_{\widetilde x_{j,k},10r_k}(P_{j,k},P(\widetilde x_{j,k},r_k))
\leq C \varepsilon.
\leqno (13.25)
$$
Let us apply Lemma 12.62 with $z=\widetilde x_{j,k}$, 
$r = {r_k \over 2}$, $\tau = C \varepsilon$,
$P_1 = P'_{j,k}$, and $P_2 = P(\widetilde x_{j,k},r_k)$.
We use (13.23) to find $\xi_l \in P'_{j,k}$ such that
$|\xi_l - z_l| \leq C \varepsilon r_k$
and use (12.5) and (13.22)
to choose $\zeta_l \in P(\widetilde x_{j,k},r_k)$ such
that $|\zeta_l - z_l| \leq C \varepsilon r_k$ (so that
(12.63) holds). Note that (12.64) follows from (13.24),
so Lemma 12.62 applies and (12.65) (with $\rho = 11 r_k$) says that
$d_{\widetilde x_{j,k},11r_k}(P'_{j,k},P(\widetilde x_{j,k},r_k))
\leq C \varepsilon$; (13.25) follows because $P_{j,k}$
is obtained from $P'_{j,k}$ by a translation by less
than $C \varepsilon r_k$ (by (13.19) and the line that follows it).

Return to the choice of $z_l$.
When (13.17) fails, we simply take $z_l = w_l$.
Recall that $w_l\in E\cap B(p_l,C\varepsilon r_k)$, so that
in this case,
$$\eqalign{
\dist(z_l,P'_{j,k}) =  \dist(w_l,P(x_{j,k},r_k))
\leq |w_l-p_l| \leq C \varepsilon r_k
}\leqno (13.26)
$$
because $P'_{j,k} = P(x_{j,k},r_k)$ and $p_l \in P(x_{j,k},r_k)$.

\ms
Let us first check (2.9). Let $j\in J_0$ and
$x\in \Sigma_0$ be such that $|x_{j,0}-x| \leq 2$.
Notice that (13.17) fails because $k=0$,
so $x_{j,0} = \widetilde x_{j,0}$ and 
$P'_{j,0} = P(x_{j,0},1)$.

Choose the $z_l$, $0 \leq l \leq d$ as above.
Observe that $z_l \in E \cap U^+$ because
$z_l \in B(x_{j,k},r_k)$ by (13.22) and
$x_{j,0} = \widetilde x_{j,0} \in E_0$
(also see the definitions (12.10) and (12.2)).
So $\dist(z_l,\Sigma_0) \leq \varepsilon$ by (12.4),
and $\dist(z_l,P_x) \leq 2\varepsilon$
by (2.4), (2.5), and also because $|x_{j,0}-x| \leq 2$
and $z_l = w_l \in B(x_{j,k},r_k)$ by (13.22).

We may now apply Lemma 12.62 with $z=x_{j,0}$, 
$r = {1 \over 2}$, $\tau = C \varepsilon$, $P_1 = P_{j,0}$, 
$P_2 = P_x$, and where $\xi_l \in P_{j,0} = P(x_{j,0},1)$ 
is chosen such that $|\xi_l - z_l| \leq C \varepsilon$ 
(using (12.5)) and $\zeta_l \in P_x$ is chosen such
that $|\zeta_l - z_l| \leq C \varepsilon$. Then
(12.63) holds, and (12.64) follows from (13.24).
So Lemma 12.62 applies, and says that 
$d_{x_{j,0},100}(P_{j,0},P_x) \leq C \varepsilon$,
as needed for (2.9).

We shall try to prove (2.8), (2.10), and (2.24) at the
same time, so let us fix $k \geq 0$ and $j\in J_k$, and 
give ourselves $m \in \{ k, k-1 \}$ and $i \in J_m$,
such that 
$$
|x_{j,k}-x_{i,m}| \leq 100 r_m.
\leqno (13.27)
$$
We want to show that $P_{i,m}$ lies close to $P_{j,k}$.
Let us  first assume that (13.17) holds for both pairs
$(j,k)$ and $(i,m)$. As before, we use Chebyshev to find
$z_l \in E \cap B(w_l,r_{k+2})$ such that
$$
\dist(z_l,P'_{j,k}) 
\leq C \beta_1(\widetilde x_{j,k},120r_{k}) r_k
\leq C \varepsilon r_k
\leqno (13.28)
$$
as in (13.23), but at the same time we use the
fact that 
$$
B(w_l,r_{k+2}) \i B(\widetilde x_{i,m},110r_{m})
\leqno (13.29)
$$
(because $|w_l-\widetilde x_{i,m}| \leq |w_l-x_{j,k}|
+|x_{j,k}-x_{i,m}|+|x_{i,m}-\widetilde x_{i,m}|
\leq {5 r_k \over 6}  + C \varepsilon r_k + 100 r_m
+ {r_m \over 3} < 102r_m$ by (13.21), (13.27), and (13.15)) 
to demand that in addition
$$\leqalignno{
\dist(z_l,P'_{i,m}) 
&\leq 2 H^d(E\cap B(w_l,r_{k+2}))^{-1}
\int_{E\cap B(w_l,r_{k+2})} \dist(y,P'_{i,m}) \, dH^d(y)
\cr&
\leq C r_{k+2}^{-d}
\int_{E \cap B(\widetilde x_{i,m},120r_k)} \dist(y,P'_{i,m}) \, dH^d(y)
%\cr&
\leq C  \beta_1(\widetilde x_{i,m},120r_{m}) r_k
\leq C \varepsilon r_k
& (13.30)
}%\leqno (13.30)
$$
by (13.18) and (13.17) for the pair $(i,m)$. 
Recall from (13.19) and the line below
that $P_{j,k}$ is obtained from $P'_{j,k}$
by a translation of $\dist(x_{j,k},P'_{j,k}) 
\leq C \beta_1(\widetilde x_{j,k},120r_{k}) r_k$,
and similarly for $P_{i,m}$. So (13.28) and (13.30)
allow us to find $\xi_l \in P_{j,k}$ and
$\zeta_l \in P_{i,m}$ such that
$$
|\xi_l-z_l| + |\zeta_l-z_l|
\leq  C r_k [\beta_1(\widetilde x_{j,k},120r_{k})
+ \beta_1(\widetilde x_{i,m},120r_{m})]
\leq C \varepsilon r_k.
\leqno (13.31)
$$
Apply Lemma 12.62 with $z=x_{j,k}$, $r = {r_k \over 2}$, 
$\tau = C [\beta_1(\widetilde x_{j,k},120r_{k})
+ \beta_1(\widetilde x_{i,m},120r_{m})]$,
$P_1 = P_{j,k}$, and $P_2 = P_{i,m}$;
(12.63) holds by (13.31) and (12.64) follows again
from (13.24). So Lemma 12.62 applies, and says that 
$$
d_{x_{j,k},\rho}(P_{j,k},P_{i,m}) 
\leq C [\beta_1(\widetilde x_{j,k},120r_{k})
+ \beta_1(\widetilde x_{i,m},120r_{m})]
\leq C \varepsilon 
\leqno (13.32)
$$
for ${r_k \over 2} \leq \rho \leq 5000 r_k$.

If (13.17) holds for $(j,k)$, but fails for $(i,m)$,
we simply select the $z_l \in E \cap B(w_l,r_{k+2})$
so that (13.28) holds. But then 
$x_{i,m} = \widetilde x_{i,m} \in E_0$
and $P_{i,m} = P(x_{i,m},r_m)$, so (12.5) and (13.29)
say that we can choose $\zeta_l \in P_{i,m}$ such that
$|\zeta_l-z_l| \leq \varepsilon r_m$. Then 
$$
|\xi_l-z_l| + |\zeta_l-z_l| \leq C \varepsilon r_k.
\leqno (13.33)
$$
by this and (13.28). Lemma 12.62 yields
$$
d_{x_{j,k},\rho}(P_{j,k},P_{i,m}) \leq C \varepsilon
\ \hbox{ for } {r_k \over 2} \leq \rho \leq 5000 r_k.
\leqno (13.34)
$$

Similarly, if (13.17) fails for $(j,k)$ and holds for $(i,m)$,
we choose the $z_l$ so that (13.30) holds, and use
the fact that $P_{j,k} = P(x_{j,k},r_k)$ to apply (12.5),
show that $\dist(z_l,P_{j,k}) \leq \varepsilon r_k$,
and get (13.33) and (13.34). Finally, when (13.17) fails
for both pairs $(j,k)$ and $(i,m)$, we choose $z_l = w_l$,
apply (12.5) to both pairs, and get (13.33) and (13.34).

\ms
Now (2.8) (with the new constant $C \varepsilon$)
follows from (13.32) or (13.34), applied when $m=k$ and
$\rho = 100 r_k$. For (2.10) (which we verify for $k-1$ 
and when $k\geq 1$), we choose $m = k-1$, and still get
(13.27) because $|x_{i,k-1}-x_{j,k}| \leq 2 r_{k-1}
=2 r_m$; then (2.10) also follows from (13.32) or (13.34),
applied with $\rho = 20 r_{k-1}$.

This completes the verification of the assumptions
of Theorem 2.15; we are now left with (2.24) to check.
For $z\in \Sigma_0$, choose $\overline z \in E_0$
such that
$$
|\overline z-f(z)| \leq 2 \dist(f(z),E_0)
\leqno (13.35)
$$
Let us check that
$$
\varepsilon'_k(f_k(z)) \leq C \beta_1(\overline z,r_{k-3})
\ \hbox{ for $z\in \Sigma_0$ and $k \geq 3$,}
\leqno (13.36)
$$
where the $\varepsilon'_k$ are defined in (2.21) and (2.22);
(2.24) will follow from this, because then
$$\eqalign{
\sum_{k \geq 0} \varepsilon'_k(f_k(z))^2
&\leq 3 + \sum_{k \geq 3} \varepsilon'_k(f_k(z))^2
\leq 3 + C \sum_{k \geq 3} \beta_1(\overline z,r_{k-3})^2
\cr&
\leq 3 + C J_1(\overline z)
\leq 3 + C M
}\leqno (13.37)
$$
for $z\in \Sigma_0$, by (13.5) and because $\overline z \in E_0$. 

By the definition (2.21)-(2.22), we just need to show that
$$
d_{x_{i,m},100r_m}(P_{j,k},P_{i,m}) \leq 
C \beta_1(\overline z,r_{k-3})
\leqno (13.38)
$$
when $j\in J_k$, $m \in \{ k-1,k \}$, and $i\in J_m$
are such that 
$$
y = f_k(z) \hbox{ lies in } 10B_{j,k} \cap 11B_{i,m}.
\leqno (13.39)
$$
Note that $|x_{j,k}-x_{i,m}| \leq 10 r_k + 11 r_m \leq 21 r_m$
because $10B_{j,k} \cap 11B_{i,m} \neq \emptyset$, so (13.27) holds.
If both pairs $(j,k)$ and $(i,m)$ satisfy (13.17), then 
$$\eqalign{
d_{x_{i,m},100r_m}(P_{j,k},P_{i,m})
&\leq 2 d_{x_{j,k},200r_m}(P_{j,k},P_{i,m}) 
\cr&\leq C [\beta_1(\widetilde x_{j,k},120r_{k})
+ \beta_1(\widetilde x_{i,m},120r_{m})]
}\leqno (13.40)
$$
by (13.32), applied with $\rho  = 200r_m < 5000 r_k$.
If (13.17) fails for at least one of the two pairs,
then we apply (13.34) instead of (13.32), and get that
$$
d_{x_{i,m},100r_m}(P_{j,k},P_{i,m}) \leq C \varepsilon.
\leqno (13.41)
$$
But since $k \geq 3$, if (13.17) fails for $(j,k)$, then
$\beta_1(\widetilde x_{j,k},120r_{k}) \geq \varepsilon$,
and (13.41) is stronger than (13.40). Similarly, if 
(13.17) fails for $(i,m)$, then
$\beta_1(\widetilde x_{i,m},120r_{m}) \geq \varepsilon$
and (13.41) is stronger than (13.40). So (13.40) holds
in all cases.

Recall from (13.39) that $y = f_k(z) \in 10B_{j,k}$; 
then
$$\eqalign{
|\overline z-y| &\leq |\overline z-f(z)|+|f(z)-f_k(z)| 
\leq 2 \dist(f(z)),E_0) +|f(z)-f_k(z)|
\cr&\leq 2 \dist(f_k(z)),E_0) + 3 |f(z)-f_k(z)|
\leq 2|f_k(z)-\widetilde x_{j,k}|+ C \varepsilon r_k
\cr& \leq 2|f_k(z)-x_{j,k}|+ {2r_k \over 3} + C \varepsilon r_k
\leq 21 r_k 
}\leqno (13.42)
$$
by (13.35), (6.9), because $\widetilde x_{j,k} \in E_0$,
and by (13.15) and (13.39). Because of this 
(and (13.15) and (13.39) again),
$$
B(\widetilde x_{j,k},120r_{k}) \cup
B(\widetilde x_{i,m},120r_{m})
\i B(\overline z,r_{k-3}).
\leqno (13.43)
$$
Let the $d$-plane $P$ minimize in the definition (1.15)
of $\beta_1(\overline z,r_{k-3})$; by (13.43) we can also
use it in the definition of $\beta_1(\widetilde x_{j,k},120r_{k})$,
and
$$\eqalign{
\beta_1(\widetilde x_{j,k},120r_{k})
&\leq (120r_{k})^{-d} 
\int_{y \in E \cap B(\widetilde x_{j,k},120r_{k})} 
{\dist(y,P) \over 120r_{k}} \, dH^d(y) 
\cr& \leq C (r_{k-3})^{-d}
\int_{y \in E \cap B(\overline z,r_{k-3})} 
{\dist(y,P) \over r_{k-3}} \, dH^d(y)
= C \beta_1(\overline z,r_{k-3}).
}\leqno (13.44)
$$
Similarly, $\beta_1(\widetilde x_{i,m},120r_{m})
\leq C \beta_1(\overline z,r_{k-3})$, and now (13.38)
follows from (13.40).

So (13.38) holds, which as we know implies (13.36) and 
then (2.24). This was the last verification that we needed
to do, in order to apply Theorems 2.15 and 2.23 as we did
for Theorem 12.3 and its other corollaries. So we completed
the proof of Corollary 13.4.
\qed

\ms
We shall now encode the sufficient condition of 
Corollary 13.4 in terms of unit normals. 
To simplify the statement, we assume that 
$n=d+1$ and $E = \d \Omega$ for some $C^1$ domain $\Omega$, 
so that we can define the outward unit normal $n(x)$ at 
$x\i E$, but with no estimates attached. We replace the
function $J_1$ with
$$
H(x) = \int_{0}^1 \Big\{ r^{-d} \int_{E\cap B(x,r)}
r^{-1}\big| \langle y-x, n_{x,r} \rangle \big| dH^d(y) \Big\}^2
\, {dr \over r}
\leqno (13.45)
$$
where we set $\dsp n_{x,r} = H^d(E\cap B(x,r))^{-1} \int_{E\cap B(x,r)}
n(y) dH^d(y)$.

\ms\proclaim Corollary 13.46. 
Suppose that $n=d+1$, $E = \d \Omega$ for some $C^1$ 
domain $\Omega$, and that $\varepsilon$, $E$, $U$,
and $\Sigma_0$ satisfy the hypotheses of Theorem 12.3.
Suppose in addition that there exist $M \geq 1$ and $C_0 \geq 1$
such that 
$$
H(x) \leq M \hbox{ and }
H^d(E\cap B(x,r)) \leq C_0 r^d %v
\ \hbox{ for }
x \in E_0 = \big\{ x\in E \, ; \, \dist(x,U) \leq 1 \big\}.
\leqno (13.47)
$$
Then there exists $K = K(n,d,C_0,M)$ such
that, if we choose the planes $P_{j,k}$ in the %v
proof of Theorem~12.3 correctly, then $g$ is $K$-bi-Lipschitz,
i.e., (12.34) holds.

\ms
Of course $K$ does not depend on any $C^1$ estimates that we may 
have on $E$. We want to deduce this from Corollary 13.4, 
so it is enough to check that for $x\in E_0$ and $0 < r < 1$,
$$
\beta_1(x,r) \leq C r^{-d} \int_{E\cap B(x,r)}
r^{-1}\big| \langle y-x, n_{x,r} \rangle \big| dH^d(y).
\leqno (13.48)
$$
%vv  added the control on disrete sums, and added numbers.
Indeed, (13.48) implies that $\int_{0}^1 \beta_1(x,r)^2
\, {dr \over r} \leq CM$ for $x\in E_0$, by (13.45) and (13.47),
and hence
$$\eqalign{
J_1(x) &= \sum_{k \geq 3} \beta_1(x,r_k)^2 
\leq  C \sum_{k \geq 3} \int_{r_{k-1}}^{r_k} \beta_1(x,r_k)^2
\, {dr \over r} 
\cr&
\leq 10^{2d} C \sum_{k \geq 3} \int_{r_{k-1}}^{r_k} \beta_1(x,r)^2
\, {dr \over r}  
\leq 10^{2d} C \int_{0}^1 \beta_1(x,r)^2 \, {dr \over r}
\leq C'M
}\leqno (13.49)
$$
by the definitions (13.5) and (13.1), which is enough
to apply Corollary 13.4.
In fact, as soon as we check that
$$
|n_{x,r}| \geq C^{-1}
\ \hbox{ for $x\in E_0$ and } 0 < r \leq 1,
\leqno (13.50)
$$
(13.48) will follow from the definition (1.15), just
because if $P$ denotes the hyperplane through $x$ orthogonal
to $n_{x,r}$, 
$\big| \langle y-x, n_{x,r} \rangle \big| =|n_{x,r}|^{-1}
\dist(y,P)$ for $y\in E\cap B(x,r)$.

So we want to check (13.50), and we apply Green's formula
to the Caccioppoli set $\Omega \cap B(x,r)$. We get that
$$
\int_{\d} \wt n(y) dH^d(y) = 0,
\leqno (13.51) %v put a tilde to distinguish
$$
where $\d$ is the reduced boundary of $\Omega \cap B(x,r)$,
and $\wt n(y)$ denotes the approximate outward unit normal to
$\Omega \cap B(x,r)$ at $y$. 
(See for instance [Gi] or [AFP].) %% les livres
On $\d_1 = \d \cap B(x,r)$,
$\d$ is the same as $E$, and $\wt n(y)$ is equal to the unit
normal to $\d \Omega$ at $y$, so
$$
\int_{\d_1} \wt n(y) dH^d(y) 
=\int_{E\cap B(x,r)} n(y) dH^d(y) 
= H^d(E\cap B(x,r)) \, n_{x,r}.
\leqno (13.52)
$$

Let $P = P(x,r)$ be the hyperplane promised by 
(12.5), and set $\d_2 = \big\{ y \in \d \cap \partial B(x,r)
\, ; \, \dist(y,P) \leq \varepsilon r \big\}$. 
See Figure 2. 
Then
$$
\big| \int_{\d_2} \wt n(y) dH^d(y) \big| 
\leq H^d(\d_2) \leq H^d(\{ y\in \partial B(x,r)
\, ; \, \dist(y,P) \leq \varepsilon r \})
\leq C \varepsilon r^d.
\leqno (13.53)
$$
Finally set $\d_3 = \d \sm (\d_1 \cup \d_2)
= \Omega \cap \big\{ y \in \partial B(x,r) %v
\, ; \, \dist(y,P) > \varepsilon r \big\}$.
Denote by $\Sigma_{+}$ and $\Sigma_-$ the two components
of $\big\{ y \in \partial B(x,r)
\, ; \, \dist(y,P) > \varepsilon r \big\}$;
we claim that $\d_3 = \Sigma_+$ or $\d_3 = \Sigma_-$. 
Indeed, since $x\in \d \Omega$, we can find points
$y_1 \in \Omega$ and $y_2 \in \R^n \sm \overline\Omega$ 
very close to $x$; then we can use Theorem 12.3 to find
a path $\gamma_1$ in $\R^n \sm E$ that goes from 
$y_1$ to $\Sigma_{+} \cup \Sigma_-$, and similarly for 
$y_2$. To be fair, Theorem 12.3 is not really needed here,
and we could construct $\gamma_1$ by hand, by concatenating
successive intervals going away from $E$, that we would draw at 
different scales.
Since $y_1$ and $y_2$ lie in different components
of $\R^n \sm E$, they are connected to different $\Sigma_\pm$,
and the claim follows.

\vglue 3.7cm %%% verifier le display page 77
\hskip 3.2cm  
\epsffile{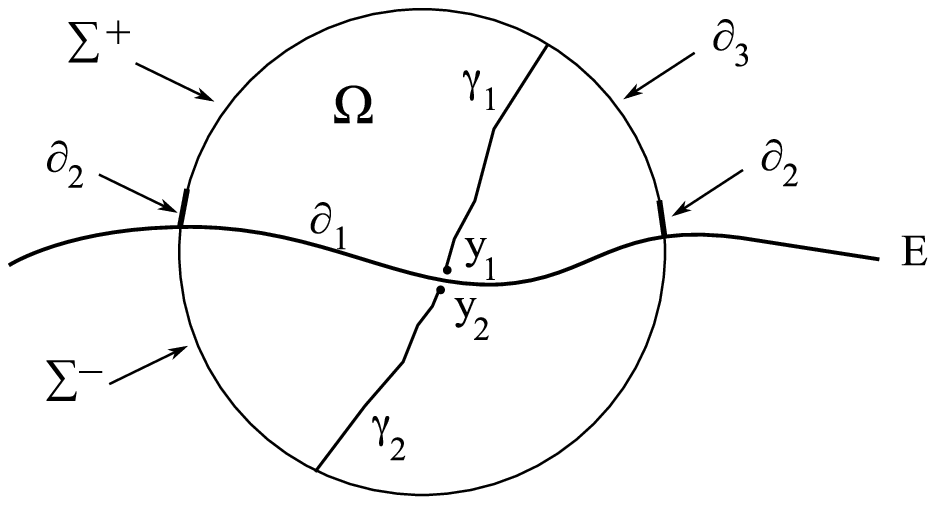} 
\medskip
\centerline{{\bf Figure 2.}}
\medskip

Suppose for definiteness that $\d_3 = \Sigma_+$,
and let $n_+$ denote the unit normal to $P$ that goes
in the direction of $\Sigma_+$. Then
$$
\int_{\d_3} \wt n(y) dH^d(y) = \omega_d \,
(1-\varepsilon^2)^{d/2} \, r^d \, n_+ \, ,
\leqno (13.54)
$$
where $\omega_d$ is the $H^d$-measure of the unit disk in $R^d$
(apply Green's formula to the convex hull of $\Sigma_+$).
Altogether,
$$\eqalign{
H^d(E\cap B(x,r)) \, |\wt n_{x,r}| &= 
\Big|\int_{\d_1} \wt n(y) dH^d(y)\Big| 
=\Big|\int_{\d_2 \cup \d_3} \wt n(y) dH^d(y)\Big| 
\cr&\geq\Big|\int_{\d_3} \wt n(y) dH^d(y)\Big|
- \Big|\int_{\d_2} \wt n(y) dH^d(y)\Big|
\geq c (1-C\varepsilon) r^d
}\leqno (13.55)
$$
by (13.51)-(13.54), which implies (13.50) 
because of (13.47). Corollary 13.46 follows.
\qed
%%% Second nice place

\msi{\bf Remark 13.56.}
We defined $H$ and stated Corollary 13.46 mostly because of
S. Semmes' characterization [Se] %% CASSC
of the Chord-Arc Surfaces with Small Constants (CASSC)
by the fact that the unit normal has a small BMO norm.
Here we allow $H$ to be large (even as a $BMO$ function), but
bounded. This boundedness assumption cannot be hoped to be 
necessary. Also see Remark 15.43 for a rapid discussion
on the available results for the various forms of uniform 
rectifiability.

It could be interesting to know whether 
$H(x)$ can be replaced with the simpler-looking
$$
\widetilde H(x) = \int_{0}^1 r^{-d} \int_{E\cap B(x,r)}
|n(y)-n_{x,r}|^2 \, {dH^d(y) dr \over r}
\leqno (13.57)
$$
in the statement of Corollary 13.46.

\msi{\bf Remark 13.58.}
In higher codimension, we could try to set
$$
H(x) = \int_{0}^1 \Big\{ r^{-d} \int_{E\cap B(x,r)}
r^{-1}\big| \pi^\perp_{x,r}(y-x) \big| dH^d(y) \Big\}^2
\, {dr \over r}
\leqno (13.59)
$$
where $\pi^\perp_{x,r}$ would be some average on 
$E \cap B(x,r)$ of $\pi^\perp(y)$. It is tempting to set
$\dsp \pi^\perp_{x,r} = H^d(E\cap B(x,r))^{-1} \int_{E\cap B(x,r)}
\pi^\perp(y) dH^d(y)$. In this case we did not 
check whether there is a $d$-plane $P(x,r)$ such that
$\dist(y,P(x,r)) \leq C \big| \pi^\perp_{x,r}(y-x) \big|$
for $y\in B(x,r)$, as we deduced from (13.50) when $d=n-1$.

\bsi
{\bf 14. Big pieces of bi-Lipschitz images and approximation
by bi-Lipschitz domains}
\ms

When $E$ is both locally Reifenberg-flat (as in Theorem 12.3),
and locally Ahlfors-regular (as in (13.2)), but the function 
$J_1$ of (13.5) is not bounded, we can still find very big
pieces of bi-Lipschitz images of $\Sigma_0$ inside $E$.

We shall also see that we can find nicer domains
$\Omega \i \R^n \sm E$, that are bi-Lipschitz images of
saw-tooth domains, and whose boundaries contain very large parts
of $E$. 

Both results will rely on the proof of Corollary 13.4, 
but also on the local uniform rectifiability of $E$,
which will be proved in the next section as a consequence
of Theorem~12.3 and the local Ahlfors-regularity of $E$.

Let us first state the result about very big bi-Lipschitz
pieces, and worry about saw-tooth domains later
(in Proposition 14.16).

\ms\proclaim Theorem 14.1.
Let $\varepsilon > 0$ be small enough, depending 
on $n$ and $d$. For each choice of $C_0 \geq 1$
and $\tau > 0$, we can find 
$K = K(n,d,C_0,\tau)$ such that
the following holds. Let $E, U, \Sigma_0 \in \R^n$ 
be given, and assume that (12.1) and (12.4) hold. 
Also assume that $E$ is closed
and that for $x\in E \cap U^+$ and $r\in (0,1]$, 
$$
H^d(E\cap B(x,r)) \leq C_0 r^d
\leqno (14.2)
$$
and there is an affine $d$-plane $P = P(x,r)$ through $x$ such that
(12.5) holds. Then there is a bi-Lipschitz mapping
$g : \R^n \to \R^n$, such that
$$
H^d(E \cap B(x,1/2) \sm g(\Sigma_0)) \leq \tau
\ \hbox{ for } x\in E \cap U,
\leqno (14.3)
$$
$$
K^{-1} |x-y| \leq |g(x)-g(y)| \leq K |x-y|
\hbox{ for } x,y \in \R^n,
\leqno (14.4)
$$
$$
|g(z)-z| \leq C \varepsilon
\ \hbox{ for } z\in \R^n,
\leqno (14.5)
$$
and 
$$
g(z)=z 
\ \hbox{ for $z\in \R^n$ such that }
\dist(z,U) \geq 13.
\leqno (14.6)
$$

\ms
We wrote (14.3) in a slightly strange localized way 
because we also want to allow the case when $E$ and $U$ 
are unbounded. Also, we did not include a lower bound
in (14.2), because it is provided by Lemma 13.6.

The proof will use the fact that there is a constant
$C_1 \geq 0$, that depends on $n$, $d$, and $C_0$, such 
that
$$
\int_{E \cap B(x,1/2)} J_1(x) \, dH^1(x) \leq C_1
\ \hbox{ for } x\in E \cap U, 
\leqno (14.7)
$$
and where $J_1$ is still as 
in (1.15) and (1.16) or (13.1) and (13.5). 
We could also have worked with some
other exponent $q < {2d \over d-2}$, but there is not
much point because $J_1$ is smaller (by (13.3)). 
We shall prove (14.7) in the next section, 
as a consequence of the local uniform rectifiability
of $E$ near $U$. In the mean time, let us check that
Theorem~14.1 follows from (14.7). Set 
$$
E'_0 = \big\{ x\in E_0 \, ; \, J_1(x) \leq \tau^{-1} C_1 \big\}
= \big\{ x\in E \, ; \, \dist(x,U) \leq 1 \hbox{ and }
J_1(x) \leq \tau^{-1} C_1 \big\}.
\leqno (14.8)
$$
Then
$$
H^d(E \cap B(x,1/2) \sm E'_0) \leq \tau
\ \hbox{ for } x\in E \cap U,
\leqno (14.9)
$$
by (14.7) and Chebyshev.

Now let us construct $g$ as we did for Corollary 13.4, 
except that for each $k \geq 0$, our original
collection of points $\widetilde x_{j,k}$, $j\in J_k$,
is now chosen to be a maximal collection in  $E'_0$ 
(instead of $E_0$) subject to the same constraint that
$|\widetilde x_{i,k}-\widetilde x_{j,k}|\geq 4r_k/3$ 
for $i \neq j$. Compare with the description above (13.15).
Other than that, we define the $x_{j,k}$ and the $P_{j,k}$
as before.

Our various estimates on the relative distances between
the planes $P_{j,k}$, all the way up to (13.34) included,
still hold now: they would even hold on the larger $E_0$.
In particular, we can use Theorem 2.15 to construct $g$. 

Let $E_\infty$ denote the limit set of (2.19). Let
us check that
$$
\overline {E'_0} = E_\infty \i g(\Sigma_0).
\leqno (14.10)
$$
Indeed the proof of (13.16) yields
$$
E'_0 \i \cup_{j\in J_k} B(x_{j,k},5r_k/3)
\leqno (14.11)
$$
for each $k \geq 0$, and so $\overline {E'_0} \i E_\infty$.
Conversely, $E_\infty \i \overline {E'_0}$ because 
$\dist(x_{j,k},E'_0) \leq |x_{j,k}-\widetilde x_{j,k}| \leq r_k/3$
for $k \geq 0$ and $j \in J_k$ (by (13.15)). 
The fact that $E_\infty \i g(\Sigma_0)$ is just 
the last statement in Theorem~2.15, so (14.10) holds.
Now (14.9) and (14.10) say that
$$
H^d(E \cap B(x,1/2) \sm g(\Sigma_0))
\leq H^d(E \cap B(x,1/2) \sm E_\infty)
\leq H^d(E \cap B(x,1/2) \sm E'_0) \leq \tau
\leqno (14.12)
$$
for $x\in E \cap U$, so (14.3) holds.
Since (14.5) and (14.6) follow from the proof of
(12.7) and (12.6), we just need to check the bi-Lipschitz
estimate (14.4).

We still want to proceed as in the proof of Corollary 13.4,
so we choose $\overline z \in E'_0$ such that
$|\overline z-f(z)| \leq 2 \dist(f(z)),E'_0)$ (as in (13.35)).
Then we follow quietly the proof of Corollary 13.4.
It is still enough to prove (13.36), because (13.37)
holds with $M = \tau^{-1} C_1$ (by (14.8)). The only
place where our current choice of $\overline z$
could make a difference is (13.42) and (13.43),
but they still hold with the same proof, because 
now $\widetilde x_{j,k} \in E'_0$. Note that $E'_0$ plays the role 
$E_0$ did in the proof of Corollary 13.4.
So (2.24) holds for the same reasons as before,
Theorem 2.23 applies, and $g$ is bi-Lipschitz.
This completes our proof of Theorem~14.1 modulo (14.7).
\qed

\msi{\bf Remark 14.13.}
Here we decided to stop when the function $J_1$ becomes too
large. We could also decide to replace $E_0$ with 
$E'_0 = \big\{ x\in E_0 \, ; \, J_\infty(x) \leq M \big\}$
for some large $M$ (compare with (14.8)), and then modify 
the proof of Corollary 12.44 as we just did for Corollary 13.4.

This would give a bi-Lipschitz mapping $g$, with the usual
properties, including (14.10). But of course this is
more interesting if we know that some part of $E$ lies 
in $E'_0$, i.e., if we have some control on the restriction
of $J_\infty$ to $E_0$. In the case of Theorem 14.1, we got 
the corresponding control from (14.7).

We could also do a similar stopping time in the context
of Corollary 12.33, and stop when $J$ gets too large.
But again this would work best if we controlled $J$.

\ms
Let us now state a result on approximating saw-tooth domains.
Keep the notations and assumptions of Theorem 14.1, set
$$
F_\infty = g^{-1}(E_\infty) \i \Sigma_0
\leqno (14.14)
$$
(by (14.10) and because $g$ is injective)
and define the saw-tooth domain 
$$
\Omega_A = [\R^n \sm V] \cup \big\{ z\in V \, ; \,
|q(z)| > A \dist(p(z),F_\infty) \big\},
\leqno (14.15)
$$
where $A \geq 1$ will be chosen soon, %v
$V = \big\{ z\in \R^n \, ; \, \dist(z,\Sigma_0) < 40
\big\}$, $p(z)\in \Sigma_0$, and $q(z)=z-p(z)$, 
are as in Lemma 10.1. See Figure 3. 

\vglue 4.9cm %%% verifier le display page 81
\hskip 1.4cm  
\epsffile{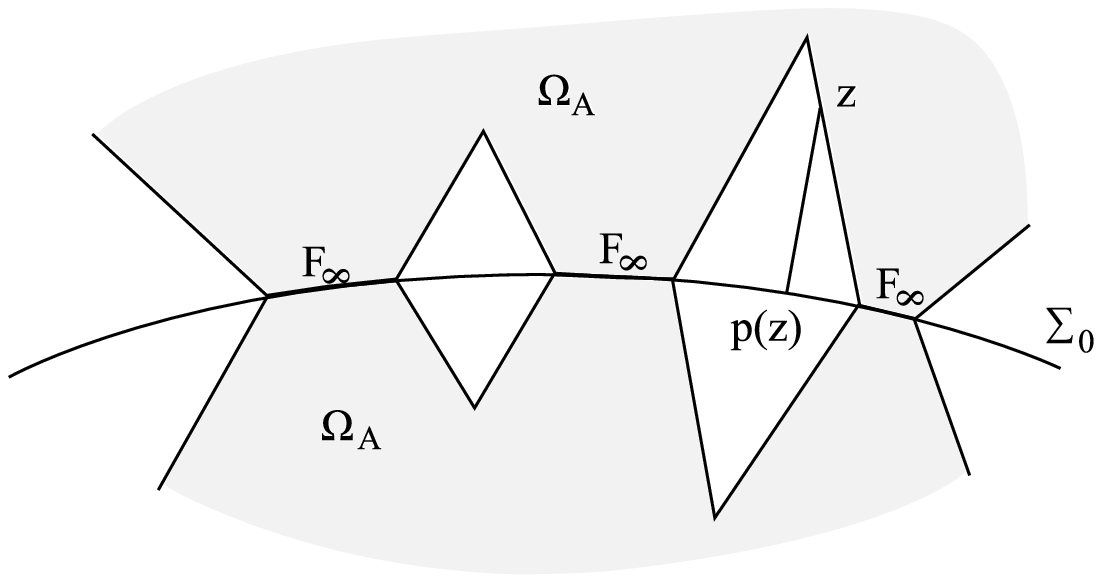} 
\medskip
\centerline{{\bf Figure 3.} The saw-tooth domain $\Omega_A$}
\medskip

Thus, if $\Sigma_0$ is a $d$-plane, $p(z)$ is the orthogonal
projection of $z$ on $P$, and $q(z)$ is the orthogonal projection
of $z$ on $P^\perp$.
If in addition $d=n-1$, $\Omega_A$ is composed of two
regions, whose common boundary is $F_\infty$, and that are
both bi-Lipschitz equivalent to a half space. When $d<n-1$,
think about a single domain which is invariant under rotations
along $\Sigma_0$ and such that 
$\d \Omega_A \cap \Sigma_0 = F_\infty$. %v
Note also that since $g$ is a bi-Lipschitz homeomorphism,
$g(\Omega_A)$ is a reasonably nice domain, whose
boundary is $g(\d \Omega_A)$. 

%%% Better location for Figure 3, if it worked

\ms\proclaim 
Proposition 14.16.
Keep the notation and assumptions of Theorem 14.1.
If $A$ is large enough, depending on $n$, $d$,
$C_0$, and $\tau$, then
$$
g(\Omega_A) \cap E \cap U= \emptyset
\ \hbox{ and }\ 
g(\d\Omega_A)\cap E\cap U = E_\infty\cap U. 
\leqno (14.17)
$$
When $n=d+1$ and $\Sigma_0$ is a plane, we even have
that $g(\d\Omega_A^{+})\cap E \cap U= 
g(\d\Omega_A^{-})\cap E \cap U = E_\infty\cap U$, 
where $\Omega_A^{+}$ and $\Omega_A^{-}$ denote the 
two connected components of $\Omega_A$.

\ms
Recall from (14.12) that 
$H^d(E \cap B(x,1/2) \sm E_\infty) \leq \tau$ for
$x\in E \cap U$, so (14.17) says that most of $E \cap U$ 
lies in $\d\Omega_A$. Thus, when $d=n-1$ and $\Sigma_0$
is a plane, for instance, Theorem 14.1 and Proposition 14.16
say that locally, $E$ contains very big pieces of bi-Lipschitz 
images of $\R^{d}$, but also give nice domains that are contained
in components of $\R^n \sm E$ and contain these big pieces %v
in their boundary. See Theorem 15.39, Remark~15.44, and 
Proposition~15.45 concerning a uniform version of this
and its relations with previous works. %v

We now prove the proposition. Let us first check that
$$
g(z) \in \R^n \sm E
\ \hbox{ when } z\in \overline\Omega_A \sm \Sigma_0
\hbox{ is such that } g(z) \in U.
\leqno (14.18)
$$
If $\dist(z,E) \geq 10^{-1}$, then 
$\dist(g(z),E) \geq \dist(z,E) -C\varepsilon >0$
by (14.5) and we are done. So assume that %v
$\dist(z,E) \leq 10^{-1}$. Choose $w\in E$ such that 
$|w-z| \leq \dist(z,E) + \varepsilon$. Then
$w\in U^+$ (because $g(z) \in U$ and
$|g(z)-z|\leq C\varepsilon$ by (14.5)), so %v
$$
\dist(z,\Sigma_0) \leq |w-z| + \dist(w,\Sigma_0)
\leq \dist(z,E) + 2\varepsilon \leq 10^{-1} + 2\varepsilon
\leqno (14.19)
$$ 
by (12.4). In particular, $z\in V$, $p(\cdot)$ 
and $q(\cdot)$ are defined near $z$, and we can find
$x\in F_\infty$ such that
$$
|p(z)-x| \leq 2 \dist(p(z),F_\infty) 
\leq {2 \over A}\, |q(z)|
\leqno (14.20)
$$ 
by (14.15) and hence %v
$$
|f(p(z))-f(x)| = |g(p(z))-g(x)| \leq {2K \over A}\, |q(z)|
\leq |q(z)|,
\leqno (14.21)
$$
because $f=g$ on $\Sigma_0$, by (14.4), %v
and if $A$ is large enough.

Let $k$ be such that $r_{k+1} \leq |q(z)| \leq r_k$;
note that $k \geq 0$ because $|q(z)| = \dist(z,\Sigma_0) <1$
by (10.12) and (14.19). Also,
$f(x) \in E_\infty = \overline{E'_0}$ 
(recall that $g=f$ on $\Sigma_0$ and use (14.14) and (14.10)),
so by (14.11) we can find $j\in J_k$ such that 
$f(x)\in 2B_{j,k}$. Then 
$$\eqalign{
|g(z)-x_{j,k}| 
&\leq |g(z)-f(p(z))|+|f(p(z))-f(x)|+|f(x)-x_{j,k}| 
\cr&
\leq |g(z)-f(p(z))|+|q(z)|+2 r_k
\leq (2+C\varepsilon)|q(z)|+2 r_k < 5r_k
}\leqno (14.22)
$$
by (14.21), (10.25) and the convention (10.19),
and by definition of $k$.

Obviously (14.18) will follow as soon as we check that 
$$
\dist(g(z),E) \geq {1 \over 2} \, |q(z)|,
\leqno (14.23)
$$
because $|q(z)| = \dist(z,\Sigma_0)>0$ by (10.12) and
because $z \notin \Sigma_0$.
Suppose that (14.23) fails, and let $\xi \in E$ be 
such that $|g(z)-\xi| \leq 2|q(z)|/3$. Notice that
$|\xi-x_{j,k}| \leq |g(z)-x_{j,k}|+2|q(z)|/3 < 6r_k$
by (14.22) and the definition of $k$, so
$$
\xi \in E\cap 6B_{j,k} \i E\cap  B(\widetilde x_{j,k},7r_k),
\leqno (14.24)
$$
because $|\widetilde x_{j,k} - x_{j,k}|< r_k$ by (13.15).
Then
$$
\dist(\xi,P(\widetilde x_{j,k}, r_k)) \leq \varepsilon r_k
\leqno (14.25)
$$
by (12.5). If (13.17) fails, we chose 
$P_{j,k}=P(\widetilde x_{j,k}, r_k)$, so 
$$
\dist(\xi,P_{j,k}) \leq C\varepsilon r_k
\leqno (14.26)
$$
as well. Otherwise, (13.25) says that
$$
d_{\widetilde x_{j,k},10r_k}(P_{j,k},P(\widetilde x_{j,k},r_k))
\leq C \varepsilon,
\leqno (14.27)
$$
and again (14.26) holds because of (14.25).

Recall from Proposition 5.4 that in $49B_{j,k}$, 
$\Sigma_k$ coincides with a $C\varepsilon$-Lipschitz
graph over $P_{j,k}$, that passes within
$C \varepsilon r_k$ of $x_{j,k}$. Thus (by (14.24)
and (14.26)) we can find $w\in \Sigma_k$ such that 
$|w-\xi| \leq C \varepsilon r_k$. Finally,
$\dist(w,\Sigma) \leq C \varepsilon r_k$ by (6.9), so
$\dist(\xi,\Sigma) \leq C \varepsilon r_k$,
and
$$
\dist(g(z),\Sigma) \leq |g(z)-\xi| + \dist(\xi,\Sigma)
\leq {2|q(z)| \over 3} \, + \dist(\xi,\Sigma)
\leq {2|q(z)| \over 3} \, + C \varepsilon r_k.
\leqno (14.28)
$$
Recall that $k$ was chosen so that $r_{k+1} \leq |q(z)| \leq r_k$,
so (14.28) implies that 
$$
\dist(g(z),\Sigma) \leq {3|q(z)|\over 4}
= {3 \dist(z,\Sigma_0)\over 4}
\leqno (14.29)
$$
by (10.12), and in contradiction with (10.27).
This proves (14.23), and (14.18) follows.

Let us now check (14.17). Suppose we can find
$w \in g(\Omega_A) \cap E \cap U$, and let $z \in \R^n$
be such that $g(z)=w$. Then $z\in \Omega_A \sm \Sigma_0$
(because $\Omega_A$ does not meet $\Sigma_0$),
so (14.18) says that $g(z) \notin E$, a contradiction.

Next suppose that 
$w \in g(\d\Omega_A) \cap E \cap U$, and again
write $w = g(z)$. If $z\in \overline\Omega_A \sm \Sigma_0$,
(14.18) gives a contradiction as before. So
$z\in \Sigma_0$. Since $z\in \d\Omega_A$, there is a
sequence $\{ z_k \}$ in $\Omega_A$ that converges
to $z$. By (14.15), $z_k \in V$ for $k$ large,
$q(z_k) = \dist(z_k,\Sigma_0)$ tends to $0$
(by (10.12)), and hence $\dist(p(z_k),F_\infty)$
tends to $0$ too. But $z_k = p(z_k) + q(z_k)$,
so $\dist(z_k,F_\infty)$ tends to $0$, and $z\in F_\infty$
because $F_\infty = g^{-1}(E_\infty)$ is closed (by
(14.10) or the definition of the limit set).
Then $w = g(z) \in E_\infty$ and we proved that
$g(\d\Omega_A)\cap E\cap U \i E_\infty\cap U$.

Conversely, if $w \in E_\infty\cap U$,
$z = g^{-1}(w)$ lies in $F_\infty \i \Sigma_0$
(by (14.14)), and (14.15) says that 
$z \in \partial \Omega_A$ (draw a line segment $L$
starting from $z$ and perpendicular to $\Sigma_0$ 
at $z$; the points of $L \sm \{ z \}$ that lie 
close enough  to $z$ all lie in $\Omega_A$). So (14.17) holds.

Finally assume that $n=d+1$ and $\Sigma_0$ is a plane,
and choose a connected component $\Omega_A^{\pm}$ of $\Omega_A$.
We already know that $g(\d\Omega_A^{\pm})\cap E \cap U \i E_\infty$,
by (14.17) and because $\d\Omega_A^{\pm} \i \d\Omega_A$.
But conversely, if $w \in E_\infty\cap U$, we know that
$z = g^{-1}(w)$ lies in $F_\infty$, and we get that
$z\in \partial \Omega_A^\pm$ by the same argument as before.
Proposition 14.16 follows.
\qed

\bsi
{\bf 15. Uniform rectifiability and Ahlfors-regular 
Reifenberg-flat sets}
\ms

The main goal of this section is to study the local
uniform rectifiability properties of locally 
Ahlfors-regular Reifenberg-flat sets of dimension $d$ in
$\R^n$. This will allow us to complete the proof of 
(14.7) and  %v
Theorem 14.1 after a small additional localization,
but for the moment we shall find it more convenient
to work with a simpler class of sets, because this will
make our statements simpler and more scale-invariant.

\ms\proclaim Definition 15.1.
Let $E \i \R^n$ be a nonempty closed set.
We say that $E$ is locally Ahlfors-regular and 
Reifenberg-flat of dimension $d$ (in short, $E$ is a 
LARRF set) if there exist constants $\varepsilon >0$ 
(always assumed to be small enough, depending on $n$ 
and $d$) and $C_0 \geq 1$ such that for each $x \in E$ 
and $0 < r \leq 1$, there is a $d$-plane
$P = P(x,r)$ through $x$ such that
$$
d_{x,r}(P,E) \leq \varepsilon,
\leqno (15.2)
$$
and also
$$
C_0^{-1} r^d \leq H^d(E\cap B(x,r)) \leq C_0 r^d
\ \hbox{ for $x\in E$ and } 0 < r \leq 1.
\leqno (15.3)
$$

\ms
We could define ``LARRF up to scale $r_0$", where we 
ask (15.2) and (15.3) to hold for $x\in E$ and
$0 < r \leq r_0$ (but we can easily reduce to $r_0=1$),
or ``globally LARRF" (when (15.2) and (15.3) hold for 
all $r$), or even give a more standard and less uniform
definition of ``local", where for each compact subset 
$K$ of $E$, there is an $r_K > 0$ such that (15.2) and (15.3) 
hold for $x\in K$ and $0 < r \leq r_K$. The results below 
would still hold with the expected modifications, because 
they are local in essence.

We decided to include the lower bound in (15.3) in the 
definition (because this is always included in the
definition of local Ahlfors-regularity),
but we know from Lemma~13.6 that it follows from the local
Reifenberg-flatness (15.2), and with a constant that depends 
only on $n$ and $d$. %v

Of course (15.2) is a rather strong condition,
so it should not be a surprise that LARRF sets are locally
uniformly rectifiable, with big pieces of Lipschitz graphs,
as in the following statement.

\ms\proclaim Theorem 15.4.
Let $E$ be a LARRF set, with $\varepsilon$ small enough
(depending only on $n$ and $d$). 
Then there exist constants $\theta > 0$ and $M > 0$, that 
depend only on $n$, $d$, and $C_0$, %v
such that for $x \in E$ and $0 < r \leq 1$, 
we can find a $d$-plane $P$ and an $M$-Lipschitz function
$F : P \to P^\perp$, such that
$$
H^d(E\cap \Gamma_F \cap B(x,r)) \geq \theta r^d,
\leqno (15.5)
$$
where $\Gamma_F$ denotes the graph of $F$ over $P$.

\ms
Recall that $P^\perp$ is the $(n-d)$-dimensional
vector space orthogonal to $P$, that $F$ is $M$-Lipschitz
when $|F(z)-F(w)| \leq M |z-w|$ for $z,w\in P$,
and that $\Gamma_F = \big\{ z+F(z) \, ; \, z\in P \big\}$.

We refer to [DS3] and its references %% UR
for general information about uniform rectifiability,
such as the fact that if $E$ contains big pieces
of Lipschitz graphs as in Theorem 15.4, it has all sort
of other interesting properties. We shall use one of them
soon, the so-called geometric lemma.

For the proof of Theorem 15.4, it is tempting to
try to prove that $E$ satisfies the bilateral weak geometric
lemma (BWGL) of Definition I.2.2 in [DS3]. %% UR
Recall that $E$ satisfies the BWGL locally if for 
each $\varepsilon >0$, (15.2) holds for all $x\in E$ and 
$0 \leq r \leq 1$, except perhaps for a Carleson set of
exceptional pairs $(x,r)$. In fact, in view of the main
theorems (such as Proposition II.2.2 in [DS3]), %% UR
we only need to check this for one small $\varepsilon >0$.
See Remark~II.2.5 in [DS3]. %% UR
However, unfortunately this small $\varepsilon$
seems to depend on the Ahlfors regularity constant $C_0$.
Thus, if we want to apply Proposition II.2.2 in [DS3] %% UR
directly to prove Theorem 15.4, we seem to be required 
to take $\varepsilon$ small enough, depending on $C_0$,
which we will not do.
Other criteria from [DS3] %% UR
seem to suffer from the same apparent defect.

To be fair, it is probable that in this simpler case, 
the proof of Proposition II.2.2 in [DS3] %% UR
goes through when $\varepsilon$ is small enough,
depending only on $n$ and $d$, but we prefer not
to check this, even though this would be a reasonable
option for the reader.

Anyway we shall deduce Theorem 15.4 from 
Proposition 3 in [Da1], %% Morceaux
so we want to check that $E$ is a local generalized
Semmes surface, as follows.

\ms\proclaim Lemma 15.6. 
Let $E$ be a LARRF set, with $\varepsilon$ small enough
(depending only on $n$ and $d$).
Then for $x\in E$ and $0 < r \leq 10^{-2}$, there is an 
$(n-d)$-plane $W$ through $x$ such that if we set
$S = W \cap \partial B(x,r/2)$, then
$$
\dist(z,E) \geq r/3 \ \hbox{ for } z \in S
\leqno (15.7)
$$
and $S$ is linked with $E$ in the sense that
there is no continuous function 
$F : [0,1] \times \R^n \to \R^n$ such that
$$
F(t,z) = z 
\ \hbox{ when $t = 0$ or } z \in \R^n \sm B(x,10r),
\leqno (15.8)
$$
$$
F(1,z) \in \R^n \sm B(x,10r)
\ \hbox{ for } z\in E,
\leqno (15.9)
$$
and 
$$
F(t,z) \in \R^n \sm S
\ \hbox{ when $z\in E$ and } 0 \leq t \leq 1.
\leqno (15.10)
$$

\ms
In other words, we cannot move $E$ away 
without crossing $S$. %v
Compare with Definition~3 in [Da1]. %% Morceaux
Here $W$ will be the plane through $x$ which
is perpendicular to $P$, where $P=P(x,20r)$ comes
from (15.2). Then (15.7) holds trivially
if $\varepsilon \leq 10^{-1}$, say, but we need
some argument to show that $S$ is linked with $E$.
And for this it will be quite pleasant to use
Theorem 12.3.

By translation and dilation invariance, we may assume
that $x=0$, $r=10$, and (15.2) holds for radii smaller than
or equal to $10^3$. Then apply Theorem 12.3 with 
$\Sigma_0 = P = P(0,200)$ and $U = B(0,102)$; 
the assumptions (12.1), (12.4), 
and (12.5) are clearly satisfied (but with the
larger constant $10^3 \varepsilon$), so we get a
bih\"{o}lder mapping $g: \R^n \to \R^n$ such that
in particular
$$
g(z) = z 
\ \hbox{ for } z \in \R^n \sm B(0,115)
\leqno (15.11)
$$
$$
|g(z)-z| \leq C \varepsilon 
\ \hbox{ for } z \in \R^n,
\leqno (15.12)
$$
and 
$$
E \cap B(x,102) = g(P) \cap B(x,102).
\leqno (15.13)
$$

Now we shall assume that we can find a homotopy $F$ as in 
(15.8)-(15.10), and derive a contradiction. Define a new 
homotopy by setting
$$
G(t,z) = F(t,g(z))
\ \hbox{ for $0 \leq t \leq 1$ and } z \in \R^n.
\leqno (15.14)
$$
Observe that
$$
|G(t,z)-z| \leq C \varepsilon
\ \hbox{ for } z \in \R^n \sm B(0,101),
\leqno (15.15)
$$
because $g(z) \in \R^n \sm B(0,100)$ by (15.12),
and then $G(t,z) = F(t,g(z)) = g(z)$ by (15.8).
Note also that
$$
|G(0,z)-z| = |g(z)-z| \leq C \varepsilon
\ \hbox{ for all } z\in \R^n, 
\leqno (15.16)
$$
by (15.8) and (15.12), and 
$$
g(z) \in E %\cap B(0,102)
\ \hbox{ for } z\in P\cap B(0,101)
\leqno (15.17)
$$
because $g(z) \in B(0,102)$ by (15.12),
and then $g(z) \in E$ by (15.13) and because $z\in P$.
Then
$$
G(1,z) \in \R^n \sm B(0,100)
\ \hbox{ for } z\in P
\leqno (15.18)
$$
because either $z\in \R^n \sm B(0,101)$
and this follows from (15.15), or else
$g(z) \in E$ by (15.17) and 
$G(1,z) = F(1,g(z)) \in \R^n \sm B(x,100)$
by (15.9). Similarly,
$$
G(t,z) \in \R^n \sm S
\ \hbox{ for } z\in P
\leqno (15.19)
$$
trivially by (15.12) when $z\in \R^n \sm B(0,101)$,
or else because $G(t,z) = F(t,g(z)) \in \R^n \sm S$
by (15.17) and (15.10).

Now the combination of (15.16), (15.18), and (15.19)
on $P$ is shocking, because it means that we can
make $P$ move away from $B(0,10)$ without crossing $S$. 
We want to find a contradiction by constructing a homotopy 
from the identity to a constant or a mapping of even degree 
on the unit sphere, but let us first modify $G$ slightly 
to make it cleaner at both ends. Since $G(1,z)$ may not 
be smooth, we choose a smooth function $G$ on $P$, such that
$|G(1,z)-G(z)| \leq 1$ for $z\in P$, and also
$$
G(z) = z 
\ \hbox{ for } z\in P \sm B(0,116).
\leqno (15.20)
$$
This last is easy to arrange, by (15.11) and (15.8). Also note that
$$
G(z) \in \R^n \sm B(0,99)
\ \hbox{ for } z\in P,
\leqno (15.21)
$$
by (15.18). We define a homotopy $\{ G_t \}, 0 \leq t \leq 1$, 
that goes from $G_0(z) = z$ to $G_1 = G$ by
$$
G_t(z) = (1-3t) z + 3t G(0,z)
\ \hbox{ for } z\in P \hbox{ and } 0 \leq t \leq 1/3,
\leqno (15.22)
$$
$$
G_t(z) = G(3t-1,z) 
\ \hbox{ for } z\in P \hbox{ and } 1/3 \leq t \leq 2/3,
\leqno (15.23)
$$
and
$$
G_t(z) = (3-3t) G(1,z) + (3t-2) G(z)
\ \hbox{ for } z\in P \hbox{ and } 2/3 \leq t \leq 1.
\leqno (15.24)
$$
Let us check that
$$
G_t(z) \in \R^n \sm S
\ \hbox{ for } z\in P \hbox{ and } 0 \leq t \leq 1.
\leqno (15.25)
$$
When $t \leq 1/3$, simply notice that $|G_t(z)-z|
\leq |G(0,z)-z| \leq C \varepsilon$ by (15.22) and
(15.16), so $G_t(z)$ is far from $S$ because 
$\dist(z,S) \geq 5$ for $z\in P$.
When $t \geq 2/3$, $|G_t(z)-G(1,z)| \leq |G(z) - G(1,z)| \leq 1$, 
so $G(z) \in \R^n \sm B(0,98)$ by (15.18)
and $G(z)$ is far from $S$. Finally, when $1/3 \leq t \leq 2/3$,
$G_t(z) = G(3t-1,z) \in \R^n \sm S$ by (15.23) and (15.19);
so (15.25) holds. Similarly,
$$
G_t(z) = z 
\ \hbox{ for $z\in P \sm B(0,116)$ and } 0 \leq t \leq 1,
\leqno (15.26)
$$
by (15.11), (15.8), and (15.20).

We want to define mappings from $\d B$ to itself,
where $\d B = \d B(0,1)$ is the unit ball of $\R^n$,
but first let us define a mapping 
$\varphi : \d B \sm P \to P \times S$. Let $v \in \d B\sm P$ 
be given. Write $v = v_1 + v_2$, with $v_1 \in P$ and 
$v_2 \in W = P^\perp$. Notice that $v_2 \neq 0$ because 
$v \notin P$. Then set
$$
\varphi_1(v) = 5v_1/|v_2| \in P
\ \hbox{ and } \ 
\varphi_2(v) = - 5 v_2 /|v_2| \in S
\leqno (15.27)
$$
(recall that $S = W \cap \d B(0,5)$ because $r=10$),
and $\varphi(v) = (\varphi_1(v),\varphi_2(v))$.
Notice that
$\varphi_1(v)-\varphi_2(v) = 5 v/|v_2|$, so
$$
v = {\varphi_1(v)-\varphi_2(v) \over
|\varphi_1(v)-\varphi_2(v)|}
\ \hbox{ for } v \in \d B \sm P.
\leqno (15.28)
$$
We define mappings $H_t : \d B \sm P \to \d B$ by
$$
H_t(v) = {G_t(\varphi_1(v))-\varphi_2(v) \over
|G_t(\varphi_1(v))-\varphi_2(v)|}
\ \hbox{ for $0 \leq t \leq 1$}
\leqno (15.29)
$$
(we move $\varphi_1(v)$ in (15.28), according to our
homotopy) and
$$
H_t(v) = {G(\varphi_1(v))-(2-t)\varphi_2(v) \over
|G(\varphi_1(v))-(2-t)\varphi_2(v)|}
\ \hbox{ for $1 \leq t \leq 2$}
\leqno (15.30)
$$
(now we contract $S$ to the origin). The denominator never vanishes:
in (15.29) because of (15.25), and in (15.30) because
of (15.21).

Next we claim that $(t,v) \to H_t(v)$ has a continuous
extension to $[0,2] \times \d B$. %v
When $v$ is close enough to $P$, for instance as soon as 
$|v_2| \leq 10^{-2}$,
$|\varphi_1(v)| = 5 |v_1|/|v_2| > 400$, and hence
$G_t(\varphi_1(v)) = \varphi_1(v)$ for all $t$, by (15.26),
and also $G(\varphi_1(v)) = \varphi_1(v)$. Thus
$$
H_t(v) = {\varphi_1(v)-\varphi_2(v) \over
|\varphi_1(v)-\varphi_2(v)|}
\ \hbox{ for $0 \leq t \leq 1$}
\leqno (15.31)
$$
and
$$
H_t(v) = {\varphi_1(v)-(2-t)\varphi_2(v) \over
|\varphi_1(v)-(2-t)\varphi_2(v)|}
\ \hbox{ for $1 \leq t \leq 2$.}
\leqno (15.32)
$$
Set $\rho(t) = 1$ for $0 \leq t \leq 1$ and
$\rho(t) = 2-t$ for $1 \leq t \leq 2$. Then 
(15.27) yields
$$
H_t(v) = { v_1 + \rho(t) v_2 \over |v_1 + \rho(t) v_2|}
\ \hbox{ for } 0 \leq t \leq 2.
\leqno (15.33)
$$
This map clearly has a continuous extension across
$\d B \cap P$ (where only $v_2$ tends to $0$).
Our continuity claim follows.

Notice that for $v\in \d B\sm P$, $H_0(v) = 
[\varphi_1(v)-\varphi_2(v)]/|\varphi_1(v)-\varphi_2(v)| =v$
by (15.29), because $G_0(z) = z$ for $z\in P$ (by (15.22)), 
and by (15.28). This is still true for $v\in \d B\cap P$,
by continuity, so $H_0(v)=v$ on $\d B$, and we shall
reach the desired contradiction as soon as we prove
that $H_2$ is not homotopic to the identity.

Apparently we need to distinguish between cases.
Let us first suppose that $d < n-1$, set $Z = H_2(\d B)$,
and check that $H^d(Z) < +\infty$. Notice that near
$P$, $H_2(v) = v_1 /|v_1|$ by (15.33), so the corresponding
part of $Z$ is contained in $P \cap \d B$. Far from
$P$, $H_2(v) = G(\varphi_1(v)/|G(\varphi_1(v)|$ by (15.30),
so the corresponding part of $Z$ is contained in the image
of a compact subset of $P$ (where $\varphi_1(v)$ lies)
by the smooth mapping $z \to G(z)/|G(z)|$ (see (15.21)
and recall that $G$ is smooth). Our claim follows.
Now $d < n-1$, so $Z$ is strictly contained in $\d B$,
which means that it omits some small ball $B'$.
But $\d B \sm B'$ can be contracted to a point (inside $\d B$), 
which implies that $G_2$ is homotopic to a constant (among
continuous functions from $\d B$ to $\d B$), a contradiction.

When $d=n-1$, $S$ is composed of two points, and it should
even be more obvious that we cannot deform $P$ across these 
points. The simplest at this point is to observe that by
the discussion above $H_2(v)$ depends only on $v_1$, i.e., 
is symmetric with respect to the hyperplane $P$; 
this forces its degree to be even (think about 
the number of inverse images at a regular point), and hence it 
is not homotopic to the identity (see [Du]). %% Dugundji
%%% Curious that this is more complicated! 
% find another mapping, or a reason why the image is small??

This contradiction with the existence of the homotopy
$F(t,z)$ above completes our proof of Lemma 15.6.
\qed

\msi{\bf Proof of Theorem 15.4.}
As we said before, Theorem 15.4 is now a consequence of
Proposition 3 in [Da1]; %%  Morceaux
the statement is not exactly the same, because
in that reference, the set $E$ is unbounded and
the assumptions and  the conclusions
both hold for all $x\in E$ and $r > 0$.
However, the proof is local (the proof of (15.5) in a 
given ball $B$ never uses information on $E \sm CB$),
and goes through in the present context.
So Theorem 15.4 follows from Lemma 15.6.
\qed

\msi{\bf Remark 15.34.}
When $n=d+1$, our proof of Theorem 15.4 could in principle be 
simplified slightly. Even though the condition of Lemma 15.6 
is supposed to be a generalization of the so-called Condition B, 
it is set up a little differently, which forced us to spend
a little more time with topology than we should have.

The global version of Condition B is that, for each choice
of $x\in E$ and $0 < r <\,{\rm diam}(E)$, we can find 
$y_1, y_2 \in B(x,r) \sm E$, that lie in different connected
components of $\R^{d+1} \sm E$, and also such that
$\dist(y_i,E) \geq C^{-1} r$ for $i=1,2$.
Thus the linking condition is a little simpler
than in Lemma 15.6, because we just need to check that
the two points of $S$ lie in different components of
$\R^{d+1} \sm E$. But the right way to localize
this is to require that for some $C \geq 10$
and all $x\in E$ and $0 < r \leq C^{-1}$, 
we can find $y_1, y_2 \in B(x,r) \sm E$, that lie in 
different connected components of $B(x,Cr) \sm E$, 
and such that $\dist(y_i,E) \geq C^{-1} r$ for $i=1,2$.
Unfortunately, most of the proofs of the fact that
every Ahlfors-regular set with Condition B is uniformly
rectifiable and contains big pieces of Lipschitz graphs
([Da1], [DJ], [DS2])   %% Morceaux ; Jerison ; Quantitativetransac
do not mention this way to localize, and the reader
would have to use Theorem~3.5 and the WNPC (weak no Poincar\'{e}
condition) in [DS4] to get a proof. %% UR-IHP
Even that way, the statement says that $C$ should
be large enough, depending on the Ahlfors-regularity
constant $C_0$, and the verification would use Theorem 12.3.

\ms
We now use Theorem 15.4 to prove that $E$ satisfies
a local form of the so-called geometric lemma.
We allow all $q < {2d \over d-2}$ in the statement,
but we shall only use $q=1$.

\ms\proclaim Corollary 15.35.
Let $E$ be a LARRF set, with $\varepsilon$ small enough,
depending only on $n$ and $d$. For each exponent $q$
such that $1 \leq q < {2d \over d-2}$ (we allow
$1 \leq q \leq +\infty$ when $d=1$), there is a constant
$C_q = C(n,d,C_0,q)$ such that
$$
\int_{y\in E \cap B(x,t)}\int_{0}^r 
\beta_q(x,r)^2 \,  {dt \over t}dH^d(y) 
\leq C_q r^d
\leqno (15.36)
$$
for $x\in E$ and $0 < r \leq 1$, where $\beta_q(x,r)$
is still defined as in (1.15).

\ms
This follows from the local version of the fact that
Condition (C6) on page 13 of [DS1] %% asterisque
implies Condition (C3) on pages 11-12 (in its 
$\beta_q$ version mentioned there). As usual, the 
local version is not mentioned in [DS1], %% asterisque
but the proof is the same. This is a magnified nonlinear
version for sets of a result of Dorronsoro [Do] % rr r OK
on the good approximation of Lipschitz functions by 
affine functions in most balls.

Note that in the derivation of Corollary 15.35
from Theorem 15.4, the reader should only expect to get 
(15.36) for $r \leq C^{-1}$; however this makes no difference
because\break
$\dsp\int_{C^{-1}r}^r \beta_q(x,t)^2 \, {dt \over t} \leq C$
anyway, just because $\beta_q(x,r) \leq C$ by (1.15) and (15.3).
\qed

\ms
We stated (15.36) with a continuous integral because
this is the way it shows up in [DS1], %% asterisque
but an easy consequence of (15.36) is that
$$
\int_{y\in E \cap B(x,r)} \sum_{ k \geq 0 \, ; \, r_k \leq r}
\beta_q(x,r_k)^2 \, dH^d(y) \leq C_q r^d
\leqno (15.37)
$$
for $x\in E$ and $0 < r \leq 1$. Indeed
$\beta_q(x,r_k) \leq 10^{1+d/q} \beta_q(x,r)$ for
$r_k \leq r \leq 10 r_k$, which gives a control
on the indices $k$ such that $r_k \leq r/10$. For
the last one, we just say that $\beta_q(x,r_k) \leq C$.

\msi{\bf Proof of (14.7) and Theorem~14.1.}  
Recall from the first part of (13.5)  
that $J_1(x) = \sum_{k \geq 0} \beta_1(x,r_k)^2$;
so (14.7) is the same as (15.37) with $q=1$ and $r=1/2$ 
(and we just add a bounded term coming from $k=0$).
To be fair, we need to localize one more time:
our set in Theorem~14.1 is not a LARRF set, but
only satisfies (15.2) for $x\in E \cap U^+$
and (15.3) for $x\in E$ such that $\dist(x,U) \leq 3/2$
(the smaller range comes from the lower bound,
which we get from Lemma 13.6).
Again, we claim that the proofs of the various theorems 
used above (namely, Proposition 3 in [Da1] %%  Morceaux
and the result of [DS1]) %% asterisque
go through in this context, with only minor
modifications (we never use large radii or
faraway points). Our proof of Theorem~14.1
is now complete.
\qed

\msi{\bf Remark 15.38.}
We announced earlier that the sufficient condition
of Corollary 13.4 has the right flavor. Indeed,
the Carleson condition (15.36) can be used to prove
that the function $J_1$ lies in $BMO_{loc}(E)$, and
by the result of [DS1] %% asterisque
it is satisfied for every locally Ahlfors-regular
uniformly rectifiable set $E$. Thus it is a necessary
condition for $E$ to be (contained in) a bi-Lipschitz 
image of $\R^d$ (or of $\Sigma_0$, since $\Sigma_0$ is smooth).

Our condition is not necessary and sufficient, because
it is easy to build bi-Lipschitz images of $\R^d$ for which $J_1$ 
is not bounded. Even for $d=1$, if $E$ is a logarithmic spiral
centered at the origin, then $J_1(0) = +\infty$ by scale invariance
but $E$ is a bi-Lipschitz image of the line.

Let us also mention that a reasonably simple example of 
Jones and Fang [Fa] shows %%% Jones aussi?
that we cannot take $q=+\infty$ in (15.36); this
is why we prefer to take $q=1$ in the statements
above and the sufficient condition in Corollary 12.44
is a little further from optimal.

We now state the regularity result for LARRF sets
that corresponds to Theorem~14.1.

\ms\proclaim Theorem 15.39.
Let $E$ be a LARRF set, with $\varepsilon$ small enough,
depending only on $n$ and $d$. For each $\tau > 0$, there
exists $K = K(n,d,C_0,\tau) \geq 1$ such that
for $x\in E$ and $0 < r \leq 10^{-1}$, there is a $d$-plane
$P$ through $x$ and a $K$-bi-Lipschitz mapping 
$g : \R^n \to \R^n$ such that 
$$
|g(z)-z| \leq C \varepsilon
\ \hbox{ for } z \in \R^n,
\leqno (15.40)
$$
$$
g(z) = z
\ \hbox{ for } z \in \R^n \sm B(x,2r)
\leqno (15.41)
$$
and
$$
H^d(E \cap B(x,r) \sm g(P)) \leq \tau r^d.
\leqno (15.42)
$$

\ms
Indeed, let $x\in E$ and $r>0$ be given. 
By translation and dilation invariance, 
we may assume that $x=0$ and $r=20$, and that the
assumptions (15.2) and (15.3) in the definition of
LARRF are satisfied for all radii $\rho \leq 200$.

Then apply Theorem 14.1 to $E$, with
$\Sigma_0 = P(x,100)$ and $U = B(0,20)$; (12.1)
is trivial and the assumptions (12.4), (12.5), and
(14.2) follow from (15.2) and (15.3). 
The theorem gives a $K$-bi-Lipschitz mapping $g$;
(15.40) and (15.41) follow from (14.6) and (14.5),
and (14.3) says that 
$H^d(E \cap B(y,1/2) \sm g(P(x,100))) \leq \tau$
for $y\in E \cap B(0,20)$. We cover 
$E \cap B(x,r)=E\cap B(0,20)$ by less than $C$ balls 
$B(y,1/2)$ and get (15.42).
\qed

\msi{\bf Remark 15.43.}
Compared to most other results that are available
for general (locally) uniformly rectifiable sets, 
Theorem 15.39 has a significantly stronger assumption,
because Reifenberg-flatness is a quite strong regularity
property in itself. Fortunately, the conclusion is
also stronger. The main result of [DS1] %% asterisque
says that the $d$-dimensional Ahlfors-regular set $E \i \R^n$
is uniformly rectifiable if and only if it has very big pieces 
of bi-Lipschitz images of $\R^d$ into $\R^m$, where
$m = \Max(n,2d+1)$. See the condition (C5) on 
page 13 of [DS1], %% asterisque
or equivalently (1.61) in Theorem I.1.57 in [DS3]. %% UR
This means that for each $\tau > 0$,
there is a $K \geq 1$, that also depends on the Ahlfors-regularity
and uniform rectifiability constants for $E$, such that
for $x\in E$ and $r>0$, there is a $K$-bi-Lipschitz mapping
$g$ from $\R^d$ to $g(\R^d) \i \R^m$ (and where we see
$\R^n$ as embedded in $\R^m$ if $m>n$) such that
$H^d(E\cap B(x,r) \sm g(\R^d)) \leq \tau r^d$.

So the difference is not so enormous, but here we do not
have to enlarge the ambient space (we can take $m=n$), and
our mapping $g$ has a bi-Lipschitz extension to $\R^n$.

There is not so much difference in the proof either.
In [DS1], %% asterisque
the main part of the construction of $g$ happens
on a stopping time region, where we have good approximations
of $E$ by $d$-planes (as we have here automatically, by
(15.2)), but also the approximating plane stays almost
parallel to an initial one (which we do not assume here).
The main point of the proof in [DS1] %% asterisque
is a control on the number of stopping-time regions
where we need to do this construction, and then a gluing 
argument to merge the different bi-Lipschitz functions
into a big one. This is the part of the argument where
some extra room may be needed, and one uses the larger
space $\R^m$.

In the present situation, we run the same sort of algorithm,
except that we do not stop when the approximating planes
turn. This is a little more unpleasant because we have
to turn with the planes, but on the other hand
we have a unique stopping time region and nothing
to glue at the end.

In an even more general situation (where $E$ is just
a closed set), P. Jones and G. Lerman [JL] proposed a more %%
comprehensive stopping-time argument that would give
some parameterization of large pieces of $E$. See the rapid 
description near the end of Chapter 2 in [Da2]. %% UR in Utah, 
To our knowledge, this result is not published yet,
but the present paper is also close to it in spirit.
% Verify- apparently this is what I say there
\ms
The situation for Chord-Arc Surfaces with Small Constant 
is a little better, because Semmes [Se1,2]  %%CASSCS
showed that they contain very big pieces of small
Lipschitz graphs (and not merely bi-Lipschitz images of $\R^d$); 
his proof, like the proofs
for general uniformly rectifiable sets, uses corona stopping 
time regions where $E$ looks like a small Lipschitz graph.

\msi {\bf Remark 15.44: Approximating domains.}
Let us also say a few words about the locally Ahlfors-regular
sets $E$ of codimension $1$ that satisfy Condition B. 
Suppose, to make things simple, that $E$ is locally
Ahlfors-regular (as in (15.3)), bounds exactly two domains
$\Omega_1$ and $\Omega_2$, and that there is a constant
$C_1 \geq 1$ such that, for $x\in E$ and $0 < r \leq 1$,
we can find $y_1 \in \Omega_1 \cap B(x,r)$ and 
$y_2 \in \Omega_2 \cap B(x,r)$ such that
$\dist(y_j,E) \geq C_1^{-1} r$ for $j = 1, 2$.

Then there exists $\theta > 0$, that depends only
on $n$, $C_0$, and $C_1$, such that
for $x\in E$, $0 < r \leq 1$, and $j = 1, 2$,
we can find a Lipschitz domain $V_j \subset \Omega_j \cap B(x,r)$
such that $H^{n-1}(E \cap \d V_j) \geq \theta r^{n-1}$.
See [DJ], where this is used to prove estimates on the %%
harmonic measure on $\Omega_j$.

For LARRF sets of codimension $1$, the combination of
Theorem 14.1 (or 15.39) and Proposition 14.16 gives the
following result.

\ms \proclaim Proposition 15.45.
Let $E$ be a LARRF set of codimension $1$,
with $\varepsilon$ small enough (depending on $n$),
and suppose that $\R^n \sm E$ has exactly two connected
components $\Omega_1$ and $\Omega_2$. 
There exist a constant $K \geq 1$ (that depends only on
$n$ and the Ahlfors-regularity constant $C_0$
in (15.3)) and, for each $\tau > 0$, $A \geq 1$
(that depends only on $n$, $C_0$, and $\tau$) such that
the following holds. For $x\in E$ and $r < 10^{-1}$, 
we can find two disjoint $A$-Lipschitz saw-tooth domain 
$\Omega_{j,A} \i \R^n$, $j=1, 2$, and a $K$-bi-Lipschitz 
mapping $g: \R^n \to \R^n$ such that, if we set
$V_j = g(\Omega_{j,A})$ for $j =1,2$,
$$
V_j \i \Omega_j
\ \hbox{ for } j =1,2,
\leqno (15.46)
$$
$$
\d V_1 \cap \d V_2 \cap B(x,r) \i E,
\leqno (15.47)
$$
and 
$$
H^{n-1} \big(E \cap B(x,r) \sm [\d V_1 \cap \d V_2] \big) 
\leq \tau r^{n-1}.
\leqno (15.48)
$$

\ms
See Figure 3 near Proposition 14.16.
Proposition 15.45 is deduced from Theorem 14.1 and Proposition 14.16
just like we deduced Theorem 15.39 from Theorem 14.1.
\qed

\ms
%Of course $\R^n \sm E$ could have more than two components
%if $E$ is not connected. If $E$ is connected, it is easy to
%see that $\R^n \sm E$ has at most two components, and the authors
%were too lazy to check that they need to be different.
%%% Alano says: Jordan-Brouwer thm, any book with homology

Compared to the result of [DJ],  %% Cond B
we get slightly uglier domains to approximate the $\Omega_j$
with (they are only bi-Lipschitz images of Lipschitz domains),
but we get very big pieces, and we get a common piece of
boundary which is accessible from both sides.
This could perhaps be useful for problems related to elliptic
PDE.

\ms
Our bi-H\"older mapping $g$ could provide a way
to approximate $\R^n \sm E$ (when $E$ is locally Reifenberg-flat)
by more regular domains contained in $\R^n \sm E$, just
by taking images of $\R^n \sm V_{\tau}$, where $V_{\tau}$ is
a tubular  neighborhood of $\Sigma_0$. We shall not pursue this
idea here. %%% Refer to [KT] ?

\bigskip
REFERENCES

\smallskip
\item {[AFP]} L. Ambrosio, N. Fusco and D. Pallara, 
\underbar {Functions of bounded variation and free disc-}  
\underbar {ontinuity problems},
Oxford Mathematical Monographs, Clarendon Press, Oxford 2000.
\smallskip
\item {[BJ1]} C. Bishop and P. Jones, Harmonic measure and arclength. 
Ann. of Math. (2) {\bf 132} (1990), no. 3, 511--547.
\smallskip
\item {[BJ2]} C. Bishop and P. Jones, Harmonic measure, $L\sp 2$ estimates 
and the Schwarzian derivative, J. Anal. Math. {\bf 62} (1994), 77--113.
\smallskip
\item {[BJ3]} C. Bishop and P. Jones, Wiggly sets and limit sets,
Ark. Mat. {\bf 35} (1997), no. 2, 201--224.
\smallskip
\item {[Da1]} G. David, Morceaux de graphes lipschitziens et int\'egrales 
singuli\`eres sur une surface, Revista Matematica Iberoamericana, 
{\bf 4}, 1 (1988), 73--114.
\smallskip
\item {[Da2]} G. David, Uniform  rectifiability, 
Lecture notes from a course in Park City (2003), 
to be published by the AMS.
\smallskip
\item {[DJ]} G. David and D. Jerison, Lipschitz approximations to 
hypersurfaces, harmonic measure, and singular integrals, 
Indiana U. Math. Journal. {\bf 39}, 3 (1990), 831-845.
\smallskip 
\item {[DDT]} G. David, T. De Pauw, and T. Toro,
A generalization of Reifenberg's theorem in $\Bbb R^3$, 
Geom. Funct. Anal. {\bf 18} (2008), 1168--1235.
\smallskip
\item {[DS1]} G. David and S. Semmes, 
\underbar {Singular integrals and rectifiable sets in $\Bbb R^n$ : 
au-del\`a des} \quad  \underbar {graphes lipschitziens},
Ast\'erisque 193, Soci\'et\'e Math\'ematique de France 1991.
\smallskip
\item {[DS2]} G. David and S. Semmes, Quantitative rectifiability and 
Lipschitz mappings, Transactions A.M.S. {\bf 337} (1993), 855--889.
\smallskip 
\item {[DS3]} G. David and S. Semmes, 
\underbar {Analysis of and on uniformly rectifiable sets}, 
A.M.S. series of Mathematical surveys and monographs, Volume 38, 1993.
\smallskip
\item {[DS4]} G. David and S. Semmes, Uniform rectifiability 
and Singular sets, Annales de l'Inst. Henri Poincar\'e, 
Analyse non lin\'eaire, {\bf 13}, N¡ 4 (1996), p. 383--443.
\smallskip
% \item {[Do]} A. Dold, \underbar {Lectures on algebraic topology}, 
% Second edition, Grundlehren der Mathematishen Wissenschaften 200, 
% Springer Verlag 1980. 
\item {[Do]} J. R. Dorronsoro, A characterization of potential spaces, 
Proc. A.M.S. {\bf 95} (1985), 21--31.
\smallskip
\item {[Du]} J. Dugundji, \underbar {Topology}, Allyn and Bacon, 
Boston, 1966.
\smallskip
\item {[Fa]} Xiang Fang, The Cauchy integral, analytic capacity
and subsets of quasicircles,
PhD. Thesis, Yale university. 
\smallskip
\item {[Fe]} H. Federer, \underbar{Geometric measure theory}, 
Grundlehren der Mathematishen Wissenschaf-ten %% conjoncturel
153, Springer Verlag 1969.
\smallskip
\item {[Gi]} E. Giusti, \underbar{Minimal surfaces and functions of bounded 
variation}, Monographs in Mathematics, 80. Birkh\"auser Verlag, Basel-Boston, 
Mass., 1984.
\smallskip
\item {[JeK]} D. Jerison, C. Kenig, Hardy spaces, $A\sb{\infty}$, 
and singular integrals on chord-arc domains,  
Math. Scand.  {\bf 50} (1982), no. 2, 221--247.
\smallskip
\item {[J1]} P. Jones, Square functions, Cauchy integrals, 
analytic capacity, and harmonic measure, 
Proc. Conf. on Harmonic Analysis and Partial Differential Equations, 
El Escorial 1987 (ed. J. Garc\'\i a-Cuerva), p. 24-68, 
Lecture Notes in Math. 1384, Springer-Verlag 1989.
\smallskip
\item {[J2]} P. Jones, Rectifiable sets and the traveling salesman problem, 
Inventiones Mathematicae {\bf 102}, 1 (1990), 1-16.
\smallskip
\item {[JL]} P. Jones and G. Lerman, Manifold-like structures of 
measures via multiscale analysis, in preparation.
\smallskip
%%% Add Kenig-Toro ??
\item {[L\'e]}	J.-C. L\'eger,  Menger curvature and rectifiability,
Ann. of Math. (2) {\bf 149} (1999), no. 3, 831--869.
\smallskip
\item {[Lr1]}  G. Lerman, Geometric transcriptions of sets and their 
applications to data analysis, PhD thesis, Yale university.
\smallskip
\item {[Lr2]} G. Lerman, Quantifying curvelike structures of measures 
by using $L_2$ Jones quantities, to appear, Comm. Pure App. Math. {\bf 56}
(2003), no.9, 1294--1365.
\smallskip
\item {[Ma]}  P. Mattila, \underbar{Geometry of sets and 
measures in Euclidean space}, Cambridge Studies in
Advanced Mathematics 44, Cambridge University Press l995.
\item {[Mo]} C. B. Morrey, \underbar{Multiple integrals in the 
calculus of variations}, Die Grundlehren der mathematischen Wissenschaften, 
Band 130 Springer-Verlag New York, Inc., New York 1966 ix+506 pp. 
% \smallskip
% \item {[Ok]} K. Okikiolu, Characterization of subsets of rectifiable curves 
% in $\Bbb R^n$, J. of the London Math. Soc. 46 (1992), 336-348.
% \smallskip
% \item {[Pa1]} H. Pajot, Un th\'{e}or\`{e}me g\'{e}om\'{e}trique du 
% ``voyageur de commerce" en dimension 2, 
% C. R. Acad Sci. Paris, t. 323, S\'{e}rie I (1996), 13-16.
\smallskip
\item {[P1]} H. Pajot, Conditions quantitatives de rectifiabilit\'{e}, 
Bulletin de la Soci\'{e}t\'{e} Math\'{e}matique de France, 
Vol. {\bf 125} (1997), 15--53.
\smallskip
\item {[P2]} H. Pajot, \underbar{Analytic capacity, rectifiability, 
Menger curvature and the Cauchy integral}, L.N. in Math. 1799,
Springer-Verlag 2002.
\smallskip
% \item {[Ne]} M. H. A. Newman,  \underbar{Elements of the topology of plane 
% sets of points}, 
% Second edition, reprinted, Cambridge University Press, New York 1961.
% \smallskip
\item {[R1]} E. R. Reifenberg, Solution of the Plateau Problem for 
$m$-dimensional surfaces of varying topological type,
Acta Math. {\bf 104}, 1960, 1--92.
\smallskip
\item {[R2]} E. R. Reifenberg,
Epiperimetric inequality related to the analyticity of minimal surfaces, 
Annals Math., {\bf 80} (1964), 1--14.
\smallskip
\item {[R3]} E. R. Reifenberg, On the analyticity of minimal
surfaces, Annals of Math., {\bf 80} (1964), 15--21.
\smallskip
\item {[Sc]} R. Schul, Analyst's traveling salesman theorems. A 
survey, In the tradition of Ahlfors-Bers. IV,  209--220, 
Contemp. Math., {\bf 432}, Amer. Math. Soc., Providence, RI, 2007.
\smallskip
\item {[Se1]} S. Semmes, Chord-arc surfaces with small constant. I,
Adv. Math. {\bf 85} (1991), no. 2, 198--223. 
\smallskip
\item {[Se2]} S. Semmes, Chord-arc surfaces with small constant. II.
Good parameterizations,  Adv. Math. 88 (1991), no. 2, 170--199.
\smallskip
\item {[Se3]} S. Semmes, Hypersurfaces in $R\sp n$ whose unit normal 
has small BMO norm, Proc. Amer. Math. Soc. {\bf 112} (1991), no. 2, 403--412. 
% \smallskip
% \item {[Se1]} S. Semmes, A criterion for the boundedness of 
% singular integrals on hypersurfaces, Trans. Amer. Math. Soc. 311 (1989), 
% no. 2, 501--513.
% \smallskip
% \item {[Se2]} S. Semmes, Differentiable function theory on 
% hypersurfaces in  $\Bbb R^n$ (without bounds on their smoothness), 
% Indiana Univ. Math. Journal 39 (1990), 985-1004.
% \smallskip
% \item {[Se3]} S. Semmes, Analysis vs. geometry on a class of rectifiable 
% hypersurfaces in $\Bbb R^n$, Indiana Univ. Math. Journal 39 (1990), 1005-1036.
\smallskip
\item {[Si]} L. Simon, \underbar{Lectures on geometric measure theory},
Proceedings of the Centre for Mathematical Analysis, 
Australian National University, 3. Australian National University, 
Centre for Mathematical Analysis, Canberra, 1983. vii+272 pp. 
ISBN: 0-86784-429-9
\smallskip
\item {[To]} T. Toro, Geometric conditions and existence of
bi-Lipschitz
 parameterizations, Duke Math.\ Journal, {\bf 77} (1995),
193--227.
\smallskip
\item {[Tu]} P. Tukia, The planar Sch\"onflies theorem for Lipschitz maps,  
Ann. Acad. Sci. Fenn. Ser. A I Math. {\bf 5} (1980), no. 1, 49--72.
\smallskip
\item {[V]} J. V\"ais\"al\"a, Quasiconformal maps of cylindrical domains,  
Acta Math. {\bf 162} (1989), no. 3-4, 201--225.

\ms\ms\ms
\noindent Guy David,  
\smallskip\noindent 
Math\'{e}matiques, B\^atiment 425,
\smallskip\noindent 
Universit\'{e} de Paris-Sud 11, 
\smallskip\noindent 
91405 Orsay Cedex, France
\smallskip\noindent 
guy.david@math.u-psud.fr

\vglue -2.9cm 
\hskip 5.7cm
\vbox{\ms\ms
\noindent Tatiana Toro, 
\smallskip\noindent 
University of Washington 
\smallskip\noindent 
Department of Mathematics, 
\smallskip\noindent 
Box 354350, Seattle, WA 98195-4350, USA
\smallskip\noindent 
toro@math.washington.edu}

\bye